\documentclass[a4papet,10pt]{article}

\usepackage{amsfonts,amssymb,mathtools}
\usepackage{graphicx}
\usepackage{amsmath}
\usepackage[latin1]{inputenc} 
\usepackage{bbm}
\usepackage{pstricks}
\usepackage{psfrag}
\usepackage{appendix}
\usepackage{amsthm}
\usepackage{dsfont}
\usepackage{stackrel}
\usepackage{eufrak}
\usepackage{mathrsfs}
\usepackage{enumerate}
\usepackage{cite}
\usepackage{chngcntr}
\counterwithin{figure}{section}
\usepackage[a4paper,top=30mm,bottom=30mm,left=35mm,right=35mm]{geometry}
\numberwithin{equation}{section}

\newcommand{\prob}{\mathbb{P}}
\newcommand{\Ex}{\mathbb{E}}
\newcommand{\Rl}{\mathbb{R}}

\newcommand{\dom}{\mbox{dom}\,}
\newcommand{\doms}{\scriptsize{\dom}}

\newcommand{\rew}{\Rl^d}

\newtheorem{theorem}{Theorem}[section]
\newtheorem{proposition}{Proposition}[section]
\newtheorem{lemma}{Lemma}[section]
\newtheorem{corollary}{Corollary}[section]
\newtheorem{example}{Example}[section]

\title{Statistical fluctuations under resetting: rigorous results}

\author{Marco Zamparo\footnote{Dipartimento di Fisica, Universit\`a degli Studi di Bari and INFN, Sezione di Bari,
    via Amendola 173, \phantom{aaz} 70126 Bari, Italy
      \newline \phantom{aaz} E-mail: \texttt{marco.zamparo@uniba.it}}}

\date{}

\begin{document}  
\maketitle

\begin{abstract}
In this paper we investigate the normal and the large fluctuations of
additive functionals associated with a stochastic process under a
general non-Poissonian resetting mechanism. Cumulative functionals of
regenerative processes are very close to renewal-reward processes and
inherit most of the properties of the latter. Here we review and use
the classical law of large numbers and central limit theorem for
renewal-reward processes to obtain same theorems for additive
functionals of a stochastic process under resetting. Then, we
establish large deviation principles for these functionals by
illustrating and applying a large deviation theory for renewal-reward
processes that has been recently developed by the author. We discuss
applications of the general results to the positive occupation time,
the area, and the absolute area of the reset Brownian motion.  While
introducing advanced tools from renewal theory, we demonstrate that a
rich phenomenology accounting for dynamical phase transitions emerges
when one goes beyond Poissonian resetting.\\


\end{abstract}

\section{Introduction}

Recent years have witnessed a growing interest of physicists in
stochastic processes under resetting, whereby a process starts anew
from its initial position at random times. Resetting mechanisms have
been proven to represent a simple way of generating non-equilibrium
stationary states \cite{Manrubia,Evans_rew} and anomalous diffusion
while keeping stationarity \cite{Padova}. They have been also verified
to speed up some searching tasks by defining intermittent search
strategies, which interspace periods of random exploration with random
returns to a starting point \cite{Evans_rew}.

Physicists have focused much of their work on the Brownian motion and
related processes under Poissonian resetting, studying deviations from
the typical behavior and addressing first-passage time problems in
connection with intermittent search strategies \cite{Evans_rew}.  In
particular, the ability of stochastic resetting to shape fluctuations
has been tested for additive functionals, such as the positive
occupation time of the Brownian motion \cite{LDP_res_2}, its area and
absolute area \cite{LDP_res_2}, the area of the fractional Brownian
motion \cite{SmithMaj}, the area of the Ornstein-Uhlenbeck process
\cite{LDP_res_1}, and the local time of the Brownian motion
\cite{LDP_res_3,LDP_res_4}. For the Brownian motion subjected to
Poissonian resetting, also positive and negative excursions have been
investigated \cite{LDP_res_exc}.  As resetting has the effect of
confining the process around the initial position, leading to the
emergence of a stationary state for the reset Brownian motion
\cite{Evans_rew}, one natural question is whether or not the resetting
mechanism can bring about a large deviation principle for an additive
functional when it does not satisfy a large deviation principle
without resetting. Poissonian resetting has been discovered to be
frequent enough to bring about a large deviation principle in some
cases, such as the positive occupation time of the Brownian motion
\cite{LDP_res_2}, but not in other cases, such as its area
\cite{LDP_res_2}.  These findings have finally prompted physicists to
inspect moderate deviations for the area of the Brownian motion and
the fractional Brownian motion under Poissonian resetting
\cite{SmithMaj}.

With this paper we drop the Poissonian hypothesis and answer the above
question when the resetting protocol is allowed to be any, and in
particular when it is allowed to be adapted to the additive functional
under consideration.  To date, there are no contributions in this
direction as far as we know.  Most likely, one of the reason that has
made Poissonian resetting very attractive is the fact that it is the
only resetting mechanism for which the moment generating function of
additive functionals can be obtained explicitly at any time
\cite{LDP_res_1}, so that their asymptotic properties can be
investigated by means of elementary mathematical tools.  Outside the
Poissonian framework, advanced mathematical arguments are needed to
get around the impossibility to work out every single step.  We resort
to the mathematical theory of renewal-reward processes to go beyond
Poissonian resetting.

From the mathematical point of view, the mechanism of resetting is
nothing but an aspect of renewal theory, which has been developed
since the 40s by Doob, Feller, Cox, and many others \cite{Smith1}. In
fact, the resetting process is a \textit{renewal process}, which is
meant to describe some event that is renewed randomly over time. In
renewal theory the occurrences of an event often involve some
broad-sense reward, and in this case one speaks of
\textit{renewal-reward systems}. For example, renewals can represent
the random arrivals of customers at a store and rewards can be the
random amount of money that each customer spends in the store.
Renewal-reward systems find classical application in queueing theory
\cite{AsmussenBook}, insurance \cite{DicksonBook}, finance
\cite{RSSTBook}, and statistical physics of polymers
\cite{GiacominBook,Frank_polymer} among others. A renewal-reward
system is characterized by the {\it renewal times} $T_1,T_2,\ldots$ at
which the event occurs and by the corresponding {\it rewards}
$X_1,X_2,\ldots$. If $S_1,S_2,\ldots$ denote the {\it waiting times}
for a new occurrence of the event, then the renewal time $T_n$ can be
expressed for each $n\ge 1$ in terms of the waiting times as
$T_n=S_1+\cdots+S_n$. Standard assumptions, which we also make, are
that the waiting time and reward pairs $(S_1,X_1),(S_2,X_2),\ldots$
form an independent and identically distributed (i.i.d.) sequence of
random vectors on a probability space $(\Omega,\mathcal{F},\prob)$,
the waiting times taking non-negative real values and the rewards
taking values in $\rew$ with some dimension $d\ge 1$.  Any dependence
between $S_n$ and $X_n$ is allowed.  To avoid trivialities we suppose
that the waiting times are not concentrated at zero, that is
$\prob[S_1=0]<1$. This implies in particular that almost all the
trajectories account for only finitely many renewals in finite
time\footnote{Since $\prob[S_1>0]>0$, there exists $\delta>0$ with the
  property that $\prob[S_1>\delta]>0$.  The second Borel-Cantelli
  lemma states that $S_n>\delta$ for infinitely many $n$ almost
  surely, so that $\lim_{n\uparrow\infty}T_n=+\infty$ $\prob$-a.s..}.
The number of renewals up to the time $t\ge 0$ is $N_t:=\sup\{n\ge
0:T_n\le t\}$ with $T_0:=0$.

A quantity of great interest in renewal-reward systems is the {\it
  cumulative reward} by the time $t\ge 0$, which is the random
variable
\begin{equation*}
  W_t:=\sum_{n=1}^{N_t}X_n
\end{equation*}
with $W_t:=0$ if $N_t=0$. The stochastic process $\{W_t\}_{t\ge 0}$ is
the so-called {\it renewal-reward process}, or compound renewal
process. Insisting on the above example, $W_t$ would be the amount of
money earned by the store up to the time $t$.  The typical value $\mu$
of the cumulative reward in the large $t$ limit is identified by the
following \textit{strong law of large numbers} (SLLN)
\cite{Gut_book}. Hereafter, $\Ex$ denotes expectation with respect to
the law $\prob$ and $\|\cdot\|$ is the Euclidean norm in $\Rl^d$:
$\|x\|:=\sqrt{x\cdot x}$ for every $x\in\Rl^d$.
\begin{theorem}
  \label{SLLN}
  If $\Ex[S_1]<+\infty$ and $\Ex[\|X_1\|]<+\infty$, then
\begin{equation*}
  \lim_{t\uparrow+\infty}\frac{W_t}{t}=\frac{\Ex[X_1]}{\Ex[S_1]}=:\mu~~~~~\prob\mbox{-a.s..}
\end{equation*}
\end{theorem}
The normal fluctuations around $\mu$ are described by the following
\textit{central limit theorem} (CLT) \cite{Gut_book}, which introduces the
asymptotic covariance matrix $\Sigma$ of $W_t$.
\begin{theorem}
  \label{CLT}
If $\Ex[S_1^2]<+\infty$ and $\Ex[\|X_1\|^2]<+\infty$, then for each
Borel set $A\subseteq\Rl^d$
\begin{equation*}
\lim_{t\uparrow+\infty}\prob\bigg[\frac{W_t-\mu t}{\sqrt{t}}\in A\bigg]=\int_A \frac{e^{-\frac{1}{2} x\,\cdot \Sigma^{-1} x}}{\sqrt{(2\pi)^d\det\Sigma}}\, dx
\end{equation*}
provided that the matrix $\Sigma$ is positive definite, with
\begin{equation*}
\Sigma:=\frac{\mbox{\upshape cov}[X_1-\mu S_1]}{\Ex[S_1]}.
\end{equation*}
\end{theorem}

Theorems \ref{SLLN} and \ref{CLT} are well-established properties of
renewal-reward processes that descend from basic tools of probability
theory, the former being a simple manipulation of the Kolmogorov's law
of large numbers and the latter being a combination of the classical
central limit theorem with the Anscombe's theorem for random sequences
with random index.  The typical fluctuations can be also described
when the variances of $S_1$ and $\|X_1\|$ are not finite by
introducing proper scaling and possibly L\'evy stable laws
\cite{Gut_book}.

While the study of the typical behaviors of the cumulative reward
$W_t$ can be considered concluded, the study of its large fluctuations
is still in progress and some definitive results have been achieved
only very recently.  In fact, many mathematicians have investigated
the large fluctuations of $W_t$ with increasing level of generality
\cite{Glynn1994,Jiang1994,Duffield1998,Macci2005,Macci2007,Lefevere2011,Borovkov2015,Borovkov2019}. At
last, the author has provided large deviation principles with optimal
hypotheses for lattice waiting-time distributions first
\cite{Models,Marco_discrete} and any waiting-time distribution later
\cite{Marco_continuous}. We shall review the findings of
\cite{Marco_continuous} in the next section.

In this paper we make use of the theory of renewal-reward processes to
investigate normal and large fluctuations of additive functionals
associated with a stochastic process under a resetting mechanism. The
reference stochastic process and the resetting mechanism are allowed
to be any. While providing general results that characterize both
typical and rare behaviors, we show that resetting at a higher
frequency than Poissonian resetting brings about a large deviation
principle where the Poissonian protocol fails.

We apply our general results to the positive occupation time, the
area, and the absolute area of the reset Brownian motion, continuing
the study initiated in \cite{LDP_res_2}. A large deviation theory
beyond Poissonian resetting allows to appreciate the emergence of
singularities in the graph of their rate functions, which are
interpreted as dynamical phase transitions. Dynamical phase
transitions have been already documented for time-averaged quantities
such as the heat exchanged after a quench \cite{DPT1,DPT2},the heat
exchanged between two thermal walls
\cite{Lefevere_heat1,Lefevere_heat2}, the active work in active matter
\cite{DPT3}, the total displacement of random walks
\cite{DPT4,LDP_res_5,DPT12,DPT555}, occupation times of a Brownian
particle \cite{DPT5,DPT6}, and the entropy production
\cite{DPT7,DPT8,DPT9,DPT10,DPT11}.  Physicists often invoke the
G\"artner-Ellis theorem to justify their calculations that lead to
dynamical phase transitions, although it is necessary for a dynamical
phase transition to occur that the hypotheses of that theorem are
violated. We do not suffer from this contradiction because we are able
to develop a theory for the large fluctuations of additive functionals
of stochastic processes under resetting that does not require the
smoothness hypotheses of the G\"artner-Ellis theorem.

The paper is organized as follows. In Section \ref{sec:LDP_W} we
review the large deviation principles for renewal-reward processes
developed by the author in \cite{Marco_continuous}. In this section we
also propose a formula for the rate function that facilitates
applications. In Section \ref{sec:rBM} we introduce general stochastic
processes under a resetting mechanism and the associated additive
functionals. For general additive functionals we prove the law of
large numbers, the central limit theorem, and a large deviation
principle. In Sections \ref{sec:applications} we apply these results
to the positive occupation time, the area, and the absolute area of
the reset Brownian motion, supplying a complete characterization of
their rate functions. In this section we also discuss how rate
functions behave in the limit of zero resetting through some examples.
Conclusions and prospects for future research are finally reported in
Section \ref{sec:conclusions}.  This paper is essentially a math work.
The most technical proofs are postponed in the appendices and are not
necessary to understand the main line of discussion.

\section{Large deviations for renewal-reward processes}
\label{sec:LDP_W}

The cumulative reward $W_t$ is said to satisfy a \textit{weak large
  deviation principle} (weak LDP) with the \textit{rate function} $I$
if there exists a lower semicontinuous function $I$ on $\Rl^d$ such
that
\begin{enumerate}[$(i)$]
\item $\liminf_{t\uparrow+\infty}\frac{1}{t}\ln\prob\big[\frac{W_t}{t}\in G\big]\ge -\inf_{w\in G}\{I(w)\}$ for each open set $G\subseteq\Rl^d$,
\item $\limsup_{t\uparrow+\infty}\frac{1}{t}\ln\prob\big[\frac{W_t}{t}\in K\big]\le -\inf_{w\in K}\{I(w)\}$ for each compact set $K\subset\Rl^d$.
\end{enumerate}
Bounds $(i)$ and $(ii)$ are called \textit{lower large deviation
  bound} and \textit{upper large deviation bound}, respectively. The
cumulative reward $W_t$ is said to satisfy a \textit{full large
  deviation principle} (full LDP) if the upper large deviation bound
holds for all closed sets $K$, and not only for compact sets. The rate
function $I$ is denoted \textit{good} if it has compact level sets,
namely if the set $\{w\in\Rl^d:I(w)\le a\}$ is compact for each $a\ge
0$. We refer to \cite{DemboBook,Frank_LDP} for the language of large
deviation theory.  In this section we review the LDPs for
renewal-reward processes established by the author in
\cite{Marco_continuous}.  For simplicity, and for the applications we
have in mind, here we assume that the limit
\begin{equation*}
 \lim_{s\uparrow+\infty}-\frac{1}{s}\ln\prob[S_1>s]=:\ell
\end{equation*}
exists as an extended real number. Notice that
$\ell\in[0,+\infty]$. This assumption is not made in
\cite{Marco_continuous}, where, besides, possible infinite-dimensional
rewards are considered.

The results of \cite{Marco_continuous} rely on Cram\'er's theory for
waiting times and rewards.  The joint Cram\'er's rate function of
waiting times and rewards is the function $J$ that maps each
$(s,w)\in\Rl\times\Rl^d$ in the extended real number
\begin{equation*}
J(s,w):=\sup_{(\zeta,k)\in\Rl\times\Rl^d}\Big\{s\zeta+w\cdot k-\ln\Ex\big[e^{\zeta S_1+k\,\cdot X_1}\big]\Big\}.
\end{equation*}
The lower-semicontinuous regularization $\Upsilon$ of
$\inf_{\gamma>0}\{\gamma J(\cdot/\gamma,\cdot/\gamma)\}$ is the
function that associates every $(\beta,w)\in\Rl\times\Rl^d$ with
\begin{equation*}
  \Upsilon(\beta,w):=\adjustlimits\lim_{\delta\downarrow 0}\inf_{s\in (\beta-\delta,\beta+\delta)}\inf_{v\in \Delta_{w,\delta}}\inf_{\gamma>0}
  \bigg\{\gamma J\bigg(\frac{s}{\gamma},\frac{v}{\gamma}\bigg)\bigg\},
\end{equation*}
$\Delta_{w,\delta}:=\{v\in\rew:\|v-w\|<\delta\}$ being the open ball
of center $w$ and radius $\delta$. We make use of $\Upsilon$ to
construct a ``rate function'' $I$ according to the formula
\begin{equation*}
  I(w):=\begin{cases}
  \inf_{\beta\in[0,1]}\{\Upsilon(\beta,w)+(1-\beta)\ell\} & \mbox{if }\ell<+\infty,\\
  \Upsilon(1,w) & \mbox{if }\ell=+\infty
  \end{cases}
\end{equation*}
for all $w\in\Rl^d$.  The following theorem states a weak and a full
LDP for the cumulative reward $W_t$ and is (part of) Theorem 1.1 of
\cite{Marco_continuous}. The function $I$ turns out to really be the
rate function.
\begin{theorem}
\label{mainth}
The following conclusions hold:
\begin{enumerate}[$(i)$]
\item the function $I$ is lower semicontinuous and convex;
\item $W_t$ satisfies a weak LDP with the rate function
  $I$. Furthermore, the upper large deviation bound holds for each
  Borel convex set $K$ whenever $\ell<+\infty$ or $I(0)<+\infty$;
\item if $\Ex[e^{\rho\|X_1\|}]<+\infty$ for some number $\rho>0$, then
  $I$ has compact level sets and $W_t$ satisfies a full LDP with the
  good rate function $I$.
\end{enumerate}
\end{theorem}

We stress that no assumption on the law of waiting times and rewards
is needed for the validity of a weak LDP.  Moreover, the upper large
deviation bound extends to convex sets if $\ell<+\infty$ or
$I(0)<+\infty$, whereas it fails in general when $\ell=+\infty$ and
$I(0)=+\infty$ as shown in \cite{Marco_continuous} by means of some
counterexamples. Regarding a full LDP with a good rate function, we
can observe that it does not require an hypothesis on waiting times
but only an exponential moment condition on rewards:
$\Ex[e^{\rho\|X_1\|}]<+\infty$ for some $\rho>0$.

We conclude the section by providing a formula for the rate function
$I$ that facilitates the use of Theorem \ref{mainth} in
applications. Let $\varphi$ be the function that maps any $k\in\Rl^d$
in the extended real number
\begin{equation*}
\varphi(k):=\sup\Big\{\zeta\in\Rl~:~\Ex\big[e^{\zeta S_1+k\,\cdot X_1}\big]\le 1\Big\},
\end{equation*}
where the supremum over the empty set is customarily interpreted as
$-\infty$. The following proposition shows that $I$ is the convex
conjugate, i.e.\ the Legendre-Fenchel transform, of $-(\varphi\wedge
\ell)$. As usual, given two extended real numbers $a$ and $b$, we
denote $\min\{a,b\}$ by $a\wedge b$ and $\max\{a,b\}$ by $a\vee b$ for
brevity. The proof of the proposition is presented in Appendix
\ref{proof:conv}.
\begin{proposition}
  \label{prop:conv}
The following conclusions hold:
\begin{enumerate}[$(i)$]
\item the function $\varphi$ is upper semicontinuous and concave;
\item $I(w)=\sup_{k\in\Rl^d}\{w\cdot k+\varphi(k)\wedge \ell\}$ for all $w\in\Rl^d$.
\end{enumerate}
\end{proposition}

\section{Stochastic processes under resetting and fluctuations}
\label{sec:rBM}

A stochastic process under resetting can be constructed from
independent copies of the reference process as follows. Let
$\{B_{1,t}\}_{t\ge 0}, \{B_{2,t}\}_{t\ge 0},\ldots$ be
i.i.d.\ $\Rl$-valued stochastic processes on a probability space
$(\Omega,\mathcal{F},\prob)$ where a renewal process $\{T_n\}_{n\ge
  0}$ is already given.  We think of these processes as fresh
realizations of the same stochastic dynamics, say $\{B_{1,t}\}_{t\ge
  0}$, and we defined the \textit{dynamics under resetting}
$\{Z_t\}_{t\ge 0}$ by pasting them together at each renewal time:
\begin{equation*}
Z_t:=B_{N_t+1,t-T_{N_t}}
\end{equation*}
for all $t\ge 0$. We remind that $N_t:=\sup\{n\ge 0: T_n\le t\}$ is
the number of renewals by the time $t$, so that $Z_t=B_{n,t-T_{n-1}}$
for $t\in[T_{n-1},T_n)$ and $n\ge 1$. The stochastic process under
  resetting $\{Z_t\}_{t\ge 0}$ is the \textit{reset Brownian motion}
  or the \textit{reset fractional Brownian motion} when
  $\{B_{1,t}\}_{t\ge 0}, \{B_{2,t}\}_{t\ge 0},\ldots$ are Brownian
  motions or fractional Brownian motions, respectively. It is the
  \textit{reset Ornstein-Uhlenbeck process} if $\{B_{1,t}\}_{t\ge 0},
  \{B_{2,t}\}_{t\ge 0},\ldots$ are Ornstein-Uhlenbeck processes.
  According to the existing literature
  \cite{Evans_rew,LDP_res_1,LDP_res_2,LDP_res_3,LDP_res_4,LDP_res_exc,SmithMaj},
  we assume that the processes $\{B_{1,t}\}_{t\ge 0},
  \{B_{2,t}\}_{t\ge 0},\ldots$ are independent of the waiting
  times. We also assume that their sample paths are measurable in
  order to enable integration over time.

An \textit{additive functional} of the process $\{Z_t\}_{t\ge 0}$ is
the random variable
  \begin{equation*}
   F_t:=\int_0^t f(Z_\tau)\,d\tau
  \end{equation*}
where $f$ is some real measurable function over $\Rl$. For instance,
the positive occupation time and the absolute area studied in
\cite{LDP_res_2} for the reset Brownian motion with Poissonian
resetting are the additive functionals corresponding to
$f(z):=\mathds{1}_{\{z>0\}}$ and $f(z):=|z|$, respectively. The area
investigated in \cite{LDP_res_2} for the reset Brownian motion with
Poissonian resetting, in \cite{SmithMaj} for the reset fractional
Brownian motion with Poissonian resetting, and in \cite{LDP_res_1} for
the reset Ornstein-Uhlenbeck process with Poissonian resetting, is
found with $f(z):=z$.  We recall that \textit{Poissonian resetting}
amounts to a pure exponential waiting time distribution:
$\prob[S_1>s]=e^{-rs}$ for all $s\ge 0$ with some resetting rate
$r>0$.

In this section we provide general results for the fluctuations of the
additive functional $F_t$ by resorting to the theory of renewal-reward
processes. The random variable $F_t$ is closely related to a
cumulative reward. In fact, if we associate the ``reward''
  \begin{equation*}
    X_n:=\int_0^{S_n}f(B_{n,\tau})\,d\tau
  \end{equation*}
with the $n$th renewal at the time $T_n$, then $F_t$ differs from the
cumulative reward $W_t:=\sum_{n=1}^{N_t}X_n$ by the contribution
$\int_0^{t-T_{N_t}}f(B_{N_t+1,\tau})\,d\tau$ of the backward
recurrence time.  The \textit{backward recurrence time}, or current
lifetime, is the time $t-T_{N_t}$ elapsed from the last renewal.
Explicitly, for all $t\ge 0$ we can write
\begin{align}
  \nonumber
  F_t&=\sum_{n=1}^{N_t}\int_{T_{n-1}}^{T_n}f(B_{n,\tau-T_{n-1}})\,d\tau+\int_{T_{N_t}}^t f(B_{N_t+1,\tau-T_{N_t}})\,d\tau\\
  &=\sum_{n=1}^{N_t}X_n+\int_0^{t-T_{N_t}}f(B_{N_t+1,\tau})\,d\tau=W_t+\int_0^{t-T_{N_t}}f(B_{N_t+1,\tau})\,d\tau.
  \label{F_WRT}
\end{align}
Such relationship suggests that some limit theorems for $F_t$ can be
obtained from limit theorems for $W_t$. This is quite immediate for
the SLLN and the CLT, where the ``incomplete reward''
$\int_0^{t-T_{N_t}}f(B_{N_t+1,\tau})\,d\tau$ disappears upon rescaling
under mild assumptions on the resetting mechanism. In other words, the
backward recurrence time cannot contribute to the typical behavior of
$F_t$. The issue of large deviations is more involved since the
backward recurrence time comes into play.  The typical behavior of the
additive functional $F_t$ is investigated in Section
\ref{sec:F_normal}, whereas the large fluctuations are addressed in
Section \ref{sec:F_large}.

\subsection{Typical behaviors of additive functionals}
\label{sec:F_normal}

Set for brevity
\begin{equation*}
  M_n:=\int_0^{S_n}|f(B_{n,\tau})|\,d\tau\ge|X_n|.
\end{equation*}
Under the hypotheses $\Ex[S_1]<+\infty$ and $\Ex[M_1^\rho]<+\infty$
with some $\rho>0$, the difference between the additive functional
$F_t$ and the cumulative reward $W_t:=\sum_{n=1}^{N_t}X_n$ satisfies
\begin{equation}
    \lim_{t\uparrow+\infty}t^{-\frac{1}{\rho}}|F_t-W_t|=0~~~~~\prob\mbox{-a.s.}.
\label{aux}
\end{equation}
In fact, for every $\epsilon>0$ we have
\begin{equation*}
  \sum_{n=1}^\infty\prob\big[n^{-\frac{1}{\rho}}M_n>\epsilon\big]=\sum_{n=1}^\infty\prob\big[M_1^\rho/\epsilon^\rho>n\big]\le
  \int_0^{+\infty}\prob\big[M_1^\rho/\epsilon^\rho>x\big]dx
  =\Ex\big[M_1^\rho/\epsilon^\rho\big]<+\infty.
\end{equation*}
Thus, the Borel-Cantelli lemma yields
$\lim_{n\uparrow\infty}n^{-\frac{1}{\rho}}M_n=0$ $\prob$-a.s.. At the
same time, Theorem \ref{SLLN} with unit rewards states that
$\lim_{t\uparrow+\infty}N_t/t=1/\Ex[S_1]>0$ $\prob$-a.s.. By combining
these two limits we realize that
$\lim_{t\uparrow\infty}t^{-\frac{1}{\rho}}M_{N_t+1}=0$
$\prob$-a.s.. At this point, (\ref{aux}) is evident since
(\ref{F_WRT}) and the inequality
$t-T_{N_t}<T_{N_t+1}-T_{N_t}=S_{N_t+1}$ give for all $t>0$
\begin{equation*}
\big|F_t-W_t\big|=\bigg|\int_0^{t-T_{N_t}}f(B_{N_t+1,\tau})\,d\tau\bigg|\le\int_0^{S_{N_t+1}}\big|f(B_{N_t+1,\tau})\big|\,d\tau=M_{N_t+1}.
\end{equation*}
Limit (\ref{aux}) has consequences that are easy to see.

Theorem \ref{SLLN} tells us that $\lim_{t\uparrow+\infty}W_t/t=\mu$
$\prob$-a.s.\ if $\Ex[S_1]<+\infty$ and $\Ex[|X_1|]<+\infty$. Under
the slightly more restrictive condition $\Ex[S_1]<+\infty$ and
$\Ex[M_1]<+\infty$ we have $\lim_{t\uparrow+\infty}t^{-1}|F_t-W_t|=0$
by (\ref{aux}).  Thus, the following SLLN for $F_t$ holds true.
\begin{theorem}
  \label{th:SLLN_F}
If $\Ex[S_1]<+\infty$ and $\Ex[M_1]<+\infty$, then
\begin{equation*}
\lim_{t\uparrow+\infty}\frac{F_t}{t}=\frac{\Ex[X_1]}{\Ex[S_1]}=:\mu~~~~~\prob\mbox{-a.s.}.
\end{equation*}
\end{theorem}  

If $\Ex[S_1^2]<+\infty$ and $\Ex[M_1^2]<+\infty$, then
$\Ex[X_1^2]<+\infty$ and, provided that $\Ex[(X_1-\mu S_1)^2]\ne 0$,
Theorem \ref{CLT} ensures us that $(W_t-\mu t)/\sqrt{vt}$ converges in
distribution as $t$ is sent to infinity to a Gaussian variable with
mean $0$ and variance $1$. The number $v$ is the asymptotic variance
of $W_t$:
\begin{equation}
v:=\frac{\Ex[(X_1-\mu S_1)^2]}{\Ex[S_1]}.
\label{variance_W}
\end{equation}
At the same time, (\ref{aux}) shows that
$\lim_{t\uparrow+\infty}t^{-\frac{1}{2}}|F_t-W_t|=0$ $\prob$-a.s..
Thus, Slutsky's theorem tells us that $(F_t-\mu t)/\sqrt{vt}$
converges in distribution to the same limit as $(W_t-\mu
t)/\sqrt{vt}$.  These arguments prove the following CLT for $F_t$.
\begin{theorem}
  \label{th:CLT_F}
If $\Ex[S_1^2]<+\infty$ and $\Ex[M_1^2]<+\infty$, then for every
$z\in\Rl$
\begin{equation*}
  \lim_{t\uparrow+\infty}\prob\bigg[\frac{F_t-\mu t}{\sqrt{vt}}\le z\bigg]=\int_{-\infty}^z \frac{e^{-\frac{1}{2}x^2}}{\sqrt{2\pi}}\,dx
\end{equation*}
provided that $v$ given by (\ref{variance_W}) is not zero.
\end{theorem}

When the function $f$ that defines the additive functional is positive
or negative, such as for the positive occupation time and the absolute
area, the condition $\Ex[X_1^2]<+\infty$ that makes the asymptotic
variance $v$ finite also is the condition $\Ex[M_1^2]<+\infty$ that
neutralizes the contribution of the backward recurrence time. In
Section \ref{sec:applications} we shall see that the conditions
$\Ex[X_1^2]<+\infty$ and $\Ex[M_1^2]<+\infty$ are tantamount even in
cases where $f$ changes sign, such as the case of the area of the
reset Brownian motion.

\subsection{Large fluctuations of additive functionals}
\label{sec:F_large}

We discuss the large fluctuations of $F_t$ by applying same
definitions of Section \ref{sec:LDP_W}. The additive functional $F_t$
satisfies a \textit{full large deviation principle} (full LDP) with
the \textit{rate function} $I$ if there exists a lower semicontinuous
function $I$ on $\Rl$ such that
\begin{enumerate}[$(i)$]
\item $\liminf_{t\uparrow+\infty}\frac{1}{t}\ln\prob\big[\frac{F_t}{t}\in G\big]\ge -\inf_{w\in G}\{I(w)\}$ for each open set $G\subseteq\Rl$,
\item $\limsup_{t\uparrow+\infty}\frac{1}{t}\ln\prob\big[\frac{F_t}{t}\in K\big]\le -\inf_{w\in K}\{I(w)\}$ for each closed set $K\subseteq\Rl$.
\end{enumerate}
The rate function $I$ is \textit{good} if it has compact level sets.
A remark is in order about the question raised in the introduction
regarding the ability of the resetting mechanism to bring about a LDP
for $F_t$ when it does not satisfy a LDP without resetting. In most
cases both $F_t$ and the additive functional without resetting, namely
$\int_0^tf(B_{1,\tau})\,d\tau$, are expected to satisfy a LDP
according to the above definition. Then, the question actually is
whether the rate function of $F_t$ turns out to be positive, so that
the large fluctuations of $F_t$ really occur with probability that is
exponentially small in the time $t$, whereby the rate function of
$\int_0^tf(B_{1,\tau})\,d\tau$ is zero. The property of goodness of
$I$ excludes that the region of its zeros extends to infinity.

As in Section \ref{sec:LDP_W}, we require that the following limit
exists:
\begin{equation*}
\ell:=\lim_{s\uparrow+\infty}-\frac{1}{s}\ln\prob[S_1>s].
\end{equation*}
Then, by Theorem \ref{mainth} and Proposition \ref{prop:conv} the
cumulative reward $W_t:=\sum_{n=1}^{N_t}X_n$ satisfies a weak LDP and
its rate function is the convex conjugate of $-(\varphi\wedge \ell)$,
$\varphi$ being the function that maps $k\in\Rl$ in
\begin{equation*}
\varphi(k):=\sup\Big\{\zeta\in\Rl~:~\Ex\big[e^{\zeta S_1+kX_1}\big]\le 1\Big\}.
\end{equation*}
When adding the contribution
$\int_0^{t-T_{N_t}}f(B_{N_t+1,\tau})\,d\tau$ of the backward
recurrence time to $W_t$ in order to form $F_t$, a putative rate
function is the lower semicontinuous convex function $I$ that
associates $w\in\Rl$ with
\begin{equation}
  I(w):=\sup_{k\in\Rl}\big\{wk+\varphi(k)\wedge \varpi(k)\big\},
  \label{def:I_F_0}
\end{equation}
where
\begin{equation}
  \varpi(k):=\liminf_{t\uparrow+\infty}-\frac{1}{t}\ln\big\{\mathcal{E}_t(k)\prob[S_1>t]\big\}
\label{def:working_hypothesis}
\end{equation}
and $\mathcal{E}_t(k):=\Ex[e^{k\int_0^tf(B_{1,\tau})\,d\tau}]$ is the
value at $k$ of the moment generating function of the additive
functional without resetting. If we drop $\mathcal{E}_t(k)$, then
$\varpi(k)=\ell$ and formula (\ref{def:I_F_0}) returns the original
rate function of $W_t$. Thus, we can get a cue for introducing
$\mathcal{E}_t(k)$ by observing that the occurrence of the event
$S_1>t$, i.e.\ $N_t=0$, entails $F_t=\int_0^tf(B_{1,\tau})\,d\tau$
whereas $W_t=0$.  The following proposition demonstrates an estimate
of the asymptotic logarithmic moment generating function of $F_t$ and
the upper large deviation bound with the rate function $I$. The proof
is presented in Appendix \ref{proof:LDP_preliminar}.  We recall that
the \textit{asymptotic logarithmic moment generating function} of
$F_t$ is the function $g$ that maps $k\in\Rl$ in
  \begin{equation*}
  g(k):=\lim_{t\uparrow+\infty}\frac{1}{t}\ln\Ex\big[e^{kF_t}\big].
  \end{equation*}
The asymptotic logarithmic moment generating function $g$ does not
exist if this limit does not exist at some point $k$. Investigation of
$g$ will allow to compare our theory with the G\"artner-Ellis theory.
\begin{proposition}
  \label{th:LDP_preliminar}
  The following conclusions hold with $I$ defined by
  (\ref{def:I_F_0}):
  \begin{enumerate}[$(i)$]
  \item  $\limsup_{t\uparrow+\infty}\frac{1}{t}\ln\Ex[e^{kF_t}]\le-\varphi(k)\wedge \varpi(k)$ for all $k\in\Rl$;
  \item $\limsup_{t\uparrow+\infty}\frac{1}{t}\ln\prob\big[\frac{F_t}{t}\in K\big]\le -\inf_{w\in K}\{I(w)\}$ for each closed set $K\subseteq\Rl$.
  \end{enumerate}
\end{proposition}

There is no possibility to promote the upper large deviation bound
stated by part $(ii)$ of Proposition \ref{th:LDP_preliminar} to a LDP
without some hypothesis on the contribution of the backward recurrence
time. In all cases concerning the reset Brownian motion that we have
considered we have found $\varpi(k)\ge\varphi(k)$ for all $k\in\Rl$,
meaning that the fluctuations of the incomplete reward
$\int_0^{t-T_{N_t}}f(B_{N_t+1,\tau})\,d\tau$ do not overcome the
fluctuations of complete rewards.  Actually, we are not able to
exhibit an example that violates $\varpi(k)\ge\varphi(k)$ for some
$k$.  Thus, we are going to assume that the condition
$\varpi(k)\ge\varphi(k)$ is fulfilled for all $k\in\Rl$. This
condition gives
\begin{equation}
  I(w):=\sup_{k\in\Rl}\big\{wk+\varphi(k)\big\}.
 \label{def:I_F}
\end{equation}
Importantly, this condition suffices to obtain a LDP for the additive
functional $F_t$, as stated by the following theorem. The proof of the
lower large deviation bound is based on Theorem \ref{mainth} and
Proposition \ref{prop:conv} and is reported in Appendix
\ref{proof:LDP}.
\begin{theorem}
  \label{th:LDP}
  The following conclusions hold with $I$ defined by (\ref{def:I_F}):
  \begin{enumerate}[$(i)$]
  \item if there exists $\rho>0$ such that
    $\Ex[e^{\rho|X_1|}]<+\infty$, then $I$ is a good rate function. If
    moreover $\ell>0$, then $I(w)=0$ if and only if
    $w=\mu:=\Ex[X_1]/\Ex[S_1]$;
  \item if $\varpi(k)\ge \varphi(k)$ for all $k\in\Rl$, then $F_t$
    satisfies a full LDP with the rate function $I$;
    \item if $\varpi(k)\ge \varphi(k)$ for all $k\in\Rl$, then the
      asymptotic logarithmic moment generating function $g$ of $F_t$
      exists and equals $-\varphi$.
  \end{enumerate}
\end{theorem}

We stress that the condition $\varpi(k)\ge\varphi(k)$ for all
$k\in\Rl$ does not mean that the contribution of the backward
recurrence time is negligible. It means that the fluctuations of the
incomplete reward $\int_0^{t-T_{N_t}}f(B_{N_t+1,\tau})\,d\tau$ are not
dominant, but they could be comparable to those of complete
rewards. In confirmation of this, formula (\ref{def:I_F}) is not the
rate function without the contribution of the backward recurrence
time, which is the rate function of $W_t$, i.e.\ the convex conjugate
of $-(\varphi\wedge\ell)$. The rate function (\ref{def:I_F}) of $F_t$
coincides with the rate function of $W_t$ only if $\varphi(k)\le\ell$
for all $k$, which holds trivially if $\ell=+\infty$. If instead
$\ell<+\infty$ and $\varphi(k)>\ell$ for some $k$, then the large
fluctuations of the incomplete reward contribute to the large
fluctuations of $F_t$.

Whatever the reference stochastic process is, we can always find a
renewal process that satisfies the condition $\varpi(k)\ge \varphi(k)$
for all $k\in\Rl$, thus bringing about a LDP for the additive
functional $F_t$. Indeed, this condition is automatically met if the
waiting times are designed to fulfill the exponential moment condition
$\Ex[e^{\rho M_1}]<+\infty$ for all $\rho>0$ with
$M_1:=\int_0^{S_1}|f(B_{1,\tau})|\,d\tau$.  As $|X_1|\le M_1$, this
condition also implies that the rate function is good. Therefore, a
full LDP with a good rate function can always be achieved by means of
a resetting protocol adapted to the additive functional $F_t$.  This
fact is finally verified in Appendix \ref{corol:proof}, where the
following corollary of Theorem \ref{th:LDP} is demonstrated.
\begin{corollary}
  \label{corol}
If $\Ex[e^{\rho M_1}]<+\infty$ for all $\rho>0$, then $F_t$ satisfies
a full LDP with the good rate function $I$ given by
(\ref{def:I_F}). Moreover, its asymptotic logarithmic moment
generating function $g$ exists and equals $-\varphi$.
\end{corollary}

There are additive functionals that fulfill $\varpi(k)\ge \varphi(k)$
for all $k\in\Rl$ for any renewal process. They have the property that
$\frac{1}{t}\ln\mathcal{E}_t(k)$ has a limit as $t$ is sent to
infinity for each $k$, meaning that the value $\mathcal{E}_t(k)$ of
the moment generating function of $\int_0^tf(B_{1,\tau})d\tau$ does
not oscillate too much at large $t$.  These additive functionals
satisfy a full LDP under any resetting mechanism, although the rate
function may not be good.  We anticipate that the positive occupation
time, the area, and the absolute area of the Brownian motion belong to
this class of additive functionals. To get insights into the issue we
observe that $\Ex[e^{\varphi(k)S_1+kX_1}]\le 1$ if
$\varphi(k)>-\infty$. In fact, by definition of $\varphi(k)$ there
exists a sequence $\{\zeta_i\}_{i\ge 1}$ with the property that
$\lim_{i\uparrow\infty}\zeta_i=\varphi(k)$ and $\Ex[e^{\zeta_i
    S_1+kX_1}]\le 1$ for all $i$, so that the Fatou's lemma implies
$\Ex[e^{\varphi(k)S_1+kX_1}]\le\liminf_{i\uparrow\infty}\Ex[e^{\zeta_iS_1+kX_1}]\le
1$. The independence between $S_1$ and $\{B_{1,t}\}_{t\ge 0}$ allows
to recast the relationship $\Ex[e^{\varphi(k)S_1+kX_1}]\le 1$ as
  \begin{equation}
    \int_{[0,+\infty)}e^{\varphi(k)s}\mathcal{E}_s(k)P(ds)=\Ex\big[e^{\varphi(k)S_1+kX_1}\big]\le 1,
\label{varphi_finito}
  \end{equation}
where $P:=\prob[S_1\in\cdot\,]$ is the probability measure induced on
$[0,+\infty)$ by $S_1$. Thus, while the inequality $\varpi(k)\ge
  \varphi(k)$ is trivial if $\varphi(k)=-\infty$, when
  $\varphi(k)>-\infty$ it reads
  $\limsup_{t\uparrow+\infty}\frac{1}{t}\ln\{e^{\varphi(k)t}\mathcal{E}_t(k)\prob[S_1>t]\}\le
  0$ conditional on (\ref{varphi_finito}). This push the idea that
  $\varpi(k)\ge\varphi(k)$ is likely to be satisfied provided that
  $\mathcal{E}_t(k)$ is not affected by too large oscillations as $t$
  goes on.  The following corollary of Theorem \ref{th:LDP} makes
  these arguments rigorous. The proof is provided in Appendix
  \ref{proof:corol_lim}.
\begin{corollary}
  \label{corol_lim}
Assume that for all $k\in\Rl$ the limit
$\lim_{t\uparrow+\infty}\frac{1}{t}\ln\mathcal{E}_t(k)$ exists and
that $\frac{1}{t}\ln\mathcal{E}_t(k)$ is eventually non-decreasing
with respect to $t$ if such limit is infinite.  Then, $F_t$ satisfies
a full LDP with the rate function $I$ given by
(\ref{def:I_F}). Moreover, its asymptotic logarithmic moment
generating function $g$ exists and equals $-\varphi$.
\end{corollary}

We conclude the section by drawing a comparison with the
G\"artner-Ellis theory. The G\"artner-Ellis theory aims to obtain a
LDP based on the sole knowledge of the asymptotic logarithmic moment
generating function $g$.  Supposing that $\varpi(k)\ge\varphi(k)$ for
all $k\in\Rl$, part $(iii)$ of Theorem \ref{th:LDP} yields
$g=-\varphi$, so that $g$ is convex and lower semicontinuous. Set
$k_-:=\inf\{k\in\Rl:g(k)<+\infty\}$ and
$k_+:=\sup\{k\in\Rl:g(k)<+\infty\}$ and notice that $g(k)<+\infty$ for
all $k\in(k_-,k_+)$ by convexity.  The G\"artner-Ellis theorem gives a
full LDP for $F_t$ with the convex conjugate of $g$ as the rate
function if the following additional conditions are met
\cite{DemboBook,Frank_LDP}:
 \begin{enumerate}[(i)]
  \item $k_-<0<k_+$;
  \item $g$ is differentiable throughout the open interval $(k_-,k_+)$;
  \item $g$ is steep, i.e.\ $\lim_{k\downarrow k_-}g'(k)=-\infty$ if
    $k_->-\infty$ and $\lim_{k\uparrow k_+}g'(k)=+\infty$ if
    $k_+<+\infty$.
 \end{enumerate}
These conditions qualify the lower semicontinuous convex function $g$
as \textit{essentially smooth}. Essentially smoothness of the
asymptotic logarithmic moment generating function is crucial in the
G\"artner-Ellis theory to deduce the lower large deviation bound via
an exponential change of measure.  Thus, the use of the
G\"artner-Ellis theorem without essentially smoothness of $g$ is
definitely incorrect, although its conclusions may be true. In fact,
Theorem \ref{th:LDP} shows that essentially smoothness of $g=-\varphi$
is not necessary at all for the conclusions of the G\"artner-Ellis
theorem to hold. Theorem \ref{th:LDP} is ultimately based on
Cram\'er's theory, whose modern construction relies on a sub-additive
argument that has nothing to do with exponential changes of measure
\cite{DemboBook,Frank_LDP}.

A general property of convex conjugation states that if $g=-\varphi$
is essentially smooth, then the rate function $I$ given by
(\ref{def:I_F}) is strictly convex on the effective domain where it is
finite (see \cite{Rockbook}, Theorem 26.3). Therefore, all models that
involve a non-strictly convex rate function cannot be faced with the
G\"artner-Ellis theory. These models are exactly the models that
exhibit a dynamical phase transition. A \textit{dynamical phase
  transition} is a point in the effective domain of $I$ that marks the
beginning of an affine stretch in its graph, thus breaking strictly
convexity of $I$. Summing up, no model that exhibits a dynamical phase
transition falls within the scope of the G\"artner-Ellis theorem.  We
shall find several dynamical phase transitions in the context of the
reset Brownian motion, which are perfectly tackled by Theorem
\ref{th:LDP}.

\section{Applications to the Brownian motion}
\label{sec:applications}

In this section we discuss applications of Theorems \ref{th:SLLN_F},
\ref{th:CLT_F}, and \ref{th:LDP} to the positive occupation time, the
area, and the absolute area of the reset Brownian motion. From now on,
$\{B_{1,t}\}_{t\ge 0},\{B_{2,t}\}_{t\ge 0},\ldots$ are assumed to be
standard Brownian motions. No restriction is imposed on the renewal
process defining the reset mechanism. The fundamental property of the
Brownian motion that allows explicit calculations is the fact that,
for each number $a\ge 0$, the process $\{B_{1,at}\}_{t\ge 0}$ is
distributed as $\{\sqrt{a}B_{1,t}\}_{t\ge 0}$. We shall tacitly make
use of this property of self-similarity many times.

Before going into specific issues, we provide an overall idea of what
assumptions on waiting times are required by Theorems \ref{th:SLLN_F}
and \ref{th:CLT_F} and by Corollary \ref{corol} to apply to an
additive functional of the reset Brownian motion with function $f$ of
at most polynomial growth. In fact, the following lemma states
conditions on the waiting times that imply $\Ex[M_1^\rho]<+\infty$ for
some $\rho>0$ or $\Ex[e^{\rho M_1}]<+\infty$ for all $\rho>0$. The
proof is presented in Appendix \ref{proof:cases}.
\begin{lemma}
  \label{lem:cases}
  Assume that $|f(z)|\le A(1+|z|^\alpha)$ for all $z\in\Rl$ with
  constants $A>0$ and $\alpha\ge 0$. The following conclusions hold:
\begin{enumerate}[$(i)$]
\item if $\Ex[S_1^{\rho(1+\alpha/2)}]<+\infty$ for some $\rho>0$, then $\Ex[M_1^\rho]<+\infty$;
\item if $\alpha\in[0,2)$ and $\limsup_{s\uparrow+\infty} s^{-\frac{2+\alpha}{2-\alpha}}\ln\prob[S_1>s]=-\infty$,
  then $\Ex[e^{\rho M_1}]<+\infty$ for all $\rho>0$.
\end{enumerate}
\end{lemma}

Suppose that $|f(z)|\le A(1+|z|^\alpha)$ for all $z\in\Rl$. Theorems
\ref{th:SLLN_F} and \ref{th:CLT_F} in combination with Lemma
\ref{lem:cases} show that the SLLN and the CLT hold for $F_t:=\int_0^t
f(Z_\tau)\,d\tau$ provided that $\Ex[S_1^{1+\alpha/2}]<+\infty$ and
$\Ex[S_1^{2+\alpha}]<+\infty$, respectively. If moreover
$\alpha\in[0,2)$ and $\limsup_{s\uparrow+\infty}
  s^{-\frac{2+\alpha}{2-\alpha}}\ln\prob[S_1>s]=-\infty$, then
  Corollary \ref{corol} of Theorem \ref{th:LDP} and Lemma
  \ref{lem:cases} imply that $F_t$ satisfies a full LDP with the good
  rate function $I$ given by (\ref{def:I_F}). The restriction
  $\alpha<2$ is due to the Gaussian tails of the Brownian motion.  We
  have $\alpha=0$ for the positive occupation time and $\alpha=1$ for
  the area and the absolute area.

\subsection{Positive occupation time of the reset Brownian motion}
\label{sec:LDP_time}

The positive occupation time of the reset Brownian motion is the
additive functional
\begin{equation*}
  F_t:=\int_0^t\mathds{1}_{\{Z_\tau>0\}}d\tau.
\end{equation*}
The associated reward
$X_1:=\int_0^{S_1}\mathds{1}_{\{B_{1,\tau}>0\}}d\tau=S_1\int_0^1\mathds{1}_{\{B_{1,S_1\tau}>0\}}d\tau$
has the property that $X_1/S_1$ is statistically independent of $S_1$
and distributed as
$\int_0^1\mathds{1}_{\{B_{1,\tau}>0\}}d\tau\in[0,1]$. The L\'evy's
arcsine law for the Wiener process \cite{Levy} gives the distribution
of $\int_0^1\mathds{1}_{\{B_{1,\tau}>0\}}d\tau$, and hence of
$X_1/S_1$: for each $x\in[0,1]$
\begin{equation*}
\prob\bigg[\frac{X_1}{S_1}\le x\bigg]=\frac{2}{\pi}\arcsin\sqrt{x}.
\end{equation*}
We find $\Ex[X_1]=\Ex[S_1]/2$, $\Ex[S_1X_1]=\Ex[S_1^2]/2$, and
$\Ex[X_1^2]=3\Ex[S_1^2]/8$. As far as the typical behavior of $F_t$ is
concerned, Theorem \ref{th:SLLN_F} tells us that the positive
occupation time $F_t$ satisfies
\begin{equation*}
\lim_{t\uparrow+\infty}\frac{F_t}{t}=\frac{1}{2}~~~~~\prob\mbox{-a.s.}
\end{equation*}
provided that $\Ex[S_1]<+\infty$. If $\Ex[S_1^2]<+\infty$, then
Theorem \ref{th:CLT_F} gives for every $z\in\Rl$
\begin{equation*}
  \lim_{t\uparrow+\infty}\prob\bigg[\frac{F_t-t/2}{\sqrt{vt}}\le z\bigg]=\int_{-\infty}^z \frac{e^{-\frac{1}{2}x^2}}{\sqrt{2\pi}}\,dx
\end{equation*}
with
\begin{equation*}
v=\frac{1}{8}\frac{\Ex[S_1^2]}{\Ex[S_1]}.
\end{equation*}

Regarding the large fluctuations, we can invoke the L\'evy's arcsine
law again to get for every $t>0$ and $k\in\Rl$
\begin{equation*}
\mathcal{E}_t(k):=\Ex\Big[e^{k\int_0^t\mathds{1}_{\{B_{1,\tau}>0\}}d\tau}\Big]=\Ex\Big[e^{kt\int_0^1\mathds{1}_{\{B_{1,\tau}>0\}}d\tau}\Big]=\int_0^1\frac{e^{tkx}}{\pi\sqrt{x(1-x)}}\,dx.
\end{equation*}
Starting from this formula, it is a simple exercise to prove that
$\lim_{t\uparrow+\infty}\frac{1}{t}\ln\mathcal{E}_t(k)=0\vee k$ for
all $k$, which states in particular that
$\frac{1}{t}\ln\mathcal{E}_t(k)$ has a finite limit when $t$ is sent
to infinity.  Thus, Corollary \ref{corol_lim} of Theorem \ref{th:LDP}
yields the following LDP.
\begin{proposition}
The positive occupation time $F_t$ of the reset Brownian motion
satisfies a full LDP with the rate function $I$ given by
(\ref{def:I_F}) for any resetting protocol. Moreover, its asymptotic
logarithmic moment generating function exists and equals $-\varphi$.
\end{proposition}

The rest of the section is devoted to characterize $I$. All the
information about the large fluctuations of the positive occupation
time $F_t$ is contained in the function $\Phi$ that associates the
pair $(\zeta,k)\in\Rl^2$ with
\begin{equation}
    \Phi(\zeta,k):=\Ex\big[e^{\zeta S_1+kX_1}\big]=\int_0^1\frac{\Ex[e^{(\zeta+kx)S_1}]}{\pi\sqrt{x(1-x)}}\,dx.
\label{def:Phi}
\end{equation}
This function determines $\varphi$ through the formula
$\varphi(k)=\sup\{\zeta\in\Rl:\Phi(\zeta,k)\le 1\}$ and $I$ by convex
conjugation of $-\varphi$.  We recall that $\varphi$ is concave and
upper semicontinuous according to Proposition \ref{prop:conv}. We
notice that if
$\ell:=\lim_{s\uparrow+\infty}-\frac{1}{s}\ln\prob[S_1>s]<+\infty$,
then $\Phi(\zeta,k)=+\infty$ for $\zeta>\ell$, so that
$\varphi(k)\le\ell$ for all $k\in\Rl$. As a consequence, $F_t$ and the
associate cumulative reward $W_t$ share the same rate function,
meaning that the backward recurrence time cannot affect the large
fluctuations of the positive occupation time on an exponential
scale. The function $\varphi$ is finite and satisfies the symmetry
$\varphi(-k)=\varphi(k)+k$ for all $k\in\Rl$. In fact, on the one hand
$\Phi(0,k)\le 1$ for $k<0$ and $\Phi(-k,k)\le 1$ for $k\ge 0$, so that
$\varphi(k)\ge-0\vee k$ for any $k$ by definition, and on the other
hand $\Phi(\zeta,-k)=\Phi(\zeta-k,k)$. The symmetry of $\varphi$
endows $I$ with the symmetry $I(1/2-w)=I(1/2+w)$ for every
$w\in\Rl$. We have
$I(w):=\sup_{k\in\Rl^d}\{wk+\varphi(k)\}\ge\sup_{k\in\Rl^d}\{wk-0\vee
k\}=+\infty$ for $w<0$ or $w>1$, which is really expected as
$F_t/t\in[0,1]$.  The values of $I$ over the interval $[0,1]$
containing its effective domain can be computed explicitly for
Poissonian resetting.

\begin{example}
  Under Poissonian resetting, i.e.\ $\prob[S_1>s]=e^{-rs}$ for all
  $s>0$ with some resetting rate $r>0$, for each pair $(\zeta,k)$ we
  find
\begin{align}
  \nonumber
  \Phi(\zeta,k)=\int_0^1\frac{\Ex[e^{(\zeta+kx)S_1}]}{\pi\sqrt{x(1-x)}}\,dx
  &=\begin{cases}
  \displaystyle{\int_0^1\frac{r}{r-\zeta-kx}\,\frac{dx}{\pi\sqrt{x(1-x)}}} & \mbox{if }r-\zeta>0\vee k,\\
  +\infty & \mbox{otherwise}
  \end{cases}\\
  \nonumber
  &=\begin{cases}
  \displaystyle{\frac{r}{\sqrt{(r-\zeta)(r-\zeta-k)}}} & \mbox{if }r-\zeta>0\vee k,\\
  +\infty & \mbox{otherwise}.
  \end{cases}
\end{align}
Thus, for all $k\in\Rl$
\begin{equation*}
  \varphi(k)=r-\frac{k+\sqrt{k^2+4r^2}}{2}
\end{equation*}
and for every $w\in[0,1]$
\begin{equation*}
  I(w)=r\Big[1-2\sqrt{w(1-w)}\Big].
\end{equation*}
The rate function $I$ depends linearly on the resetting rate.  This
rate function has been previously determined in \cite{LDP_res_2} by
direct inspection at finite times $t$ of the distribution of the
positive occupation time. In fact, an explicit formula exists for this
distribution in the special case of Poissonian resetting.
\end{example}

In order to describe the rate function $I$ beyond Poissonian resetting
we need to examine in detail the properties of the function $\varphi$.
When $0<\ell<+\infty$ these properties are characterized by the two
extended real numbers
\begin{equation*}
\Lambda:=\Ex\bigg[\frac{e^{\ell S_1}}{\sqrt{1+S_1}}\bigg]
\end{equation*}
and
\begin{equation*}
\Xi:=\Ex\big[\sqrt{S_1}e^{\ell S_1}\big].
\end{equation*}
We let the following lemma to present the relevant features of
$\varphi$. The proof is reported in Appendix
\ref{proof:varphi_time}. The lemma shows in particular that the
asymptotic logarithmic moment generating function $g=-\varphi$ is not
differentiable in the interesting case $0<\ell<+\infty$,
$\Lambda<+\infty$, and $\Xi<+\infty$. In this case, the large
fluctuations of the positive occupation time of the reset Brownian
motion cannot be faced with the G\"artner-Ellis theory.
\begin{lemma}
  \label{lem:varphi_time}
  The following conclusions hold:
  \begin{enumerate}[$(i)$]
   \item if $\ell=+\infty$ or $0<\ell<+\infty$ and $\Lambda=+\infty$,
     then $\varphi(k)$ solves the equation $\Ex[e^{\varphi(k)
         S_1+kX_1}]=1$ for every $k\in\Rl$.  The function $\varphi$ is
     analytic throughout $\Rl$ with
     \begin{equation*}
       -\varphi'(k)=\frac{\Ex[X_1e^{\varphi(k)S_1+kX_1}]}{\Ex[S_1e^{\varphi(k)S_1+kX_1}]},
     \end{equation*}
     and the limits $\lim_{k\downarrow-\infty}-\varphi'(k)=0$ and
     $\lim_{k\uparrow+\infty}-\varphi'(k)=1$ hold true;
\item if $0<\ell<+\infty$ and $\Lambda<+\infty$, then there exists
  $\lambda>0$ such that $\Ex[e^{\ell S_1-\lambda X_1}]=1$ and
  $\varphi(k)$ solves the equation $\Ex[e^{\varphi(k) S_1+kX_1}]=1$
  for $|k|<\lambda$, whereas $\varphi(k)=\ell\wedge(\ell-k)$ for
  $|k|\ge\lambda$. The function $\varphi$ is analytic on the open
  interval $(-\lambda,\lambda)$ with
\begin{equation*}
       -\varphi'(k)=\frac{\Ex[X_1e^{\varphi(k)S_1+kX_1}]}{\Ex[S_1e^{\varphi(k)S_1+kX_1}]}.
     \end{equation*}
If $\Xi=+\infty$, then $\varphi$ is differentiable throughout
$\Rl$. If $\Xi<+\infty$, then $\varphi$ is not differentiable at
$-\lambda$ and $\lambda$ with finite left derivatives
\begin{equation*}
  -\varphi_-'(k)=\begin{cases}
  0 & \mbox{for }k=-\lambda,\\
  \lim_{k\uparrow\lambda}-\varphi'(k)=\frac{\Ex[X_1e^{\varphi(k)S_1+kX_1}]}{\Ex[S_1e^{\varphi(k)S_1+kX_1}]}<1& \mbox{for }k=\lambda
  \end{cases}
\end{equation*}
and finite right derivatives
\begin{equation*}
  -\varphi_+'(k)=\begin{cases}
  \lim_{k\downarrow\lambda}-\varphi'(k)=\frac{\Ex[X_1e^{\varphi(k)S_1+kX_1}]}{\Ex[S_1e^{\varphi(k)S_1+kX_1}]}>0 & \mbox{for }k=-\lambda,\\
  1& \mbox{for }k=\lambda;
  \end{cases}
\end{equation*}
\item  if $\ell=0$, then $-\varphi(k)=0\vee k$.
  \end{enumerate}  
\end{lemma}

We are now in the position to discuss the graph of $I$. We remind that
if $w$ is a subgradient of the convex function $-\varphi$ at the point
$k$, then $wh+\varphi(h)\le wk+\varphi(k)$ for all $h\in\Rl$ and
$I(w):=\sup_{h\in\Rl}\{wh+\varphi(h)\}\le wk+\varphi(k)$, which
implies $I(w)=wk+\varphi(k)$. We refer to \cite{Rockbook} for general
properties of convex functions that sometimes we use.  If
$\ell=+\infty$ or $0<\ell<+\infty$ and $\Lambda=+\infty$, then for
each given $w\in(0,1)$ there exists $k\in\Rl$ such that
$w=-\varphi'(k)$ as part $(i)$ of Lemma \ref{lem:varphi_time} states
that $-\varphi'$ increases continuously from $0$ at $k=-\infty$ to $1$
at $k=+\infty$. Thus, $w$ is a subgradient of $-\varphi$ at $k$ and we
have $I(w)=wk+\varphi(k)$. Then, in the case $\ell=+\infty$ or
$0<\ell<+\infty$ and $\Lambda=+\infty$, $I$ turns out to be analytic
on the open interval $(0,1)$ by the analytic implicit function
theorem. The values of $I$ at the points $w=0$ and $w=1$ are
determined by the limits $I(0)=\lim_{w\downarrow 0}I(w)$ and
$I(1)=\lim_{w\uparrow 1}I(w)$, which hold because $I$ is convex and
lower semicontinuous (see \cite{Rockbook}, Corollary 7.5.1). The same
conclusions are valid when $0<\ell<+\infty$, $\Lambda<+\infty$, and
$\Xi=+\infty$. In fact, in this case $-\varphi'$ increases
continuously from $0$ at $k=-\lambda$ to $1$ at $k=\lambda$ by part
$(ii)$ of Lemma \ref{lem:varphi_time}. The rate function is again
analytic on $(0,1)$ and, moreover, $I(0)=\varphi(-\lambda)$ as $0$ is
a subgradient of $-\varphi$ at $k=-\lambda$ and
$I(1)=\lambda+\varphi(\lambda)$ as $1$ is a subgradient of $-\varphi$
at $k=\lambda$. We stress that Poissonian resetting falls in the class
$0<\ell=r<+\infty$ and $\Lambda=+\infty$.

The graph of the rate function $I$ is characterized by two
singularities and two affine stretches when $0<\ell<+\infty$,
$\Lambda<+\infty$, and $\Xi<+\infty$. For instance, this occurs with
the waiting times described by the laws
\begin{equation*}
\prob[S_1>s]=\frac{e^{-s}}{1+s^\alpha}
\end{equation*}
and
\begin{equation*}
\prob[S_1>s]=e^{-s^\beta-s}
\end{equation*}
for all $s\ge 0$ with $\alpha>3/2$ and $\beta\in(0,1)$.  In fact,
according to part $(ii)$ of Lemma \ref{lem:varphi_time}, in this case
all points in $[0,w_-]$ with $w_-:=-\varphi_+'(-\lambda)$ are
subgradient of $-\varphi$ at $-\lambda$ and all points in $[w_+,1]$
with $w_+:=-\varphi_-'(\lambda)$ are subgradient of $-\varphi$ at
$\lambda$ (see \cite{Rockbook}, Theorem 23.2). Thus,
$I(w)=-w\lambda+\varphi(-\lambda)$ for all $w\in[0,w_-]$ and
$I(w)=w\lambda+\varphi(\lambda)$ for all $w\in[w_+,1]$, so that two
affine stretches emerge. Moreover, as $-\varphi'$ increases
continuously from $w_-$ at $k=-\lambda$ to $w_+$ at $k=\lambda$, for
each $w\in(w_-,w_+)$ there exists $k\in(-\lambda,\lambda)$ such that
$w=-\varphi'(k)$, and $I(w)=wk+\varphi(k)$ follows. The rate function
$I$ is analytic over the open interval $(w_-,w_+)$ by the analytic
implicit function theorem. It is clear that $w_-<\mu=1/2$ and
$w_+>\mu$ as $I(\mu)=0$ by Theorem \ref{th:LDP}.  In conclusion, the
points $w_-$ and $w_+$ are singular points for $I$ that mark some form
of dynamical phase transition.

Dynamical phase transitions under resetting have been previously
documented for the total displacement of certain random walks
\cite{LDP_res_5}. The same type of singular behavior has been found
for the macroscopic observables of renewal models of statistical
mechanics, such as the pinning model of polymers or the
Poland-Scheraga model, under equilibrium conditions
\cite{Models,Giambattista_bj}.  All these models share a common
regenerative structures, whether they are kinetic models with few
degrees of freedom or equilibrium models with many degrees of
freedom. Although the characterization of these phase transitions is
beyond the scope of the present work, we can try to give an
explanation by drawing an analogy with sums of independent and
identical distributed random variables and the phenomenon of
condensation of their fluctuations
\cite{Condensation_1,Condensation_2,Condensation_3,Condensation_4}. After
all, $F_t$ is close to the cumulative reward $W_t$, which is a random
sum of independent and identical distributed random variables.  If
$F_t$ were exactly a sum of independent and identical distributed
random variables, then we would conclude that a fluctuation $w$ of
$F_t/t$ above $w_+>\mu$ is realized by the combination of two
different mechanisms: many small deviations of the waiting times all
in the same direction for a partial excursion up to $w_+$, and a big
jump of only one of the waiting times for the remaining part $w-w_+$
of the fluctuation. The same would occur for a fluctuation of $F_t/t$
below $w_-<\mu$. We believe that this picture also describes the
dynamical phase transitions of the positive occupation time of the
reset Brownian motion.

Finally, the case $\ell=0$ corresponds to a heavy-tailed waiting time
and gives $I(w)=0$ for all $w\in[0,1]$. This rate function simply
tells us that the decay with time of the probability of a large
fluctuation of the positive occupation time is slower than exponential
for heavy-tailed waiting times. In this case, the description of the
large fluctuations requires more precise asymptotic results for
specific classes of waiting times. Still regarding renewal-reward
processes as primary tools, we mention that precise large deviation
principles for renewal-reward processes have been established when
rewards are independent of waiting times for several types of
heavy-tailed waiting times and rewards \cite{L1,L2,L3,L4,L5,L6}. Some
forms of dependence between rewards and waiting times have been
considered in \cite{L8,critical}.

\subsection{Area of the reset Brownian motion}
\label{sec:LDP_area}

The area of the reset Brownian motion is the additive functional
\begin{equation*}
  F_t:=\int_0^tZ_\tau\,d\tau.
\end{equation*}
The corresponding reward $X_1:=\int_0^{S_1}B_{1,\tau}\,d\tau$ enjoys
the property that $\sqrt{3/S_1^3}X_1$ is statistically independent of
$S_1$ and distributed as $\sqrt{3}\int_0^1B_{1,\tau}\,d\tau$, which in
turn is distributed as a Gaussian random variable with mean $0$ and
variance $1$ \cite{Janson}.  In order to characterize the typical
behavior of the area $F_t$ we observe that $\Ex[X_1]=0$ if
$\Ex[S_1^{3/2}]<+\infty$ and $\Ex[X_1^2]=\Ex[S_1^3]/3$ if
$\Ex[S_1^3]<+\infty$.  Furthermore, since
$\Ex\big[\int_0^1|B_{1,\tau}|\,d\tau\big]=\sqrt{8/9\pi}$ and
$\Ex\big[\{\int_0^1|B_{1,\tau}|\,d\tau\}^2\big]=3/8$ \cite{Janson}, we
realize that
$\Ex[M_1]=\Ex\big[\int_0^{S_1}|B_{1,\tau}|\,d\tau\big]=\Ex\big[S_1^{3/2}\int_0^1|B_{1,\tau}|\,d\tau\big]=\sqrt{8/9\pi}\,\Ex[S_1^{3/2}]$
and, similarly, that $\Ex[M_1^2]=(3/8)\,\Ex[S_1^3]$. Thus, Theorem
\ref{th:SLLN_F} states that $F_t$ satisfies
\begin{equation*}
\lim_{t\uparrow+\infty}\frac{F_t}{t}=0~~~~~\prob\mbox{-a.s.}
\end{equation*}
whenever $\Ex[S_1^{3/2}]<+\infty$.  Under the hypothesis
$\Ex[S_1^3]<+\infty$, Theorem \ref{th:CLT_F} shows that for every
$z\in\Rl$
\begin{equation*}
  \lim_{t\uparrow+\infty}\prob\bigg[\frac{F_t}{\sqrt{vt}}\le z\bigg]=\int_{-\infty}^z \frac{e^{-\frac{1}{2}x^2}}{\sqrt{2\pi}}\,dx
\end{equation*}
with
\begin{equation*}
v=\frac{1}{3}\frac{\Ex[S_1^3]}{\Ex[S_1]}.
\end{equation*}
As mentioned at the end of Section \ref{sec:F_normal}, the hypothesis
$\Ex[M_1^2]<+\infty$ is the same as $\Ex[X_1^2]<+\infty$ although the
function $f(z):=z$ that defines the area takes both positive and
negative values.

Moving to the large fluctuations of $F_t$, we observe that for all
$t\ge 0$ and $k\in\Rl$
\begin{equation*}
  \mathcal{E}_t(k):=\Ex\Big[e^{k\int_0^tB_{1,\tau}\,d\tau}\Big]=\Ex\Big[e^{kt^{3/2}\int_0^1B_{1,\tau}\,d\tau}\Big]=e^{\frac{1}{6}k^2t^3}.
\end{equation*}
Thus, $\frac{1}{t}\ln\mathcal{E}_t(k)=\frac{(kt)^2}{6}$ is
non-decreasing with respect to $t$ and Corollary \ref{corol_lim} of
Theorem \ref{th:LDP} gives the following result.
\begin{proposition}
The area $F_t$ of the reset Brownian motion satisfies a full LDP with
the rate function $I$ given by (\ref{def:I_F}) for any resetting
protocol. Moreover, its asymptotic logarithmic moment generating
function exists and equals $-\varphi$.
\end{proposition}

The function $\varphi$ that defines $I$ by convex conjugation now
reads
\begin{equation*}
\varphi(k):=\sup\Big\{\zeta\in\Rl~:~\Ex\big[e^{\zeta S_1+kX_1}\big]\le 1\Big\}=\sup\Big\{\zeta\in\Rl~:~\Ex\big[e^{\zeta S_1+\frac{1}{6}k^2S_1^3}\big]\le 1\Big\}
\end{equation*}
and satisfies the symmetry $\varphi(-k)=\varphi(k)$. This symmetry
gives $I(-w)=I(w)$ for all $w\in\Rl$.  We point out that if
$\ell<+\infty$, then $\varphi(k)=0$ for $k=0$ and $\varphi(k)=-\infty$
for $k\ne 0$. Thus, as in the case of the positive occupation time, we
have $\varphi(k)\le\ell$ for all $k\in\Rl$, showing that $F_t$ and
$W_t$ satisfies a full LDP with the same rate function. However, if
$\ell<+\infty$, then $I(w)=0$ for every $w$ and the probability of a
large fluctuation decays slower than exponential in time.

We want to characterize the rate function $I$, and a preliminary study
of $\varphi$ is needed to this aim. The relevant rate is not $\ell$,
as for the positive occupation time, but rather $r$ defined by
\begin{equation*}
  r:=\liminf_{s\uparrow+\infty}-\frac{1}{s^3}\ln\prob[S_1>s].
\end{equation*}
When $r$ is finite, we set for brevity $\xi:=\varphi(\sqrt{6r})$. When
$r$ and $\xi$ are finite, we introduce the real number
\begin{equation*}
  \Lambda:=\Ex\big[e^{\xi S_1+rS_1^3}\big]=\Ex\big[e^{\varphi(\sqrt{6r}) S_1+\sqrt{6r}X_1}\big]\le 1
\end{equation*}
and the extended real number
\begin{equation*}
  \Xi:=\Ex\big[S_1^3e^{\xi S_1+rS_1^3}\big].
\end{equation*}
The following lemma, which is proved in Appendix
\ref{proof:varphi_area}, collects the main features of the function
$\varphi$ and puts $\xi$, $\Lambda$, and $\Xi$ into context. It
demonstrates in particular that the asymptotic logarithmic moment
generating function $g=-\varphi$ is not steep in the interesting case
$0<r<+\infty$, $\xi>-\infty$, $\Lambda=1$, and $\Xi<+\infty$. The
G\"artner-Ellis theory does not apply in this case.
\begin{lemma}
  \label{lem:varphi_area}
  The following conclusions hold:
  \begin{enumerate}[$(i)$]
   \item if $r=+\infty$, then $\varphi(k)$ solves the equation
     $\Ex[e^{\varphi(k) S_1+\frac{1}{6}k^2S_1^3}]=1$ for every
     $k\in\Rl$.  The function $\varphi$ is analytic throughout $\Rl$
     with
     \begin{equation*}
       -\varphi'(k)=\frac{k}{3}\frac{\Ex[S_1^3e^{\varphi(k)S_1+\frac{1}{6}k^2S_1^3}]}{\Ex[S_1e^{\varphi(k)S_1+\frac{1}{6}k^2S_1^3}]},
     \end{equation*}
     and the limits $\lim_{k\downarrow-\infty}-\varphi'(k)=-\infty$ and
     $\lim_{k\uparrow+\infty}-\varphi'(k)=+\infty$ hold true;
\item if $0<r<+\infty$, then $\varphi(k)$ solves the equation
  $\Ex[e^{\varphi(k) S_1+\frac{1}{6}k^2S_1^3}]=1$ for $|k|<\sqrt{6r}$,
  whereas $\varphi(k)=-\infty$ for $|k|>\sqrt{6r}$. The function
  $\varphi$ is analytic on the open interval $(-\sqrt{6r},\sqrt{6r})$
  with
\begin{equation*}
       -\varphi'(k)=\frac{k}{3}\frac{\Ex[S_1^3e^{\varphi(k)S_1+\frac{1}{6}k^2S_1^3}]}{\Ex[S_1e^{\varphi(k)S_1+\frac{1}{6}k^2S_1^3}]}.
\end{equation*}
The derivative $\varphi'$ has the limits
\begin{equation*}
  \lim_{k\downarrow-\sqrt{6r}}-\varphi'(k)=
  \begin{cases}
     -\sqrt{\frac{2r}{3}}\frac{\Ex[S_1^3e^{\xi S_1+rS_1^3}]}{\Ex[S_1e^{\xi S_1+rS_1^3}]} & \mbox{if $\xi>-\infty$, $\Lambda=1$, and $\Xi<+\infty$},\\
    -\infty & \mbox{otherwise}
    \end{cases}
\end{equation*}
and
\begin{equation*}
  \lim_{k\uparrow\sqrt{6r}}-\varphi'(k)=
  \begin{cases}
    \sqrt{\frac{2r}{3}}\frac{\Ex[S_1^3e^{\xi S_1+rS_1^3}]}{\Ex[S_1e^{\xi S_1+rS_1^3}]} & \mbox{if $\xi>-\infty$, $\Lambda=1$, and $\Xi<+\infty$},\\
    +\infty & \mbox{otherwise};
    \end{cases}
\end{equation*}

\item if $r=0$, then $\varphi(k)=0$ for $k=0$ and $\varphi(k)=-\infty$
  for $k\ne 0$.
  \end{enumerate}  
\end{lemma}

We can now discuss the graph of $I$. If $r=+\infty$, or $0<r<+\infty$
and $\xi=-\infty$, or $0<r<+\infty$ and $\xi>-\infty$ and $\Lambda<1$,
or $0<r<+\infty$ and $\xi>-\infty$ and $\Lambda=1$ and $\Xi=+\infty$,
then parts $(i)$ and $(ii)$ of Lemma \ref{lem:varphi_area} state that
the convex function $-\varphi$ is differentiable with continuously
increasing derivative from $-\infty$ at $k=-\sqrt{6r}$ to $+\infty$ at
$k=\sqrt{6r}$. We agree that $-\sqrt{6r}=-\infty$ and
$\sqrt{6r}=+\infty$ when $r=+\infty$.  As a consequence, for each
$w\in\Rl$ there exists $k\in(-\sqrt{6r},\sqrt{6r})$ such that
$w=-\varphi'(k)$, and $I(w)=wk+\varphi(k)$ follows. The rate function
$I$ turns out to be analytic throughout $\Rl$ and good in all these
cases.

Singularities emerge when $0<r<+\infty$, $\xi>-\infty$, $\Lambda=1$,
and $\Xi<+\infty$, such as for the waiting times of Example
\ref{example:area} below. By setting
\begin{equation*}
w_\pm:=\pm\sqrt{\frac{2r}{3}}\frac{\Ex[S_1^3e^{\xi S_1+rS_1^3}]}{\Ex[S_1e^{\xi S_1+rS_1^3}]}
\end{equation*}
and by invoking part $(ii)$ of Lemma \ref{lem:varphi_area}, we have
that for each $w\in(w_-,w_+)$ there exists
$k\in(-\sqrt{6r},\sqrt{6r})$ such that $w=-\varphi'(k)$, in such a way
that $I(w)=wk+\varphi(k)$. The rate function is analytic on the open
interval $(w_-,w_+)$. For $w\ge w_+$, the function that maps
$k\in[-\sqrt{6r},\sqrt{6r}]$ in $wk+\varphi(k)$ is non-decreasing,
whereas $\varphi(k)=-\infty$ for
$k\notin[-\sqrt{6r},\sqrt{6r}]$. Thus,
$I(w)=\sup_{k\in\Rl}\{wk+\varphi(k)\}=w\sqrt{6r}+\varphi(\sqrt{6r})=w\sqrt{6r}+\xi$
for all $w\ge w_+$. Similarly, $I(w)=-w\sqrt{6r}+\xi$ for all $w\le
w_-$. The graph of $I$ is now characterized by two affine stretches
with extremes $w_-$ and $w_+$, which are associated with some
dynamical phase transitions. As for the positive occupation time, we
believe that the mechanism underlying these dynamical phase
transitions is the phenomenon of condensation of fluctuations. We
leave a quantitative description as an open problem to be addressed in
future research.

To conclude, we observe that $I(w)=0$ for all $w\in\Rl$ if $r=0$.
This case, which includes Poissonian resetting as already observed in
\cite{LDP_res_2}, needs more precise asymptotic results to be devised
for specific classes of waiting times.

\begin{example}
  \label{example:area}
  The simplest waiting time distribution that brings about a full LDP
  with good rate function for the area is the super-exponential
  distribution
\begin{equation*}
  \prob[S_1>s]=e^{-rs^3}
\end{equation*}
for all $s\ge 0$ with some ``resetting rate'' $r>0$.  We discuss this
example in detail. While $\varphi(k)=-\infty$ for $|k|>\sqrt{6r}$, for
$k\in[-\sqrt{6r},\sqrt{6r}]$ it turns out that $\varphi(k)$ solves the
equation
\begin{equation}
 \Ex\big[e^{\varphi(k) S_1+\frac{1}{6}k^2S_1^3}\big]=\int_0^{+\infty}3rs^2e^{\varphi(k) s+(\frac{1}{6}k^2-r)s^3}ds=1.
\label{example:area_1}
\end{equation}
We find $\xi:=\varphi(\sqrt{6r})=-(6r)^{\frac{1}{3}}>-\infty$,
$\Lambda=1$, and $\Xi=10\,r^{-2}<+\infty$. Thus, two singularities
and two affine stretches occur. In fact, for every $w\in\Rl$ the rate
function reads
\begin{equation*}
  I(w)=\begin{cases}
  -w\sqrt{6r}+\xi & \mbox{if }w\le w_-,\\
  wk+\varphi(k) & \mbox{if }w\in(w_-,w_+),\\
   w\sqrt{6r}+\xi & \mbox{if }w\ge w_+
  \end{cases}
\end{equation*}
with $w_\pm=\pm\frac{20}{3}(6r)^{-\frac{1}{6}}$ and
$k\in(-\sqrt{6r},\sqrt{6r})$ such that $w+\varphi'(k)=0$.  

This example offers a possibility to understand how the rate function
behaves in the limit of zero resetting. Let us write $\varphi_r$ and
$I_r$ in place of $\varphi$ and $I$ to stress the dependence on
$r$. We already know that $I_0=0$. The change of variable $s\mapsto
r^{-\frac{1}{3}}s$ in (\ref{example:area_1}) shows that
$\varphi_r(k)=r^{\frac{1}{3}}\varphi_1(r^{-\frac{1}{2}}k)$ for all
$k$. It follows that for every $w\in\Rl$
\begin{equation*}
I_r(w):=\sup_{k\in\Rl}\big\{wk+\varphi_r(k)\big\}=r^{\frac{1}{3}}\sup_{k\in\Rl}\big\{r^{\frac{1}{6}}wk+\varphi_1(k)\big\}=r^{\frac{1}{3}}I_1(r^{\frac{1}{6}}w).
\end{equation*}
The analysis of $I_1$ in a neighborhood of the origin is a simple
exercise, which allows one to conclude that there exists a positive
constant $C$ such that
\begin{equation*}
  \bigg|I_1(w)-\frac{\Gamma(1/3)}{2}w^2\bigg|\le Cw^4
\end{equation*}
for any $w$, $\Gamma$ being the Euler gamma function. Thus, for all
$r>0$ and $w\in\Rl$ we find
\begin{equation*}
  \bigg|I_r(w)-\frac{\Gamma(1/3)}{2}r^{\frac{2}{3}}w^2\bigg|\le Crw^4, 
\end{equation*}
which states that $I_r(w)$ approaches the value $0$ with speed
$r^{\frac{2}{3}}$ when $r$ goes to 0.
\end{example}

\subsection{Absolute area of the reset Brownian motion}

The absolute area of the reset Brownian motion is the additive
functional
\begin{equation*}
  F_t:=\int_0^t|Z_\tau|\,d\tau
\end{equation*}
with reward $X_1:=\int_0^{S_1}|B_{1,\tau}|\,d\tau$. Repeating some of
the arguments of the beginning of Section \ref{sec:LDP_area}, we find
$\Ex[X_1]=\sqrt{8/9\pi}\,\Ex[S_1^{3/2}]$,
$\Ex[X_1^2]=(3/8)\,\Ex[S_1^3]$, and
$\Ex[S_1X_1]=\sqrt{8/9\pi}\,\Ex[S_1^{5/2}]$.  Then, regarding the
typical behavior of the absolute area $F_t$, Theorem \ref{th:SLLN_F}
states that
\begin{equation*}
\lim_{t\uparrow+\infty}\frac{F_t}{t}=\mu=\sqrt{\frac{8}{9\pi}}\frac{\Ex[S_1^{3/2}]}{\Ex[S_1]}~~~~~\prob\mbox{-a.s.}
\end{equation*}
if $\Ex[S_1^{3/2}]<+\infty$. When $\Ex[S_1^3]<+\infty$, Theorem
\ref{th:CLT_F} tells us that for all $z\in\Rl$
\begin{equation*}
  \lim_{t\uparrow+\infty}\prob\bigg[\frac{F_t-\mu t}{\sqrt{vt}}\le z\bigg]=\int_{-\infty}^z \frac{e^{-\frac{1}{2}x^2}}{\sqrt{2\pi}}\,dx
\end{equation*}
with
\begin{equation*}
  v=\frac{3}{8}\frac{\Ex[S_1^3]}{\Ex[S_1]}-\mu\sqrt{\frac{32}{9\pi}}\frac{\Ex[S_1^{5/2}]}{\Ex[S_1]}+
  \mu^2\frac{\Ex[S_1^2]}{\Ex[S_1]}.
\end{equation*}

The study of the large fluctuations of $F_t$ requires to introduce
some further mathematical details. The \textit{Airy function}
$\mbox{Ai}$ is the analytic function that maps $x\in\Rl$ in
\begin{equation*}
  \mbox{Ai}(x):=\frac{1}{\pi}\int_0^{+\infty} \cos\bigg(xu+\frac{u^3}{3}\bigg)du.
\end{equation*}
The derivative $\mbox{Ai}'$ of $\mbox{Ai}$ has countably many simple
zeros $0>z_1>z_2>\cdots$ on the negative real axis \cite{Janson,Airy}. For
each $i\ge 1$ we set $\nu_i:=2^{-1/3}|z_i|$ and
\begin{equation*}
  c_i:=-\frac{1+3\int_{z_i}^0\mbox{Ai}(x)dx}{3\mbox{Ai}''(z_i)}=\frac{1+3\int_{z_i}^0\mbox{Ai}(x)dx}{3|z_i|\mbox{Ai}(z_i)}.
\end{equation*}
We have $\nu_1=0.80861\ldots$ and $c_1=1.48257\ldots$ \cite{Airy}. The
asymptotics of $\nu_i$ and $c_i$ is characterized by the limits
$\lim_{i\uparrow\infty}(2\pi i)^{-2/3}\nu_i=2$ and
$\lim_{i\uparrow\infty}(-1)^i\sqrt{3i/2}\,c_i=-1$ \cite{Airy}.  Our
interest in the Airy function stems from the fact that the moment
generating function of $\int_0^1|B_{1,\tau}|\,d\tau$ takes the value
\cite{Janson,Airy}
\begin{equation}
\Ex\Big[e^{k\int_0^1|B_{1,\tau}|\,d\tau}\Big]=\sum_{i=1}^\infty c_ie^{-\nu_i|k|^{2/3}}
\label{Airy_exp}
\end{equation}
for $k<0$.  The following lemma is proved in Appendix
  \ref{proof:auxxx} and completes the picture for $k\ge 0$.
\begin{lemma}
  \label{lem:auxxx}
  There exists a positive constant $L$ such that for all $k\ge 0$
  \begin{equation*}
    e^{\frac{1}{6}k^2}\le\Ex\Big[e^{k\int_0^1|B_{1,\tau}|\,d\tau}\Big]\le L\,e^{\frac{1}{6}k^2}.
    \end{equation*}
\end{lemma}

Formula \ref{Airy_exp} shows that for all $t>0$ and $k<0$
\begin{equation*}
\mathcal{E}_t(k):=\Ex\Big[e^{k\int_0^t|B_{1,\tau}|\,d\tau}\Big]=\Ex\Big[e^{kt^{3/2}\int_0^1|B_{1,\tau}|\,d\tau}\Big]=\sum_{i=1}^\infty c_ie^{-t\nu_i|k|^{2/3}},
\end{equation*}
so that
$\lim_{t\uparrow+\infty}\frac{1}{t}\ln\mathcal{E}_t(k)=-\nu_1|k|^{2/3}$
exists and is finite for $k\le 0$. Lemma \ref{lem:auxxx} entails
$\lim_{t\uparrow+\infty}\frac{1}{t}\ln\mathcal{E}_t(k)=+\infty$ for
every $k>0$.  We claim that $\frac{1}{t}\ln\mathcal{E}_t(k)$ is
non-decreasing with respect to $t$ for $k>0$. To show this, fix $k>0$
and $q>1$. H\"older's inequality gives for all $t>0$
\begin{equation*}
  \big\{\mathcal{E}_t(k)\big\}^q=\bigg\{\Ex\Big[e^{k\int_0^t|B_{1,\tau}|d\tau}\Big]\bigg\}^q\le\Ex\Big[e^{qk\int_0^t|B_{1,\tau}|d\tau}\Big]=
  \Ex\Big[e^{k\int_0^{q^{2/3}t}|B_{1,\tau}|d\tau}\Big]=\mathcal{E}_{q^{2/3}t}(k).
\end{equation*}
Choosing $q=(s/t)^{3/2}$ with $s>t$ and observing that
$\ln\mathcal{E}_s(k)\ge 0$ we find
\begin{equation*}
  \frac{1}{t}\ln\mathcal{E}_t(k)\le\sqrt{\frac{t}{s}}\frac{1}{s}\ln\mathcal{E}_{s}(k)\le\frac{1}{s}\ln\mathcal{E}_{s}(k).
\end{equation*}
We have thus verified the hypothesis of Corollary \ref{corol_lim} of
Theorem \ref{th:LDP}, deducing the following LDP.
\begin{proposition}
The absolute area $F_t$ of the reset Brownian motion satisfies a full
LDP with the rate function $I$ given by (\ref{def:I_F}) for any
resetting protocol. Moreover, its asymptotic logarithmic moment
generating function exists and equals $-\varphi$.
\end{proposition}

In contrast to the cases of the positive occupation time and the area,
the rate functions of $F_t$ and $W_t$ can differ for the absolute
area, meaning that the backward recurrence time can affect the large
fluctuations of $F_t$. For instance, this occur for Poissonian
resetting, where $\varphi(k)>\ell=r$ for all $k<0$ sufficiently large
in magnitude. Indeed, in the case of Poissonian resetting, formula
(\ref{Airy_exp}) shows that for $k<0$ and $\zeta<r+\nu_1|k|^{2/3}$
\begin{align}
  \nonumber
  \Ex\big[e^{\zeta S_1+kX_1}\big]&=\int_0^{+\infty} \Ex\Big[e^{ks^{3/2}\int_0^1|B_{1,\tau}|\,d\tau}\Big] re^{\zeta s-rs}ds\\
  \nonumber
  &=\int_0^{+\infty}\sum_{i=1}^\infty rc_ie^{(\zeta-r-\nu_i|k|^{2/3})s}ds=\sum_{i=1}^\infty \frac{rc_i}{r+\nu_i|k|^{2/3}-\zeta}.
\end{align}
The last equality is justified by the fact that $|c_i|/\nu_i$ goes to
zero as $i^{-7/6}$ when $i$ is sent to infinity. Thus, $\Ex[e^{\zeta
    S_1+kX_1}]<+\infty$ for $k<0$ and $\zeta<r+\nu_1|k|^{2/3}$ with
the limit $\lim_{\zeta\uparrow r+\nu_1|k|^{2/3}}\Ex[e^{\zeta
    S_1+kX_1}]=+\infty$, and $\Ex[e^{r
    S_1+kX_1}]=|k|^{-2/3}\sum_{i=1}^\infty(rc_i/\nu_i)<1$ for $k<0$
sufficiently large in magnitude. It follows that $\varphi(k)$ is for
each $k<0$ the unique real number $\zeta$ that solves the equation
$\Ex[e^{\zeta S_1+kX_1}]=1$, and that $\varphi(k)>r$ for all $k<0$
sufficiently large in magnitude.

The rate function $I$ of the absolute area satisfies $I(w)=+\infty$
for $w<0$ because the absolute area takes non-negative values. For a
formal proof it is enough to observe that $\varphi(k)\ge 0$ for all
$k\le 0$ by definition, so that
$I(w):=\sup_{k\in\Rl}\{wk+\varphi(k)\}\ge\sup_{k\le 0}\{wk\}=+\infty$
for $w<0$. On the positive semiaxis the function $\varphi$ shares the
same properties of the corresponding function associated with the
area. Consequently, above the mean the rate function $I$ has the same
features of the rate function of the area. To describe them we set
\begin{equation*}
  r:=\liminf_{s\uparrow+\infty}-\frac{1}{s^3}\ln\prob[S_1>s].
\end{equation*}
When $r$ is finite, we write $\xi:=\varphi(\sqrt{6r})$ for
brevity. When $r$ and $\xi$ are finite, we introduce the real number
\begin{equation*}
  \Lambda:=\Ex\big[e^{\xi S_1+\sqrt{6r}X_1}\big]\le 1
\end{equation*}
and the extended real number
\begin{equation*}
  \Xi:=\Ex\big[S_1^3e^{\xi S_1+rS_1^3}\big].
\end{equation*}
The following lemma characterizes the function $\varphi$ through the
values of $\xi$, $\Lambda$, and $\Xi$. The proof is presented in
Appendix \ref{proof:varphi_Aarea}.  Similarly to the area, the
asymptotic logarithmic moment generating function $g=-\varphi$ is not
steep when $0<r<+\infty$, $\xi>-\infty$, $\Lambda=1$, and
$\Xi<+\infty$. This case is beyond the scope of the G\"artner-Ellis
theory.
\begin{lemma}
  \label{lem:varphi_Aarea}
  The following conclusions hold:
  \begin{enumerate}[$(i)$]
   \item if $r=+\infty$, then $\varphi(k)$ solves the equation
     $\Ex[e^{\varphi(k) S_1+kX_1}]=1$ for every
     $k\in\Rl$.  The function $\varphi$ is analytic throughout $\Rl$
     with
     \begin{equation*}
       -\varphi'(k)=\frac{\Ex[X_1e^{\varphi(k)S_1+kX_1}]}{\Ex[S_1e^{\varphi(k)S_1+kX_1}]},
     \end{equation*}
     and the limits $\lim_{k\downarrow-\infty}-\varphi'(k)=0$ and
     $\lim_{k\uparrow+\infty}-\varphi'(k)=+\infty$ hold true;
\item if $0<r<+\infty$, then $\varphi(k)$ solves the equation
  $\Ex[e^{\varphi(k) S_1+kX_1}]=1$ for $k<\sqrt{6r}$, whereas
  $\varphi(k)=-\infty$ for $k>\sqrt{6r}$. The function $\varphi$ is
  analytic on $(-\infty,\sqrt{6r})$ with
\begin{equation*}
       -\varphi'(k)=\frac{\Ex[X_1e^{\varphi(k)S_1+kX_1}]}{\Ex[S_1e^{\varphi(k)S_1+kX_1}]}.
\end{equation*}
The derivative $\varphi'$ has the limits
$\lim_{k\downarrow-\infty}-\varphi'(k)=0$ and
\begin{equation*}
  \lim_{k\uparrow\sqrt{6r}}-\varphi'(k)=
  \begin{cases}
    \frac{\Ex[X_1e^{\xi S_1+\sqrt{6r}X_1}]}{\Ex[S_1e^{\xi S_1+\sqrt{6r}X_1}]} & \mbox{if $\xi>-\infty$, $\Lambda=1$, and $\Xi<+\infty$},\\
    +\infty & \mbox{otherwise};
    \end{cases}
\end{equation*}

\item if $r=0$ and $\ell=+\infty$, or $\ell<+\infty$ and $\Ex[e^{\ell
    S_1}]=+\infty$, or $\ell<+\infty$ and $\Ex[e^{\ell S_1}]<+\infty$
  and $\Ex[e^{(\ell+\nu_1|k|^{2/3})S_1+kX_1}]>1$ for all $k<0$, then
  $\varphi(k)$ solves the equation $\Ex[e^{\varphi(k)
      S_1+kX_1}]=1$ for $k\le 0$, whereas $\varphi(k)=-\infty$ for
  $k>0$. The function $\varphi$ is analytic on $(-\infty,0)$ with
\begin{equation*}
       -\varphi'(k)=\frac{\Ex[X_1e^{\varphi(k)S_1+kX_1}]}{\Ex[S_1e^{\varphi(k)S_1+kX_1}]}.
\end{equation*}
The derivative $\varphi'$ has the limits
$\lim_{k\downarrow-\infty}-\varphi'(k)=0$ and
\begin{equation*}
    \lim_{k\uparrow 0}-\varphi'(k)=\begin{cases}
    \mu=\sqrt{\frac{8}{9\pi}}\frac{\Ex[S_1^{3/2}]}{\Ex[S_1]} & \mbox{if }\Ex[S_1^{3/2}]<+\infty,\\
    +\infty & \mbox{if }\Ex[S_1^{3/2}]=+\infty.
    \end{cases}
\end{equation*}
  \end{enumerate}
\end{lemma}

We point out that part $(iii)$ of Lemma \ref{lem:varphi_Aarea} does
not consider the case $\ell<+\infty$, $\Ex[e^{\ell S_1}]<+\infty$, and
$\Ex[e^{(\ell+\nu_1|k|^{2/3})S_1+kX_1}]\le 1$ for some
$k<0$. Actually, it seems that this circumstance cannot occur. In
fact, numerical simulations suggest that the function that maps $s>0$
in
\begin{equation*}
\Ex\Big[e^{\nu_1s^{2/3}-s\int_0^1|B_{1,\tau}|\,d\tau}\Big]=\sum_{i=1}^\infty c_ie^{-(\nu_i-\nu_1)s^{2/3}}
\end{equation*}
has limit $1$ when $s$ approaches $0$ and is strictly increasing. It
follows that
\begin{equation*}
\Ex\Big[e^{(\ell+\nu_1|k|^{2/3})s+k\int_0^s|B_{1,\tau}|\,d\tau}\Big]\ge\Ex\Big[e^{\nu_1|k|^{2/3}s-|k|s^{3/2}\int_0^1|B_{1,\tau}|\,d\tau}\Big]>1
\end{equation*}
for all $k<0$ and $s>0$, which implies
$\Ex[e^{(\ell+\nu_1|k|^{2/3})S_1+kX_1}]>1$ for every
$k<0$. Unfortunately, we are not able to prove rigorously that
$\sum_{i=1}^\infty c_ie^{-(\nu_i-\nu_1)s^{2/3}}$ is strictly
increasing with respect to $s>0$, and we state this fact as a
conjecture.

Let us discuss briefly the graph of the rate function $I$ on the
region $[0,+\infty)$ containing its effective domain. If $r=+\infty$,
  or $0<r<+\infty$ and $\xi=-\infty$, or $0<r<+\infty$ and
  $\xi>-\infty$ and $\Lambda<1$, or $0<r<+\infty$ and $\xi>-\infty$
  and $\Lambda=1$ and $\Xi=+\infty$, then parts $(i)$ and $(ii)$ of
  Lemma \ref{lem:varphi_Aarea} show that the convex function
  $-\varphi$ is differentiable with continuously increasing derivative
  from $0$ at $k=-\infty$ to $+\infty$ at $k=\sqrt{6r}$. We agree that
  $\sqrt{6r}=+\infty$ when $r=+\infty$.  As a consequence, for each
  $w\in(0,+\infty)$ there exists $k\in(-\infty,\sqrt{6r})$ such that
  $w=-\varphi'(k)$ and $I(w)=wk+\varphi(k)$. The rate function $I$
  turns out to be analytic on $(0,+\infty)$ and good in all these
  cases. We have $I(0)=\lim_{w\downarrow 0}I(w)$ as $I$ is convex and
  lower semicontinuous (see \cite{Rockbook}, Corollary 7.5.1).

A singularity emerges when $0<r<+\infty$, $\xi>-\infty$, $\Lambda=1$,
and $\Xi<+\infty$, as for the area.  By setting
\begin{equation}
w_+:=\frac{\Ex[X_1e^{\xi S_1+\sqrt{6r}X_1}]}{\Ex[S_1e^{\xi S_1+\sqrt{6r}X_1}]}
\label{def:wplus}
\end{equation}
and by invoking part $(ii)$ of Lemma \ref{lem:varphi_area}, we have
that for each $w\in(0,w_+)$ there exists $k\in(-\infty,\sqrt{6r})$
such that $w=-\varphi'(k)$, in such a way that
$I(w)=wk+\varphi(k)$. The rate function is analytic on the open
interval $(0,w_+)$. As before, $I(0)=\lim_{w\downarrow 0}I(w)$. For
$w\ge w_+$, the function that maps $k\in(-\infty,\sqrt{6r}]$ in
$wk+\varphi(k)$ is non-decreasing, whereas $\varphi(k)=-\infty$ for
$k>\sqrt{6r}$. Thus,
$I(w)=\sup_{k\in\Rl}\{wk+\varphi(k)\}=w\sqrt{6r}+\varphi(\sqrt{6r})=w\sqrt{6r}+\xi$
for all $w\ge w_+$. The graph of $I$ is now characterized by an affine
stretch and a singularity at the point $w_+$. Example \ref{example:AA}
below discusses a situation where $0<r<+\infty$, $\xi>-\infty$,
$\Lambda=1$, and $\Xi<+\infty$.

Finally, we consider the case $r=0$, which differs from the
corresponding case of the area because the fluctuations of the
absolute area below the mean are characterized by a positive rate
function. This fact has been already found in \cite{LDP_res_2} when
studying Poissonian resetting.  Assuming that
$\Ex[e^{(\ell+\nu_1|k|^{2/3})S_1+kX_1}]>1$ for all $k<0$ when
$\ell<+\infty$ and $\Ex[e^{\ell S_1}]<+\infty$ and setting
$\mu:=+\infty$ if $\Ex[S_1^{3/2}]=+\infty$ for brevity, part $(iii)$
of Lemma \ref{lem:varphi_Aarea} states that the convex function
$-\varphi$ is differentiable on $(-\infty,0)$ with continuously
increasing derivative from $0$ at $k=-\infty$ to $\mu$ at $k=0$,
whereas $\varphi(k)=-\infty$ for $k>0$.  Thus, for each $w\in(0,\mu)$
there exists $k\in(-\infty,0)$ such that $w=-\varphi'(k)$ and
$I(w)=wk+\varphi(k)$. The rate function turns out to be analytic on
$(0,\mu)$ and $I(0)=\lim_{w\downarrow 0}I(w)$. We stress that these
conclusions are valid even for heavy-tailed waiting times, for which
$\ell=0$.  If $\mu<+\infty$, then $I(w)=0$ for $w\ge\mu$ since the
function that maps $k\in(-\infty,0]$ in $wk+\varphi(k)$ is
non-decreasing in this case.

\begin{example}
    \label{example:AA}
The rate function $I$ of the absolute area of the reset Brownian
motion with Poissonian resetting was studied in \cite{LDP_res_2}. Here
we discuss in detail the case
\begin{equation*}
  \prob[S_1>s]=e^{-rs^3}
\end{equation*}
for all $s\ge 0$ with some ``resetting rate'' $r>0$. We have
$\varphi(k)=-\infty$ for $k>\sqrt{6r}$ and
\begin{equation}
 \Ex\big[e^{\varphi(k) S_1+kX_1}\big]=\int_0^{+\infty}3rs^2e^{\varphi(k) s-rs^3}\Ex\Big[e^{ks^{\frac{3}{2}}\int_0^1|B_{1,\tau}|d\tau}\Big]ds=1
\label{example:absarea_1}
\end{equation}
for $k<\sqrt{6r}$. Unlike the area, the case $k=\sqrt{6r}$ is not
immediate.  Lemma \ref{lem:auxxx} tells us that for every $\zeta<0$
\begin{equation*}
\frac{3r}{\zeta^2}=3r\int_0^{+\infty} s^2e^{\zeta s}ds\le\Ex\big[e^{\zeta S_1+\sqrt{6r}X_1}\big]\le 3Lr\int_0^{+\infty} s^2e^{\zeta s}ds\le\frac{3Lr}{\zeta^2}
\end{equation*}
with some positive constant $L$.  Thus, the function that maps
$\zeta\in\Rl$ in $\Ex[e^{\zeta S_1+\sqrt{6r}X_1}]$ is finite and
continuous for $\zeta<0$ and goes to infinity when $\zeta$ goes to 0
from below. This gives $\xi:=\varphi(\sqrt{6k})>-\infty$ and
$\Lambda:=\Ex[e^{\xi S_1+\sqrt{6r}X_1}]=1$. Since $\xi<0$, we find
$\Xi:=\Ex[S_1^3e^{\xi S_1+rS_1^3}]=72\,r|\xi|^{-5}<+\infty$.
According to the above discuss, we finally realize that the rate
function takes the value
\begin{equation*}
  I(w)=\begin{cases}
  +\infty & \mbox{if }w\le0,\\
  wk+\varphi(k) & \mbox{if }w\in(0,w_+),\\
   w\sqrt{6r}+\xi & \mbox{if }w\ge w_+
  \end{cases}
\end{equation*}
with $k<\sqrt{6r}$ such that $w+\varphi'(k)=0$ and $w_+$ as in
(\ref{def:wplus}). This rate function exhibits one singularity and
one affine stretch.

Once again, we can make use of this example to discuss the limit of
zero resetting. As for the area, writing $\varphi_r$ and $I_r$ in
place of $\varphi$ and $I$ to stress the dependence on $r$, the change
of variable $s\mapsto r^{-\frac{1}{3}}s$ in (\ref{example:absarea_1})
shows that $\varphi_r(k)=r^{\frac{1}{3}}\varphi_1(r^{-\frac{1}{2}}k)$
for all $k$ and $I_r(w)=r^{\frac{1}{3}}I_1(r^{\frac{1}{6}}w)$ for all
$w$.  The analysis of $I_1$ in a neighborhood of the origin allows to
conclude that there exists a positive constant $C$ such that
\begin{equation*}
  \bigg|I_1(w)-\frac{4\nu_1^3}{27w^2}\bigg|\le C(1+w)
\end{equation*}
for any $w>0$. Thus, for all $r>0$ and $w>0$ we obtain
\begin{equation*}
  \bigg|I_r(w)-\frac{4\nu_1^3}{27w^2}\bigg|\le C\big(r^{\frac{1}{3}}+r^{\frac{1}{2}}w\big).
\end{equation*}
As $r$ goes to zero, $I_r(w)$ approaches a non-trivial limit with
speed $r^{\frac{1}{3}}$.  The result is consistent with the known
asymptotics of the absolute area without resetting \cite{Janson}:
\begin{equation*}
  \lim_{t\uparrow+\infty}\frac{1}{t}\ln\prob\bigg[\frac{1}{t}\int_0^t|B_{1,\tau}|d\tau\le w\bigg]
  =\lim_{t\uparrow+\infty}\frac{1}{t}\ln\prob\bigg[\int_0^1|B_{1,\tau}|d\tau\le \frac{w}{\sqrt{t}}\bigg]=
  -\frac{4\nu_1^3}{27w^2}.
\end{equation*}
\end{example}

\section{Conclusions}
\label{sec:conclusions}

The mathematical theory of renewal-reward processes boasts a long
tradition.  In this paper we have taken advantage of this theory to
investigate the fluctuations, both normal and large, of additive
functionals associated with a stochastic process under a general
non-Poissonian resetting mechanism. While providing the law of large
numbers, the central limit theorem, and a large deviation principle,
we have demonstrated that a suitable resetting protocol can always
bring about a large deviation principle with good rate function. In
other words, it is always possible to find a resetting mechanism that
makes the probability of a large fluctuation of these functionals
exponentially small in time. This confirms the ability of resetting to
confine the process around the initial position.

We have used our general results to characterize the fluctuations of
the positive occupation time, the area, and the absolute area of the
reset Brownian motion under any resetting mechanism.  Our general
results constitute a solid background for researchers that aim to
analyze the fluctuations of other additive functionals and other
stochastic processes than the Brownian motion. We stress that there
are no restrictions to the resetting protocol that can be
implemented.

Our study of the positive occupation time, the area, and the absolute
area of the reset Brownian motion has shown that a rich phenomenology
accounting for dynamical phase transitions emerges beyond Poissonian
resetting.  This paper opens to investigation of the mechanisms that
originate such dynamical phase transitions. Although we leave to
future research the task of providing precise results, we suggest that
these phase transitions are due to a phenomenon of condensation of
fluctuations similarly to sums of independent and identically
distributed random variables
\cite{Condensation_1,Condensation_2,Condensation_3,Condensation_4}.
After all, an additive functional $F_t$ of the reset Brownian motion
is close to the associated cumulative reward $W_t$, which is a random
sum of i.i.d.\ rewards. According to this interpretation, a
fluctuation of $F_t/t$ that crosses a singular point of the rate
function would be realized by the combination of two different
mechanisms: many small deviations of the rewards all in the same
direction for a partial excursion up to the singular point, and a big
jump of only one of the rewards for the remaining part of the
fluctuation. We point out that even if in some cases the fluctuations
of rewards are solely determined by the waiting times, such as for the
positive occupation time of the reset Brownian motion, in general they
are the result of the interplay between the fluctuations of waiting
times and the fluctuations of the reference process, such as for the
area and the absolute area of the reset Brownian motion.


\appendix

\section{Proof of Proposition \ref{prop:conv}}
\label{proof:conv}

The proof of the proposition requires a little of convex analysis and
we refer to \cite{Rockbook} for the needed results. It $\psi$ is a
proper convex function on $\Rl^d$ taking extended real values, we
denote by $\dom\psi:=\{x\in\Rl^d:\psi(x)<+\infty\}$ the effective
domain of $\psi$. The convex conjugate of $\psi$ is the function
$\psi^\star$ that maps $x^\star\in\Rl^d$ in
$\psi^\star(x^\star):=\sup_{x\in\Rl^d}\{x^\star\cdot x-\psi(x)\}$. The
function $\psi^\star$ is convex and lower semicontinuous and
$\psi^{\star\star}$ turns out to be the lower semicontinuous hull of
$\psi$ (see \cite{Rockbook}, Theorem 12.2), i.e.\ the largest lower
semicontinuous function not greater than $\psi$. Consequently,
$\psi^{\star\star}=\psi$ when $\psi$ is lower semicontinuous.

If $x^\star$ is a subgradients of $\psi$ at $x$, then
$\psi^\star(x^\star)=x^\star\cdot x-\psi(x)$ (see \cite{Rockbook},
Theorem 23.5). If $\psi$ is lower semicontinuous, then for each
$x^\star$ in the relative interior of $\dom\psi^\star$ there exists
$x\in\dom\psi$ such that $x^\star\in\partial\psi(x)$ (see
\cite{Rockbook}, Theorems 23.4 and 23.5). The relative interior of
$\dom\psi^\star$ is the interior that results when $\dom\psi^\star$ is
regarded as a subset of its affine hull.

\subsection{Proof of part $\boldsymbol{(i)}$}

If $k\in\Rl^d$ is such that $\varphi(k)>-\infty$, then
$\Ex[e^{\varphi(k)S_1+k\,\cdot X_1}]\le 1$. Indeed, by definition of
supremum, there exists a sequence $\{\zeta_i\}_{i\ge 1}$ of real
numbers that satisfies $\Ex[e^{\zeta_i S_1+k\,\cdot X_1}]\le 1$ for
all $i$ and $\lim_{i\uparrow\infty}\zeta_i=\varphi(k)$. Fatou's lemma
gives $\Ex[e^{\varphi(k) S_1+k\cdot
    X_1}]\le\liminf_{i\uparrow\infty}\Ex[e^{\zeta_i S_1+k\cdot
    X_1}]\le 1$. We shall use this property of $\varphi$ many times
along the paper.

Let us prove concavity of $\varphi$ by verifying that
$\varphi(ah+bk)\ge a\varphi(h)+b\varphi(k)$ for any given $h\in\Rl^d$,
$k\in\Rl^d$, $a\in(0,1)$, and $b\in(0,1)$ such that
$\varphi(h)>-\infty$, $\varphi(k)>-\infty$, and $a+b=1$. Since
$\Ex[e^{\varphi(h)S_1+h\,\cdot X_1}]\le 1$ and
$\Ex[e^{\varphi(k)S_1+k\,\cdot X_1}]\le 1$, H\"older's inequality
gives
\begin{equation*}
\Ex\Big[e^{\{a\varphi(h)+b\varphi(k)\}S_1+\{ah+bk\}\cdot X_1}\Big]\le\Big\{\Ex\big[e^{\varphi(h)S_1+h\,\cdot X_1}\big]\Big\}^a\Big\{\Ex\big[e^{\varphi(k)S_1+k\,\cdot X_1}\big]\Big\}^b\le 1.
\end{equation*}
We deduce $\varphi(ah+bk)\ge a\varphi(h)+b\varphi(k)$ from here by the
definition of $\varphi(ah+bk)$. Thus, $-\varphi$ is convex with
effective domain $\dom(-\varphi)=\{k\in\Rl^d:\varphi(k)>-\infty\}$.

The function $\varphi$ is upper semicontinuous if the level set
$\mathcal{L}_a:=\{k\in\Rl^d:\varphi(k)\ge a\}$ is closed for any
$a\in\Rl$. Given $a\in\Rl$, let us show that if $\{k_i\}_{i\ge
  1}\subseteq\mathcal{L}_a$ is a sequence converging to $k$, then
$k\in\mathcal{L}_a$. This demonstrates that $\mathcal{L}_a$ is
closed. As $\varphi(k_i)\ge a>-\infty$, we have $\Ex[e^{aS_1+k_i\cdot
    X_1}]\le\Ex[e^{\varphi(k_i)S_1+k_i\cdot X_1}]\le 1$ for all
$i$. Then, Fatou's lemma yields $\Ex[e^{a S_1+k\cdot
    X_1}]\le\liminf_{i\uparrow\infty}\Ex[e^{a S_1+k_i\cdot X_1}]\le
1$. This bound shows that $\varphi(k)\ge a$.

\subsection{Proof of part $\boldsymbol{(ii)}$}

Basically, the proof of part $(ii)$ amounts to verify that for all
$\beta\ge 0$ and $w\in\Rl^d$
\begin{equation}
\Upsilon(\beta,w)=\sup_{k\in\doms(-\varphi)}\big\{w\cdot k+\beta\varphi(k)\big\}.
\label{conv_0}
\end{equation}
In fact, if (\ref{conv_0}) is valid and $\ell=+\infty$, then part
$(ii)$ follows by setting $\beta=1$ in (\ref{conv_0}). If
(\ref{conv_0}) is valid and $\ell<+\infty$, then part $(ii)$ follows
by applying Sion's minimax theorem to the function that, for a given
$w\in\Rl^d$, maps $(\beta,k)\in[0,1]\times\dom(-\varphi)$ in $w\cdot
k+\beta\varphi(k)+(1-\beta)\ell$.  This function is convex and
continuous with respect to $\beta$ for each fixed $k\in\dom(-\varphi)$
and concave and upper semicontinuous with respect to $k$ for each
fixed $\beta\in[0,1]$.  Thus, Sion's minimax theorem allows us to
exchange an infimum with respect to $\beta$ and a supremum with
respect to $k$ by stating that
\begin{align}
  \nonumber
  I(w):=\inf_{\beta\in[0,1]}\big\{\Upsilon(\beta,w)+(1-\beta)\ell\big\}&=
  \adjustlimits\inf_{\beta\in[0,1]}\sup_{k\in\doms(-\varphi)}\big\{w\cdot k+\beta\varphi(k)+(1-\beta)\ell\big\}\\
  \nonumber
  &=\adjustlimits\sup_{k\in\doms(-\varphi)}\inf_{\beta\in[0,1]}\big\{w\cdot k+\beta\varphi(k)+(1-\beta)\ell\big\}\\
  \nonumber
  &=\sup_{k\in\doms(-\varphi)}\big\{w\cdot k+\varphi(k)\wedge \ell\big\}.
\end{align}
We stress that Sion's minimax theorem applies because $\beta$ lies in
a compact set.

Let us move to the proof of identity (\ref{conv_0}). The lower
bound
\begin{equation}
\Upsilon(\beta,w)\ge\sup_{k\in\doms(-\varphi)}\big\{w\cdot k+\beta\varphi(k)\big\}
\label{conv_l}
\end{equation}
is immediate for all $\beta\ge 0$ and $w\in\Rl^d$. In fact, if
$(\zeta,k)\in\Rl\times\Rl^d$ is such that $\Ex[e^{\zeta S_1+k\,\cdot
    X_1}]\le 1$, then $J(s,w)\ge s\zeta+w\cdot k$ by definition for
every $(s,w)\in\Rl\times\Rl^d$. In particular, we have $J(s,w)\ge
w\cdot k+s\varphi(k)$ if $k\in\dom(-\varphi)$.  It follows that
$\Upsilon(\beta,w)\ge w\cdot k+\beta\varphi(k)$ for all $\beta\ge 0$,
$w\in\Rl^d$, and $k\in\dom(-\varphi)$.

The upper bound
\begin{equation}
\Upsilon(\beta,w)\le\sup_{k\in\doms(-\varphi)}\big\{w\cdot k+\beta\varphi(k)\big\}
\label{conv_1}
\end{equation}
is much more involved, and we need to split the proof in several
steps. We solve the case $\beta=0$ first, then we address the case
$\beta>0$ by means of a truncation argument.  The function that maps
$(\zeta,k)\in\Rl\times\Rl^d$ in $\ln\Ex[e^{\zeta S_1+k\,\cdot X_1}]$
is proper convex, so that there exist $(\beta_o,w_o)\in\Rl\times\Rl^d$
and $c\in\Rl$ such that $\ln\Ex[e^{\zeta S_1+k\,\cdot X_1}]\ge
\beta_o\zeta+w_o\cdot k-c$ for all $(\zeta,k)\in\Rl\times\Rl^d$ (see
\cite{Rockbook}, Corollary 12.1.2). Let
$\mathcal{D}:=\{(\zeta,k)\in\Rl\times\Rl^d:\Ex[e^{\zeta S_1+k\,\cdot
    X_1}]<+\infty\}$ be its effective domain. For every $w\in\Rl^d$,
$\delta>0$, and $\gamma_o\in(0,\delta)$ such that
$\gamma_o|\beta_o|<\delta$ and $\gamma_o\|w_o\|<\delta$ we have
\begin{align}
  \nonumber
  \adjustlimits\inf_{s\in (-\delta,\delta)}\inf_{v\in \Delta_{w,\delta}}\inf_{\gamma>0}
  \bigg\{\gamma J\bigg(\frac{s}{\gamma},\frac{v}{\gamma}\bigg)\bigg\}&\le
  \inf_{\gamma>0}\bigg\{\gamma J\bigg(\frac{\gamma_o\beta_o}{\gamma},\frac{w+\gamma_ow_o}{\gamma}\bigg)\bigg\}\\
\nonumber
  &\le\gamma_o J\bigg(\beta_o,\frac{w}{\gamma_o}+w_o\bigg)\\
  \nonumber
  &=\sup_{(\zeta,k)\in\mathcal{D}}\Big\{\gamma_o\beta_o\zeta+(w+\gamma_ow_o)\cdot k-
  \gamma_o\ln\Ex\big[e^{\zeta S_1+k\,\cdot X_1)}\big]\Big\}\\
  \nonumber
  &\le\sup_{(\zeta,k)\in\mathcal{D}}\big\{w\cdot k\big\}+\gamma_o c\le\sup_{(\zeta_o,k)\in\mathcal{D}}\big\{w\cdot k\big\}+\delta c.
\end{align}
On the other hand, if $(\zeta_o,k)\in\mathcal{D}$, then
$k\in\dom(-\varphi)$. In fact, if $\Ex[e^{\zeta_o S_1+k\,\cdot
    X_1}]<+\infty$, then $\lim_{\zeta\downarrow-\infty}\Ex[e^{\zeta
    S_1+k\,\cdot X_1}]=0$ by the dominated convergence theorem, so
that the set $\{\zeta\in\Rl:\Ex[e^{\zeta S_1+k\,\cdot X_1}]\le 1\}$ is
non-empty. This way, we find for all $w\in\Rl^d$
\begin{equation*}
  \Upsilon(0,w):=\adjustlimits\lim_{\delta\downarrow 0}\inf_{s\in (-\delta,\delta)}\inf_{v\in \Delta_{w,\delta}}\inf_{\gamma>0}
  \bigg\{\gamma J\bigg(\frac{s}{\gamma},\frac{v}{\gamma}\bigg)\bigg\}\le \sup_{(\zeta_o,k)\in\mathcal{D}}\{w\cdot k\}\le\sup_{k\in\doms(-\varphi)}\{w\cdot k\}.
\end{equation*}
This bound is (\ref{conv_1}) for $\beta=0$.

The proof of (\ref{conv_1}) for $\beta>0$ relies on a truncation
argument that allows us to deal with differentiable approximants of
$\varphi$.  For each integer $i\ge 1$, the approximant of $\varphi$ we
consider is the real function $\varphi_i$ that associates any
$k\in\Rl^d$ with the unique solution $\varphi_i(k)$ of the equation
\begin{equation*}
\Ex\Big[e^{\varphi_i(k)S_1+k\,\cdot X_1}\mathds{1}_{\{S_1\vee\|X_1\|\le i\}}\Big]=1.
\end{equation*}
The function $\varphi_i$ is finite throughout $\Rl^d$ and analytic by
the analytic implicit function theorem. The gradient
$\nabla\varphi_i(k)$ and the Hessian matrix $H_i(k)$ of $\varphi_i$ at
$k$ satisfies for every $k\in\Rl^d$ and $u\in\Rl^d$ the relationships
\begin{equation}
  \nabla\varphi_i(k)=-\frac{\Ex\big[X_1e^{\varphi_i(k)S_1+k\,\cdot X_1}\mathds{1}_{\{S_1\vee\|X_1\|\le i\}}\big]}
               {\Ex\big[S_1e^{\varphi_i(k)S_1+k\,\cdot X_1}\mathds{1}_{\{S_1\vee\|X_1\|\le i\}}\big]}
\label{grad}
\end{equation}
and
\begin{equation*}
u\,\cdot H_i(k)\,u=-\frac{\Ex\big[\{u\cdot \nabla\varphi_i(k)S_1+u\,\cdot X_1\}^2e^{\varphi_i(k)S_1+k\,\cdot X_1}\mathds{1}_{\{S_1\vee\|X_1\|\le i\}}\big]}
  {\Ex\big[S_1e^{\varphi_i(k)S_1+k\,\cdot X_1}\mathds{1}_{\{S_1\vee\|X_1\|\le i\}}\big]}\le 0.
\end{equation*}
The latter shows that $\varphi_i$ is concave.  

Clearly, the sequence $\{\varphi_i(k)\}_{i\ge 1}$ is non-increasing for
any $k\in\Rl^d$. Let $\varphi_\infty$ be the function that maps $k$ in
$\varphi_\infty(k):=\lim_{i\uparrow\infty}\varphi_i(k)=\inf_{i\ge
  1}\{\varphi_i(k)\}$.  We claim that $\varphi_\infty=\varphi$ as
desired. In fact, the bound $\varphi_\infty(k)\ge \varphi(k)$ is
non-trivial only for $k\in\dom(-\varphi)$, and if $k\in\dom(-\varphi)$,
then it follows from the fact that $\Ex[e^{\varphi(k)S_1+k\,\cdot
    X_1}\mathds{1}_{\{S_1\vee\|X_1\|\le i\}}]\le 1$ for all $i$ by
construction of $\varphi$. The bound $\varphi_\infty(k)\le \varphi(k)$
is non-trivial only when $\varphi_\infty(k)>-\infty$ and is due to
Fatou's lemma, which shows that $\Ex[e^{\varphi_\infty(k)S_1+k\,\cdot
    X_1}]\le\liminf_{i\uparrow\infty}\Ex[e^{\varphi_i(k)S_1+k\,\cdot
    X_1}\mathds{1}_{\{S_1\vee\|X_1\|\le i\}}]=1$. The latter
inequality entails that $\varphi_\infty(k)\le \varphi(k)$ by
definition of $\varphi(k)$.

We shall verify that for all $\beta>0$ and $w\in\Rl^d$
\begin{equation}
\Upsilon(\beta,w)\le\adjustlimits\inf_{i\ge 1}\sup_{k\in\Rl^d}\big\{w\cdot k+\beta\varphi_i(k)\big\}=\inf_{i\ge 1}\,(-\beta\varphi_i)^\star(w).
\label{conv_3}
\end{equation}
This bound gives the upper bound (\ref{conv_1}) as follows. Fix
$\beta_o>0$ and consider the function $\vartheta$ that maps
$w\in\Rl^d$ in $\vartheta(w):=\inf_{i\ge
  1}\,(-\beta_o\varphi_i)^\star(w)$. The limit
$\vartheta(w)=\lim_{i\uparrow\infty}(-\beta_o\varphi_i)^\star(w)$
holds true because $\{(-\beta_o\varphi_i)^\star(w)\}_{i\ge 1}$ is a
non-increasing sequence, and the function $\vartheta$ turns out to be
convex as it is the pointwise limit of a sequence of convex
functions. Since $\vartheta^{\star\star}$ is the largest lower
semicontinuous function not greater than $\vartheta$ and since
$\Upsilon(\beta_o,\cdot\,)$ is lower semicontinuous by construction
and $\Upsilon(\beta_o,w)\le\vartheta(w)$ for all $w\in\Rl^d$ by
(\ref{conv_3}), we have
$\Upsilon(\beta_o,w)\le\vartheta^{\star\star}(w)$ for every
$w\in\Rl^d$. Let us show that
$\vartheta^{\star\star}(w)\le(-\beta_o\varphi)^\star(w)$ for all
$w\in\Rl^d$. The fact that
$\vartheta(w)\le(-\beta_o\varphi_i)^\star(w)$ for all $w\in\Rl^d$ and
$i\ge 1$ yields $\vartheta^\star(k)\ge
(-\beta_o\varphi_i)^{\star\star}(k)=-\beta_o\varphi_i(k)$ for all
$k\in\Rl^d$. We have exploited the identity
$(-\beta_o\varphi_i)^{\star\star}=-\beta_o\varphi_i$, which is valid
because $-\beta_o\varphi_i$ is convex and continuous. By sending $i$
to infinity we find $\vartheta^\star(k)\ge-\beta_o\varphi(k)$ for all
$k\in\Rl^d$. Thus,
$\vartheta^{\star\star}(w)\le(-\beta_o\varphi)^\star(w)$ for each
$w\in\Rl^d$.

Let us move to demonstrate (\ref{conv_3}).  To this aim, for each
integer $i\ge 1$ we introduce ``truncated'' functions $J_i$ and
$\Upsilon_i$, which are defined for every $(s,w)\in\Rl\times\Rl^d$ and
$(\beta,w)\in\Rl\times\Rl^d$ by
\begin{equation*}
J_i(s,w):=\sup_{(\zeta,k)\in\Rl\times\Rl^d}\bigg\{s\zeta+w\cdot k-\ln\Ex\Big[e^{\zeta S_1+k\,\cdot X_1}\mathds{1}_{\{S_1\vee\|X_1\|\le i\}}\Big]\bigg\}
\end{equation*}
and 
\begin{equation*}
  \Upsilon_i(\beta,w):=\adjustlimits\lim_{\delta\downarrow 0}\inf_{s\in (\beta-\delta,\beta+\delta)}\inf_{v\in \Delta_{w,\delta}}\inf_{\gamma>0}
  \bigg\{\gamma J_i\bigg(\frac{s}{\gamma},\frac{v}{\gamma}\bigg)\bigg\}.
\end{equation*}
We have $J(s,w)\le J_i(s,w)$ and $\Upsilon(\beta,w)\le
\Upsilon_i(\beta,w)$ as a consequence. Bound (\ref{conv_3}) is due to
the identity
\begin{equation}
\Upsilon_i(\beta,w)=\sup_{k\in\Rl^d}\big\{w\cdot k+\beta\varphi_i(k)\big\}
\label{conv_4}
\end{equation}
valid for all $\beta>0$, $w\in\Rl^d$, and $i\ge 1$, which we are going
to verify by exploiting the differentiability of $\varphi_i$.

By repeating the arguments we used to deduce the lower bound
(\ref{conv_l}), we realize that
$\Upsilon_i(\beta,w)\ge\sup_{k\in\Rl^d}\{w\cdot k+\beta\varphi_i(k)\}$
for all $\beta>0$, $w\in\Rl^d$, and $i\ge 1$. The main task to achieve
(\ref{conv_4}) is then to prove that
$\Upsilon_i(\beta,w)\le\sup_{k\in\Rl^d}\{w\cdot
k+\beta\varphi_i(k)\}$. The differentiability of $\varphi_i$ comes
into play at this point. Fix $\beta_o>0$ and $i\ge 1$. The problem is
to demonstrate that $\Upsilon_i(\beta_o,w)\le
(-\beta_o\varphi_i)^\star(w)$ for all
$w\in\dom(-\beta_o\varphi_i)^\star$. Let $D_o$ be the relative
interior of $\dom(-\beta_o\varphi_i)^\star$. We prove at first that
$\Upsilon_i(\beta_o,w)\le (-\beta_o\varphi_i)^\star(w)$ for all $w\in
D_o$, and then we extend this bound from $D_o$ to
$\dom(-\beta_o\varphi_i)^\star$. We know from the general properties
of convex functions mentioned at the beginning that if $w\in D_o$,
then there exists $k_o\in\Rl^d$ such that
$w=-\beta\nabla\varphi_i(k_o)$ and
$(-\beta_o\varphi_i)^\star(w)=w\cdot k_o+\beta_o\varphi_i(k_o)$. Set
$\zeta_o:=\varphi_i(k_o)$ and
\begin{equation*}
\gamma_o:=\frac{\beta_o}{\Ex[S_1e^{\zeta_oS_1+k_o\cdot X_1}\mathds{1}_{\{S_1\vee\|X_1\|\le i\}}]}.
\end{equation*}
By definition we have
\begin{equation*}
  \Upsilon_i(\beta_o,w)\le\inf_{\gamma>0}\bigg\{\gamma J_i\bigg(\frac{\beta_o}{\gamma},\frac{w}{\gamma}\bigg)\bigg\}
  \le\gamma_o J_i\bigg(\frac{\beta_o}{\gamma_o},\frac{w}{\gamma_o}\bigg).
\end{equation*}
On the other hand, $J_i$ is the convex conjugate of the function that
associates the pair $(\zeta,k)\in\Rl\times\Rl^d$ with $\ln\Ex[e^{\zeta
    S_1+k\,\cdot X_1}\mathds{1}_{\{S_1\vee\|X_1\|\le i\}}]$, and the
gradient of this function at $(\zeta_o,k_o)$ is exactly
$(\beta_o/\gamma_o,w/\gamma_o)$ as one can easily verify by appealing
to (\ref{grad}). Then, bearing in mind that $\Ex[e^{\zeta_o
    S_1+k_o\,\cdot X_1}\mathds{1}_{\{S_1\vee\|X_1\|\le i\}}]=1$, we
get
\begin{align}
  \nonumber
\Upsilon_i(\beta_o,w)&\le\gamma_o J_i\bigg(\frac{\beta_o}{\gamma_o},\frac{w}{\gamma_o}\bigg)=\beta_o\zeta_o+w\cdot k_o-\gamma_o\ln\Ex\Big[e^{\zeta_o S_1+k_o\,\cdot
    X_1}\mathds{1}_{\{S_1\vee\|X_1\|\le i\}}\Big]\\
\nonumber
&=w\cdot k_o+\beta_o\varphi_i(k_o)=(-\beta_o\varphi_i)^\star(w).
\end{align}
Thus, $\Upsilon_i(\beta_o,w)\le(-\beta_o\varphi_i)^\star(w)$ for any
$w\in D_o$. In order to extend this bound from $D_o$ to
$\dom(-\beta_o\varphi_i)^\star$, pick $w_o\in D_o$ and any
$w\in\dom(-\beta_o\varphi_i)^\star$. The half-open line segment
$\{\lambda w+(1-\lambda)w_o:\lambda\in[0,1)\}$ lies in $D_o$ (see
  \cite{Rockbook}, Theorem 6.1), so that we can state
  \begin{equation}
  \Upsilon_i\big(\beta_o,\lambda w+(1-\lambda)w_o\big)\le (-\beta_o\varphi_i)^\star\big(\lambda w+(1-\lambda)w_o\big)
\label{D_to_dom}
  \end{equation}
for all $\lambda\in[0,1)$. We have $\lim_{\lambda\uparrow
    1}(-\beta_o\varphi_i)^\star(\lambda
  w+(1-\lambda)w_o)=(-\beta_o\varphi_i)^\star(w)$ as
  $(-\beta_o\varphi_i)^\star$ is convex and lower semicontinuous (see
  \cite{Rockbook}, Corollary 7.5.1). At the same time, we have
  $\liminf_{\lambda\uparrow 1}\Upsilon_i(\beta_o,\lambda
  w+(1-\lambda)w_o)\ge \Upsilon_i(\beta_o,w)$ since
  $\Upsilon_i(\beta_o,\cdot\,)$ is lower semicontinuous by
  construction. Thus, we find $\Upsilon_i(\beta_o,w)\le
  (-\beta_o\varphi_i)^\star(w)$ by sending $\lambda$ to one from below
  in (\ref{D_to_dom}). Resuming, we have shown that
  $\Upsilon_i(\beta_o,w)\le(-\beta_o\varphi_i)^\star(w)$ for all
  $w\in\dom(-\beta_o\varphi_i)^\star$, and hence for all
  $w\in\Rl^d$. Since $\beta_o>0$ and $i\ge 1$ are arbitrary, the proof
  of (\ref{conv_4}) is concluded.

\section{Proof of Proposition \ref{th:LDP_preliminar}}
\label{proof:LDP_preliminar}

Let $\varpi$ be the function that maps $k\in\Rl$ in $\varpi(k)$ and
set $\vartheta:=\varphi\wedge\varpi$ for brevity.

\subsection{Proof of part $\boldsymbol{(i)}$}

Fix $k\in\Rl$.  Let us demonstrate the bound
\begin{equation}
  \limsup_{t\uparrow+\infty}\frac{1}{t}\ln\Ex\big[e^{kF_t}\big]\le -\vartheta(k).
\label{boundsup_cumul}
\end{equation}
Assume that $\varphi(k)>-\infty$ and $\varpi(k)>-\infty$, otherwise
$\vartheta(k)=\varphi(k)\wedge\varpi(k)=-\infty$ and the bound is
trivial. Pick numbers $\epsilon>0$ and $\varrho<\varpi(k)$ and denote
by $P:=\prob[S_1\in\cdot\,]$ the probability measure induced by $S_1$
over $[0,+\infty)$. By definition of $\varpi(k)$, there exists a
  positive constant $C$ such that for all $t\ge 0$
\begin{equation*}
  \mathcal{E}_t(k)\prob[S_1>t]=\int_{(t,+\infty)}\mathcal{E}_t(k)P(ds)\le Ce^{\epsilon t-\varrho t}
\end{equation*}
with $\mathcal{E}_t(k):=\Ex[e^{k\int_0^tf(B_{1,\tau})\,d\tau}]$.
Moreover, we have
\begin{equation*}
  a:=\Ex\big[e^{\varphi(k)S_1-\epsilon S_1}\mathcal{E}_{S_1}(k)\big]=\Ex\big[e^{\varphi(k)S_1-\epsilon S_1+kX_1}\big]<1
\end{equation*}
by the independence between $S_1$ and $\{B_{1,t}\}_{t\ge 0}$ and the
definition of $\varphi(k)$.  Then, invoking again the independence
between the waiting times and the processes $\{B_{1,t}\}_{t\ge
  0},\{B_{2,t}\}_{t\ge 0},\ldots$ we can write for every $t$
\begin{align}
  \nonumber
  \Ex\big[e^{kF_t}\big]&=\Ex\bigg[e^{\sum_{i=1}^{N_t}k\int_0^{S_i}f(B_{i,\tau})\,d\tau+k\int_0^{t-T_{N_t}}f(B_{N_t+1,\tau})\,d\tau}\bigg]\\
  \nonumber
  &=\sum_{n=0}^\infty\Ex\bigg[\prod_{i=1}^n \mathcal{E}_{S_i}(k)\mathds{1}_{\{T_n\le t\}}\mathcal{E}_{t-T_n}(k)\mathds{1}_{\{S_{n+1}>t-T_n\}}\bigg]\\
  \nonumber
  &=\sum_{n=0}^\infty\Ex\Bigg[\prod_{i=1}^n \mathcal{E}_{S_i}(k)\mathds{1}_{\{T_n\le t\}}\int_{(t-T_n,+\infty)}\mathcal{E}_{t-T_n}(k)P(ds)\Bigg]\\
  &\le C\sum_{n=0}^\infty\Ex\Bigg[\prod_{i=1}^n \mathcal{E}_{S_i}(k)\mathds{1}_{\{T_n\le t\}}e^{\epsilon(t-T_n)-\varrho(t-T_n)}\Bigg].
  \label{moment_bound_0}
\end{align}
Hereafter we make the usual convention that an empty sum is $0$ and an
empty product is $1$.  If $\varphi(k)\le\varrho$, then bound
(\ref{moment_bound_0}) gives
\begin{align}
  \nonumber
  \Ex\big[e^{kF_t}\big]&\le C\sum_{n=0}^\infty\Ex\Bigg[\prod_{i=1}^n \mathcal{E}_{S_i}(k)\mathds{1}_{\{T_n\le t\}}e^{\epsilon(t-T_n)-\varphi(k)(t-T_n)}\Bigg]\\
  \nonumber
  &\le Ce^{\epsilon t-\varphi(k)t}\sum_{n=0}^\infty\Ex\Bigg[\prod_{i=1}^n e^{\varphi(k)S_i-\epsilon S_i}\mathcal{E}_{S_i}(k)\Bigg]\\
  \nonumber
  &=Ce^{\epsilon t-\varphi(k)t}\sum_{n=0}^\infty a^n=\frac{Ce^{\epsilon t-\varphi(k)t}}{1-a}.
\end{align}
If $\varphi(k)>\varrho$, then bound (\ref{moment_bound_0}) gives
\begin{align}
  \nonumber
  \Ex\big[e^{kF_t}\big]&\le Ce^{\epsilon t-\varrho t}\sum_{n=0}^\infty\Ex\Bigg[\prod_{i=1}^n e^{\varphi(k)S_i-\epsilon S_i}\mathcal{E}_{S_i}(k)\Bigg]
  =\frac{Ce^{\epsilon t-\varrho t}}{1-a}.
\end{align}
In conclusion, we find
\begin{equation*}
  \limsup_{t\uparrow+\infty}\frac{1}{t}\ln\Ex\big[e^{kF_t}\big]\le -\varphi(k)\wedge\varrho+\epsilon,
\end{equation*}
which implies (\ref{boundsup_cumul}) thanks to the arbitrariness of
$\epsilon>0$ and $\varrho<\varpi(k)$.

\subsection{Proof of part $\boldsymbol{(ii)}$}

We are going to use part $(i)$ of the proposition to demonstrate the
upper large deviation bound for closed sets with rate function $I$
equal to the convex conjugate of $-\vartheta$.  To begin with, we
observe that the function $\varpi$ is concave, so that $\vartheta$ is
concave. Similarly to part $(i)$ of Proposition \ref{prop:conv}, this
fact follows from the H\"older's inequality and we do not repeat the
proof. Then, we notice that $\inf_{w\in\Rl}\{I(w)\}=0$. Indeed, since
$\varphi(0)=0$ and $\varpi(0)\ge 0$ by definition, we have $I(w)\ge
\vartheta(0)=\varphi(0)\wedge\varpi(0)=0$ for all $w\in\Rl$. At the
same time, denoting by $(-\varphi)^\star$ the convex conjugate of
$-\varphi$ and by $(-\varphi)^{\star\star}$ the convex conjugate of
$(-\varphi)^\star$ we find
\begin{align}
  \nonumber
  \inf_{w\in\Rl}\{I(w)\}=\adjustlimits\inf_{w\in\Rl}\sup_{k\in\Rl}\big\{wk+\vartheta(k)\big\}&\le
  \adjustlimits\inf_{w\in\Rl}\sup_{k\in\Rl}\big\{wk+\varphi(k)\big\}\\
  \nonumber
  &=\inf_{w\in\Rl}\big\{(-\varphi)^\star(w)\big\}=-(-\varphi)^{\star\star}(0).
\end{align}
This bound is $\inf_{w\in\Rl}\{I(w)\}\le 0$ since
$(-\varphi)^{\star\star}(0)=\varphi(0)=0$ as $-\varphi$ is convex and
lower semicontinuous.

We are ready to prove the upper large deviation bound for the additive
functional $F_t$.  Fix a closed set $K$ such that $\inf_{w\in
  K}\{I(w)\}>0$, otherwise the upper large deviation bound trivially
holds.  If $\vartheta(k)=-\infty$ for $k\ne 0$, then $I(w)=0$ for all
$w\in\Rl$ and $\inf_{w\in K}\{I(w)\}=0$. Thus, we are actually
considering a model where $\vartheta(k)$ is finite for some $k\ne
0$. We need to distinguish three different situations:
$\vartheta(k)=-\infty$ for all $k<0$, $\vartheta(k)=-\infty$ for all
$k>0$, and $\vartheta(k)>-\infty$ for $k$ in an open neighborhood of
the origin.

If $\vartheta(k)=-\infty$ for all $k<0$, then $I(w)=\sup_{k\ge 0}\{wk+\vartheta(k)\}$
is non-decreasing with respect to $w$ and
$\lim_{w\downarrow-\infty}I(w)=0$ since
$\inf_{w\in\Rl}\{I(w)\}=0$. Thus, the assumption $\inf_{w\in
  K}\{I(w)\}>0$ implies $v:=\inf\{K\}>-\infty$. As $K$ is a closed
set, $v\in K$ and consequently $I(v)\ge\inf_{w\in K}\{I(w)\}$. Since
$K\subseteq[v,+\infty)$, the Chernoff bound and limit
  (\ref{boundsup_cumul}) imply for every $k\ge 0$
\begin{align}
    \nonumber
  \limsup_{t\uparrow+\infty}\frac{1}{t}\ln\prob\bigg[\frac{F_t}{t}\in K\bigg]&\le
  \limsup_{t\uparrow+\infty}\frac{1}{t}\ln\prob\bigg[\frac{F_t}{t}\ge v\bigg]\\
  \nonumber
  &\le  \limsup_{t\uparrow+\infty}\frac{1}{t}\ln\prob\big[e^{kF_t-tvk}\big]=-vk-\vartheta(k).
\end{align}
 By optimizing with respect to $k$ we find the upper large deviation
 bound for the closed set $K$:
  \begin{equation*}
    \limsup_{t\uparrow+\infty}\frac{1}{t}\ln\prob\bigg[\frac{F_t}{t}\in K\bigg]
    \le -\sup_{k\ge 0}\big\{vk+\vartheta(k)\big\}=-I(v)\le -\inf_{w\in K}\{I(w)\}.
\end{equation*}

The case $\vartheta(k)=-\infty$ for all $k>0$ is tackled in a similar
way. Let us move to the case where $\vartheta(k)>-\infty$ for $k$ in
an open neighborhood of the origin. In this case, the subdifferential
$\partial$ of the convex function $-\vartheta$ at the origin is a
non-empty bounded closed interval (see \cite{Rockbook}, Theorem
23.4). We have $I(w)=0$ if and only if $w\in \partial$ (see
\cite{Rockbook}, Theorem 23.5). Given any point $w_o\in \partial$, we
can state that $I(w)=\sup_{k\le 0}\{wk+\vartheta(k)\}$ for $w\le w_o$
and $I(w)=\sup_{k\ge 0}\{wk+\vartheta(k)\}$ for $w\ge w_o$. In fact,
$w_ok+\vartheta(k)\le I(w_o)=0$ for all $k\in\Rl$, so that
$wk+\vartheta(k)=(w-w_o)k+w_ok+\vartheta(k)\le 0$ for all $k\ge 0$ if
$w\le w_o$ and $wk+\vartheta(k)=(w-w_o)k+w_ok+\vartheta(k)\le 0$ for
all $k\le 0$ if $w\ge w_o$.

Since $\inf_{w\in K}\{I(w)\}>0$, $\partial$ belongs to the complement
$K^c$ of $K$. Let $(v_-,v_+)$ be the union of all the open intervals
in $K^c$ that contain $\partial$. Notice that either $v_-$ or $v_+$
must be finite since $K$ is non-empty. Clearly, $v_-\le w_o$ and
$v_+\ge w_o$, so that $I(v_-)=\sup_{k\le 0}\{v_-k+\vartheta(k)\}$ and
$I(v_+)=\sup_{k\ge 0}\{v_+k+\vartheta(k)\}$.  If $v_+<+\infty$, then
$v_+\in K$ as $K$ is closed and consequently $I(v_+)\ge\inf_{w\in
  K}\{I(w)\}$. As before, the Chernoff bound and limit
(\ref{boundsup_cumul}) give
\begin{equation}
\limsup_{t\uparrow+\infty}\frac{1}{t}\ln\prob\bigg[\frac{F_t}{t}\ge v_+\bigg]\le -\sup_{k\ge 0}\big\{v_+k+\vartheta(k)\big\}=-I(v_+)\le-\inf_{w\in K}\{I(w)\}.
\label{upper_1001}
\end{equation}
Similarly, if $v_->-\infty$, then $I(v_-)\ge\inf_{w\in K}\{I(w)\}$ and
the Chernoff bound combined with limit (\ref{boundsup_cumul}) yields
\begin{equation}
\limsup_{t\uparrow+\infty}\frac{1}{t}\ln\prob\bigg[\frac{F_t}{t}\le v_-\bigg]\le -I(v_-)\le-\inf_{w\in K}\{I(w)\}.
\label{upper_1002}
\end{equation}
Bounds (\ref{upper_1001}) and (\ref{upper_1002}) prove the upper large
deviation bound for the closed set $K$ since
$K\subseteq(-\infty,w_-]\cup[w_+,+\infty)$.

\section{Proof of Theorem \ref{th:LDP}}
\label{proof:LDP}

\subsection{Proof of part $\boldsymbol{(i)}$}

The rate function $I=(-\varphi)^\star$ has compact level sets if some
minorant, such as the function that maps $w$ in
$\sup_{k\in\Rl}\{wk+\varphi(k)\wedge\ell\}\le I(w)$, has compact level
sets. By Theorem \ref{mainth} and Proposition \ref{prop:conv}, the
latter has compact level sets if $\Ex[e^{\rho|X_1|}]<+\infty$ for some
$\rho>0$.

As $I$ is the convex conjugate of the convex and lower semicontinuous
function $-\varphi$ that takes value zero at the origin, $I(w)=0$ if
and only if $w$ is a subgradient of $-\varphi$ at the origin (see
\cite{Rockbook}, Theorem 23.5). Assume that $\ell>0$ and that there
exists $\rho>0$ such that $\Ex[e^{\rho|X_1|}]<+\infty$. We claim that
$\varphi$ is differentiable at the origin with
$-\varphi'(0)=\mu$. Thus, $\mu$ is the only subgradient of $-\varphi$
at the origin. To begin with, we pick a finite number
$\ell_o\in(0,\ell)$ such that $\Ex[e^{\ell_oS_1}]>1$. The dominated
converge theorem shows that there exists $\rho_o\in(0,\rho)$ such that
  $1<\Ex[e^{\ell_o S_1+kX_1}]<+\infty$ for $|k|<\rho_o$. This way, the
  function $\Phi$ that maps the pair
  $(\zeta,k)\in(-\infty,\ell_o)\times(-\rho_o,\rho_o)$ in
  $\Phi(\zeta,k):=\Ex[e^{\zeta S_1+kX_1}]$ is analytic and
  $\Phi(\ell_o,k)>1$ for $|k|<\rho_o$. As
  $\varphi(k):=\sup\{\zeta\in\Rl:\Phi(\zeta,k)\le 1\}$, for each
  $k\in(-\rho_o,\rho_o)$ the number $\varphi(k)$ turns out to be the
  unique solution $\zeta$ of the equation $\Phi(\zeta,k)=1$, and
  $\varphi$ turns out to be analytic on $(-\rho_o,\rho_o)$ by the
  analytic implicit function theorem. The identity
  $\Phi(\varphi(k),k)=1$ for $|k|<\rho_o$ and the fact that
  $\varphi(0)=0$ give
  \begin{equation*}
    -\varphi'(0)=\frac{\frac{\partial\Phi}{\partial k}(0,0)}{\frac{\partial\Phi}{\partial\zeta}(0,0)}=\frac{\Ex[X_1]}{\Ex[S_1]}=:\mu.
\end{equation*}

\subsection{Proof of part $\boldsymbol{(ii)}$}

Assume that $\varpi(k)\ge\varphi(k)$ for all $k\in\Rl$. According to
part $(ii)$ of Proposition \ref{th:LDP_preliminar}, $F_t$ satisfies the
upper large deviation bound with the rate function
$I=(-\varphi)^\star$. Let us show that a lower large deviation bound
holds with the same rate function.  We prove such bound by a
truncation argument.  For each integer $i\ge 1$, let $\varphi_i$ be
the function that maps $k\in\Rl$ in
\begin{equation*}
\varphi_i(k):=\sup\Big\{\zeta\in\Rl ~:~ \Ex\big[e^{\zeta S_1+kX_1}\mathds{1}_{\{S_1\vee M_1\le i\}}\big]\le 1\Big\}.
\end{equation*}
By repeating some arguments of the proof of Proposition
\ref{prop:conv}, one can easily verify that the function $\varphi_i$
is concave and upper semicontinuous. One can also verify that, for any
$k$, the sequence $\{\varphi_i(k)\}_{i\ge 1}$ is non-increasing and
converges to $\varphi(k)$. It follows in particular that
$(-\varphi_i)^\star(w):=\sup_{k\in\Rl}\{wk+\varphi_i(k)\}$ is
non-increasing for every $w$, so that
\begin{equation}
\lim_{i\uparrow\infty}\,(-\varphi_i)^\star(w)=:J(w)
\label{J_limit}
\end{equation}
exists and $J(w)\le (-\varphi_i)^\star(w)$ for all $i$. We shall show
that for each open set $G\subseteq\Rl$ and point $w\in G$
\begin{equation}
  \liminf_{t\uparrow+\infty}\frac{1}{t}\ln\prob\bigg[\frac{F_t}{t}\in G\bigg]\ge -J(w).
\label{lower_F_00}
\end{equation}
This bound implies the lower large deviation bound as follows.  The
function $J$ that associates each $w$ with $J(w)$ is convex since it
is the pointwise limit of a sequence of convex functions. Convexity
entails that the convex conjugate
$J^{\star\star}(w):=\sup_{k\in\Rl}\{wk-J^\star(k)\}$ of the convex
conjugate $J^\star(k):=\sup_{w\in\Rl}\{kw-J(w)\}$ of $J$ equals the
lower-semicontinuous regularization of $J$ (see \cite{Rockbook},
Theorem 12.2), which is the function that maps $w$ in
$\lim_{\delta\downarrow 0}\inf_{v\in \Delta_{w,\delta}}\{J(v)\}$. Fix an
open set $G$ and a point $w\in G$. Since there exists $\delta_o>0$
such that $\Delta_{w,\delta_o}\subseteq G$, we have from (\ref{lower_F_00})
\begin{equation*}
  \liminf_{t\uparrow+\infty}\frac{1}{t}\ln\prob\bigg[\frac{F_t}{t}\in G\bigg]\ge -\inf_{v\in \Delta_{w,\delta_o}}\{J(v)\}\ge -\inf_{v\in \Delta_{w,\delta}}\{J(v)\}
\end{equation*}
for all $\delta\in(0,\delta_o)$, in such a way that
\begin{equation}
  \liminf_{t\uparrow+\infty}\frac{1}{t}\ln\prob\bigg[\frac{F_t}{t}\in G\bigg]\ge
  -\adjustlimits\lim_{\delta\downarrow 0}\inf_{v\in \Delta_{w,\delta}}\{J(v)\}=-J^{\star\star}(w).
\label{lower_F_01}
\end{equation}
At the same time, $J(w)\le(-\varphi_i)^\star(w)$ for all $w$ and $i$
gives $J^\star(k)\ge (-\varphi_i)^{\star\star}(k)=-\varphi_i(k)$ for
all $k$ and $i$ since $-\varphi_i$ is convex and lower
semicontinuous. By sending $i$ to infinity we realize that
$J^\star(k)\ge -\varphi(k)$ and $J^{\star\star}(w)\le
(-\varphi)^\star(w)=I(w)$. Thus, (\ref{lower_F_01}) implies
\begin{equation}
  \liminf_{t\uparrow+\infty}\frac{1}{t}\ln\prob\bigg[\frac{F_t}{t}\in G\bigg]\ge-I(w),
\end{equation}
which yields the lower large deviation bound after optimizing over
$w$.

Let us now verify (\ref{lower_F_00}) by introducing a truncated
probability distribution for waiting times and rewards. Fix once and
for all an open set $G\subseteq\Rl$ and a point $w\in G$. Assume that
$J(w)<+\infty$, otherwise there is nothing to prove. Let $s>0$ be a
number such that $m_1:=\frac{1}{2}\Ex[S_1\wedge s]>0$ and
$m_2:=2\Ex[(S_1\wedge s)^2]>0$. By the monotone convergence theorem
and (\ref{J_limit}) there exists $i_o\ge 1$ such that for all $i\ge
i_o$
\begin{equation}
  \begin{cases}
    \prob[S_1\vee M_1\le i]>0,\\
    \Ex[(S_1\wedge s)|S_1\vee M_1\le i]> m_1,\\
    \Ex[(S_1\wedge s)^2|S_1\vee M_1\le i]< m_2,\\
     \sup_{i\ge i_o}\{(-\varphi_i)^\star(w)\}<+\infty.
  \end{cases}
  \label{i_vincoli}
\end{equation}
For each $i\ge i_o$ we construct the truncated probability distribution
$Q_i$ for $(S_1,X_1)$ by conditioning on the event $S_1\vee M_1\le i$:
\begin{equation*}
Q_i(A):=\prob\big[(S_1,X_1)\in A \big| S_1\vee M_1\le i\big]=\frac{\prob[(S_1,X_1)\in A,\,S_1\vee M_1\le i]}{\prob[S_1\vee M_1\le i]}
\end{equation*}
for each Borel set $A\subseteq\Rl^2$. The Kolmogorov extension theorem
tells us that there exists a probability space
$(\Omega_i,\mathcal{F}_i,\prob_i)$ that hosts infinitely many
i.i.d.\ waiting time and reward pairs $(S_1,X_1),(S_2,X_2),\ldots$
distributed according to the law $Q_i$.  Pay attention to not confuse
$X_n$ on $(\Omega_i,\mathcal{F}_i,\prob_i)$ with
$X_n:=\int_0^{S_n}f(B_{n,\tau})\,d\tau$ on
$(\Omega,\mathcal{F},\prob)$.  Denoting expectation with respect to
$\prob_i$ by $\Ex_i$, set for $k\in\Rl$
\begin{equation*}
  \psi_i(k):=\sup\Big\{\zeta\in\Rl ~:~ \Ex_i\big[e^{\zeta S_1+kX_1}\big]\le 1\Big\}.
\end{equation*}
Since $\lim_{s\uparrow+\infty}-\frac{1}{s}\ln\prob_i[S_1>s]=+\infty$,
Theorem \ref{mainth} and Proposition \ref{prop:conv} state that $W_t$
satisfies a weak LDP with the rate function $J_i=(-\psi_i)^\star$ with
respect to the model $(\Omega_i,\mathcal{F}_i,\prob_i)$.

To proceed we observe that, by the last of (\ref{i_vincoli}), there
exists a large number $\eta\ge 1/m_1$ such that for all $i\ge i_o$
\begin{equation}
  \frac{(m_1\eta-1)^2}{2m_2\eta}\ge (-\varphi_i)^\star(w)+1.
\label{eta_prop}
\end{equation}
Fix $i\ge i_o$ and let $q_t$ be the integer part of $1+\eta t$.  As
$G$ is open, there exists $\delta>0$ such that
$(w-2\delta,w+2\delta)\subset G$. Bearing in mind that
$|F_t-W_t|\le\int_0^{S_{N_t+1}}|f(B_{N_t+1,\tau})|\,d\tau= M_{N_t+1}$,
for $i\ge i_o$ and $t>i/\delta$ we find $F_t/t\in G$ if
$W_t/t\in(w-\delta,w+\delta)$ and $M_{N_t+1}<i$. Thus, for $i\ge i_o$
and $t>i/\delta$ we have
\begin{align}
  \nonumber
\prob\bigg[\frac{F_t}{t}\in G\bigg]&\ge \prob\bigg[\frac{F_t}{t}\in G,\,N_t<q_t\bigg]\\
\nonumber
&\ge \prob\bigg[\frac{F_t}{t}\in G,\,N_t<q_t,\,\max\big\{S_1\vee M_1,\ldots,S_{q_t}\vee M_{q_t}\big\}\le i\bigg]\\
\nonumber
&\ge \prob\bigg[\frac{W_t}{t}\in(w-\delta,w+\delta),\,N_t<q_t,\,\max\big\{S_1\vee M_1,\ldots,S_{q_t}\vee M_{q_t}\big\}\le i\bigg]\\
\nonumber
&=\prob_i\bigg[\frac{W_t}{t}\in(w-\delta,w+\delta),\,N_t<q_t\bigg]\Big\{\prob\big[S_1\vee M_1\le i\big]\Big\}^{q_t}\\
&\ge\Big\{\prob\big[S_1\vee M_1\le i\big]\Big\}^{q_t}\bigg\{\prob_i\bigg[\frac{W_t}{t}\in(w-\delta,w+\delta)\bigg]
-\prob_i\big[N_t\ge q_t\big]\bigg\}.
\label{lower_F_11}
\end{align}

The lower large deviation bound for $W_t$ with respect to the model
$(\Omega_i,\mathcal{F}_i,\prob_i)$ gives
\begin{equation}
  \liminf_{t\uparrow+\infty}\frac{1}{t}\ln\prob_i\bigg[\frac{W_t}{t}\in(w-\delta,w+\delta)\bigg]\ge-J_i(w)=-(-\psi_i)^\star(w).
\label{psi_lower}
\end{equation}
We now point out that for each $k$
\begin{align}
  \nonumber
  \psi_i(k)&:=\sup\Big\{\zeta\in\Rl ~:~ \Ex_i\big[e^{\zeta S_1+kX_1}\big]\le 1\Big\}\\
  \nonumber
  &=\sup\Big\{\zeta\in\Rl ~:~ \Ex\big[e^{\zeta S_1+kX_1}\mathds{1}_{\{S_1\vee M_1\le i\}}\big]\le \prob\big[S_1\vee M_1\le i\big]\Big\}\\
  \nonumber
   &\le\sup\Big\{\zeta\in\Rl ~:~ \Ex\big[e^{\zeta S_1+kX_1}\mathds{1}_{\{S_1\vee M_1\le i\}}\big]\le 1\Big\}=:\varphi_i(k),
\end{align}
so that (\ref{psi_lower}) entails
\begin{equation}
  \liminf_{t\uparrow+\infty}\frac{1}{t}\ln\prob_i\bigg[\frac{W_t}{t}\in (w-\delta,w+\delta)\bigg]\ge-(-\psi_i)^\star(w)\ge-(-\varphi_i)^\star(w).
\label{lower_F_12}
\end{equation}

In order to estimate $\prob_i[N_t\ge q_t]$ we make use of the Chernoff
bound, which yields for each $\lambda\ge 0$
\begin{equation*}
  \prob_i\big[N_t\ge q_t\big]=\prob_i\big[T_{q_t}\le t\big]\le e^{\lambda t}\,\Ex_i\big[e^{-\lambda T_{q_t}}\big]
  =e^{\lambda t}\Big(\Ex_i\big[e^{-\lambda S_1}\big]\Big)^{q_t}\le e^{\lambda t}\Big(\Ex_i\big[e^{-\lambda (S_1\wedge s)}\big]\Big)^{q_t},
\end{equation*}
so that
\begin{equation*}
  \limsup_{t\uparrow+\infty}\frac{1}{t}\ln\prob_i\big[N_t\ge q_t\big]\le \lambda+\eta\ln\Ex_i\big[e^{-\lambda (S_1\wedge s)}\big].
\end{equation*}
On the other hand, by combining the inequalities $e^{-x}\le
1+x+\frac{1}{2}x^2$ valid for $x\ge 0$ and $1+x\le e^x$ valid for any $x$
with the second and third of (\ref{i_vincoli}), we realize that
\begin{align}
  \nonumber
  \Ex_i\big[e^{-\lambda (S_1\wedge s)}\big]&=\Ex\big[e^{-\lambda (S_1\wedge s)}\big|S_1\vee M_1\le i\big]\\
  \nonumber
  &\le 1-\lambda\,\Ex\big[(S_1\wedge s)\big|S_1\vee M_1\le i\big]+\frac{\lambda^2}{2}
  \Ex\big[(S_1\wedge s)^2\big|S_1\vee M_1\le i\big]\\
  \nonumber
  &\le 1-\lambda m_1+\frac{\lambda^2}{2}m_2\le e^{-\lambda m_1+\frac{\lambda^2}{2}m_2}.
\end{align}
Thus, we find
\begin{equation*}
  \limsup_{t\uparrow+\infty}\frac{1}{t}\ln\prob_i\big[N_t\ge q_t\big]\le (1-m_1\eta)\lambda+m_2\eta\frac{\lambda^2}{2}.
\end{equation*}
By optimizing over $\lambda$, i.e.\ by taking
$\lambda:=(m_1\eta-1)/(m_2\eta)\ge 0$, and by recalling
(\ref{eta_prop}) we finally obtain
\begin{equation}
\limsup_{t\uparrow+\infty}\frac{1}{t}\ln\prob_i\big[N_t\ge q_t\big]\le -\frac{(m_1\eta-1)^2}{2m_2\eta}\le -(-\varphi)^\star(w)-1.
\label{lower_F_13}
\end{equation}

Bounds (\ref{lower_F_12}) and (\ref{lower_F_13}) show that
$\prob_i[N_t\ge q_t]$ decays at a faster rate than $\prob_i[W_t/t\in
  (w-\delta,w+\delta)]$. Then, (\ref{lower_F_11}) gives
\begin{align}
  \nonumber
  \liminf_{t\uparrow+\infty}\frac{1}{t}\ln\prob\bigg[\frac{F_t}{t}\in G\bigg]
  \ge\eta\ln\prob\big[S_1\vee M_1\le i\big]-(-\varphi_i)^\star(w).
\end{align}
We get (\ref{lower_F_00}) from here by sending $i$ to infinity as
$\lim_{i\uparrow\infty}\prob[S_1\vee M_1\le i]=1$.

\subsection{Proof of part $\boldsymbol{(iii)}$}

Suppose that $\varpi(k)\ge\varphi(k)$ for all $k\in\Rl$. Part $(i)$ of
Proposition \ref{th:LDP_preliminar} shows that for every $k$
\begin{equation*}
\limsup_{t\uparrow+\infty}\frac{1}{t}\ln\Ex\big[e^{kF_t}\big]\le-\varphi(k).
\end{equation*}
Let  us verify that for every $k$ we also have
\begin{equation}
  \label{boundinf_cumul}
\liminf_{t\uparrow+\infty}\frac{1}{t}\ln\Ex\big[e^{kF_t}\big]\ge-\varphi(k).
\end{equation}

Pick $k\in\Rl$ and fix numbers $w\in\Rl$ and $\epsilon>0$. As $kF_t\ge
twk-t\epsilon|k|$ when $F_t/t\in(w-\epsilon,w+\epsilon)$, for each $t$
we find the bound
\begin{equation*}
  \Ex\big[e^{kF_t}\big]\ge  \Ex\bigg[e^{kF_t}\mathds{1}_{\big\{\frac{F_t}{t}\in(w-\epsilon,w+\epsilon)\big\}}\bigg]\ge
  e^{twk-t\epsilon|k|}\prob\bigg[\frac{F_t}{t}\in(w-\epsilon,w+\epsilon)\bigg].
\end{equation*}
Thus, the lower large deviation bound gives
\begin{equation*}
  \liminf_{t\uparrow+\infty}\frac{1}{t}\ln\Ex\big[e^{kF_t}\big]\ge wk-I(w)-\epsilon|k|,
\end{equation*}
and the arbitrariness of $w$ and $\epsilon$ implies
\begin{equation*}
  \liminf_{t\uparrow+\infty}\frac{1}{t}\ln\Ex\big[e^{kF_t}\big]\ge \sup_{w\in\Rl}\big\{wk-I(w)\big\}.
\end{equation*}
This result is (\ref{boundinf_cumul}) because the convex conjugate of
$I$ is $-\varphi$ since $\varphi$ is upper semicontinuous.

\section{Proof of Corollary \ref{corol}}
\label{corol:proof}

Assume that $\Ex[e^{\rho M_1}]<+\infty$ for every $\rho>0$. In order
to prove this corollary of Theorem \ref{th:LDP} we show that the
condition $\varpi(k)\ge\varphi(k)$ is satisfied for all $k\in\Rl$, with
$\varpi(k)$ defined by (\ref{def:working_hypothesis}). Fix $k\in\Rl$ such
that $\varphi(k)>-\infty$ and pick $\epsilon\in(0,1)$.  H\"older's
inequality gives for all $t>0$
\begin{align}
  \nonumber
  \mathcal{E}_t(k)\prob[S_1>t]=\Ex\Big[e^{k\int_0^tf(B_{1,\tau})\,d\tau}\mathds{1}_{\{S_1>t\}}\Big]&\le \Ex\Big[e^{|k|\int_0^{S_1}|f(B_{1,\tau})|\,d\tau}\mathds{1}_{\{S_1>t\}}\Big]\\
  \nonumber
  &\le \Big\{\Ex\big[e^{(|k|/\epsilon)M_1}\big]\Big\}^\epsilon\big\{\prob[S_1>t]\big\}^{1-\epsilon},
\end{align}
so that
\begin{equation*}
\varpi(k):=\liminf_{t\uparrow+\infty}-\frac{1}{t}\ln\Big\{\mathcal{E}_t(k)\prob[S_1>t]\Big\}\ge \ell(1-\epsilon).
\end{equation*}
By sending $\epsilon$ to zero we get $\varpi(k)\ge\ell$.  Thus,
$\varpi(k)\ge\varphi(k)$ is verified if we demonstrate that
$\varphi(k)\le\ell$, which trivially holds when $\ell=+\infty$.
Although it is hard to exhibit a function $f$ for which $\Ex[e^{\rho
    M_1}]<+\infty$ for all $\rho>0$ when $\ell<+\infty$, we are not
able to exclude that such a function exists.

Suppose that $\ell<+\infty$.  Set $f_\pm:=(\pm f)\vee 0$ and
$g_\pm(t):=\Ex[\int_0^tf_\pm(B_{1,\tau})\,d\tau]\ge 0$ for all $t\ge
0$. We show at first that $\limsup_{t\uparrow+\infty}g_\pm(t)/t\le
0$. Pick $\rho>0$ and denote by $P:=\prob[S_1\in\cdot\,]$ the
probability measure induced on $[0,+\infty)$ by $S_1$. Jensen's
  inequality gives $e^{\rho g_\pm(s)}\le \Ex\big[e^{\rho\int_0^s
      f_\pm(B_{1,\tau})\,d\tau}\big]\le
  \Ex\big[e^{\rho\int_0^s|f(B_{1,\tau})|\,d\tau}\big]$ for all $s\ge 0$. We
  use this bound and the fact that $g_\pm$ is non-decreasing to state
  that for every $t>0$
  \begin{align}
    \nonumber
    e^{\rho g_\pm(t)}\,\prob[S_1>t]&=\int_{(t,+\infty)}e^{\rho g_\pm(t)}P(ds)\le\int_{[0,+\infty)}e^{\rho g_\pm(s)}P(ds)\\
      \nonumber
      &\le\int_{[0,+\infty)}\Ex\Big[e^{\rho\int_0^s|f(B_{1,\tau})|\,d\tau}\Big]P(ds)\\
        \nonumber
        &=\Ex\Big[e^{\rho\int_0^{S_1}|f(B_{1,\tau})|\,d\tau}\Big]=\Ex\big[e^{\rho M_1}\big]<+\infty.
    \end{align}
It follows that
 \begin{equation*}
    \rho\limsup_{s\uparrow+\infty}\frac{g_\pm(t)}{t}-\ell=\limsup_{t\uparrow+\infty}\frac{1}{t}\ln \Big\{e^{\rho g_\pm(t)}\prob[S_1>t]\Big\}\le 0,
\end{equation*}
and the arbitrariness of $\rho$ implies
$\limsup_{t\uparrow+\infty}g_\pm(t)/t\le 0$.

We are now ready to prove that $\varphi(k)\le\ell$.  Fix any
$\epsilon>0$.  Since $\limsup_{t\uparrow+\infty}g_\pm(t)/t\le 0$,
there exists $t>0$ with the property that $0\le g_\pm(s)\le\epsilon s$
for all $s>t$, in such a way that
$\Ex\big[k\int_0^sf(B_{1,\tau})\,d\tau\big]=-kg_-(s)+kg_+(s)\ge-\epsilon|k|s$
for $s>t$. Once again, we appeal to Jensen's inequality to get
\begin{align}
  \nonumber
  \Ex\Big[e^{\{\varphi(k)-\epsilon|k|\}S_1}\mathds{1}_{\{S_1>t\}}\Big]&\le\Ex\Big[e^{\varphi(k) S_1-kg_-(S_1)+kg_+(S_1)}\Big]\\
  \nonumber
  &\le\Ex\Big[e^{\varphi(k) S_1+k\int_0^{S_1}f(B_{1,\tau})\,d\tau}\Big]=\Ex\big[e^{\varphi(k) S_1+kX_1}\big]\le 1.
  \end{align}
This bound demonstrates that $\varphi(k)-\epsilon|k|\le\ell$ since
$\Ex[e^{\zeta S_1}\mathds{1}_{\{S_1>t\}}]=+\infty$ for
$\zeta>\ell$. The arbitrariness of $\epsilon$ gives
$\varphi(k)\le\ell$.

\section{Proof of Corollary \ref{corol_lim}}
\label{proof:corol_lim}

According to Theorem \ref{th:LDP}, it suffices to verify that
$\varpi(k)\ge\varphi(k)$ for any given $k$ such that
$\varphi(k)>-\infty$. Recall that
$\int_{[0,+\infty)}e^{\varphi(k)s}\mathcal{E}_s(k)P(ds)\le 1$ for such
  a $k$, where $P:=\prob[S_1\in\cdot\,]$ is the probability measure
  induced on $[0,+\infty)$ by $S_1$.

Set $\lambda:=\lim_{t\uparrow+\infty}\frac{1}{t}\ln\mathcal{E}_t(k)$,
which exists by assumption. If $\lambda=-\infty$ or $\ell=+\infty$,
then $\varpi(k)=+\infty$ and the inequality $\varpi(k)\ge\varphi(k)$
trivially holds. Suppose that $\lambda>-\infty$ and that
$\ell<+\infty$.  We address first the case $\lambda<+\infty$.  In this
case we have $\varpi(k)=\ell-\lambda$ and $\mathcal{E}_t(k)\ge
e^{(\lambda-\epsilon)t}$ for an arbitrary number $\epsilon>0$ and all
sufficiently large $t$. Thus, for all sufficiently large $t$ we find
\begin{equation*}
  \Ex\big[e^{\varphi(k)S_1+(\lambda-\epsilon)S_1}\mathds{1}_{\{S_1>t\}}\big]=\int_{(t,+\infty)}e^{\varphi(k)s+(\lambda-\epsilon)s}P(ds)
  \le\int_{(t,+\infty)}e^{\varphi(k)s}\mathcal{E}_s(k)P(ds)\le 1.
\end{equation*}
This bound implies $\varphi(k)+(\lambda-\epsilon)\le\ell$ since
$\Ex[e^{\zeta S_1}\mathds{1}_{\{S_1>t\}}]=+\infty$ for
$\zeta>\ell$. The arbitrariness of $\epsilon$ gives
$\ell-\lambda\ge\varphi(k)$, i.e.\ $\varpi(k)\ge\varphi(k)$.

Let us consider now the case $\lambda=+\infty$. In this case
$\frac{1}{t}\ln\mathcal{E}_t(k)$ is eventually non-decreasing with
respect to $t$ by assumption. We show that
$e^{\varphi(k)t}\mathcal{E}_t(k)\prob[S_1>t]\le 1$ for all
sufficiently large $t$, which gives $\varpi(k)\ge\varphi(k)$. In fact,
let $t$ be a sufficiently large time such that
$\frac{1}{t}\ln\mathcal{E}_t(k)\le\frac{1}{s}\ln\mathcal{E}_s(k)$ for
every $s>t$. If $\varphi(k)+\frac{1}{t}\ln\mathcal{E}_t(k)\le 0$, then
the bound $e^{\varphi(k)t}\mathcal{E}_t(k)\prob[S_1>t]\le 1$ is
trivial. If $\varphi(k)+\frac{1}{t}\ln\mathcal{E}_t(k)>0$, then
$\varphi(k)t+\ln\mathcal{E}_t(k)\le\varphi(k)s+\frac{s}{t}\ln\mathcal{E}_t(k)\le\varphi(k)s+\ln\mathcal{E}_s(k)$
for each $s>t$, so that
\begin{equation*}
  e^{\varphi(k)t}\mathcal{E}_t(k)\prob[S_1>t]=\int_{(t,+\infty)}e^{\varphi(k)t+\ln\mathcal{E}_t(k)}P(ds)
  \le\int_{(t,+\infty)}e^{\varphi(k)s}\mathcal{E}_s(k)P(ds)\le 1.
\end{equation*}

\section{Proof of Lemma \ref{lem:cases}}
\label{proof:cases}

The lemma relies on the fact that the uniform norm
$\|B_1\|_\infty:=\sup_{\tau\in[0,1]}\{|B_{1,\tau}|\}$ of the Brownian
motion has a Gaussian tail:
\begin{equation}
  \lim_{x\uparrow+\infty}\frac{1}{x^2}\ln\prob\big[\|B_1\|_\infty>x\big]=-\lambda
\label{tail_norm}
\end{equation}
with some $\lambda>0$ (see \cite{Talagrand}, Corollary 3.2).

\subsection{Proof of part $\boldsymbol{(i)}$}

Part $(i)$ is proved as follows.  Pick $\rho>0$ and denote by
$P:=\prob[S_1\in\cdot\,]$ the probability measure induced on
$[0,+\infty)$ by $S_1$. Since $|f(z)|\le A(1+|z|^\alpha)$ for all
  $z\in\Rl$ with $A>0$ and $\alpha\ge 0$, we have
  $M_1:=\int_0^{S_1}|f(B_{1,\tau})|\,d\tau\le
  AS_1+A\int_0^{S_1}|B_{1,\tau}|^\alpha\,d\tau$. Since
  $\int_0^s|B_{1,\tau}|^\alpha=s\int_0^1|B_{1,s\tau}|^\alpha$ is
  distributed as $s^{1+\frac{\alpha}{2}}\int_0^1|B_{1,\tau}|^\alpha$
  for each $s\ge 0$, we can state that
\begin{align}
  \nonumber
  \Ex\big[M_1^\rho\big]&\le A^\rho \int_{[0,+\infty)}\Ex\bigg[\bigg(s+\int_0^s|B_{1,\tau}|^\alpha\,d\tau\bigg)^\rho\,\bigg]P(ds)\\
  \nonumber
  &= A^\rho \int_{[0,+\infty)}\Ex\bigg[\bigg(s+s^{1+\frac{\alpha}{2}}\int_0^1|B_{1,\tau}|^\alpha\,d\tau\bigg)^\rho\,\bigg]P(ds)\\
    \nonumber
    &\le A^\rho \int_{[0,+\infty)}(1\vee s)^{\rho(1+\frac{\alpha}{2})}\,\Ex\Big[\big(1+\|B_{1,\tau}\|_\infty^\alpha\big)^\rho\Big]P(ds)\\
      \nonumber
     &=A^\rho\,\Ex\Big[(1\vee S_1)^{\rho(1+\frac{\alpha}{2})}\Big]\Ex\Big[\big(1+\|B_{1,\tau}\|_\infty^\alpha\big)^\rho\Big].
\end{align}
Thus, $\Ex[S_1^{\rho(1+\frac{\alpha}{2})}]<+\infty$ implies
$\Ex[M_1^\rho]<+\infty$ since
$\Ex[(1+\|B_1\|_\infty^\alpha)^\rho]<+\infty$ for each $\alpha\ge 0$
and $\rho>0$.

\subsection{Proof of part $\boldsymbol{(ii)}$}

Regarding part $(ii)$, we suppose that $\alpha\in[0,2)$, we fix
  $\rho>0$, and we write
  \begin{align}
    \nonumber
    \Ex\big[e^{\rho M_1}\big]&\le\int_{[0,+\infty)}\Ex\Big[e^{\rho As+\rho A s^{\frac{\alpha+2}{2}}\int_0^1|B_{1,\tau}|^\alpha\,d\tau}\Big]P(ds)\\
      \nonumber
    &\le\int_{[0,+\infty)}\Ex\Big[e^{\rho As+\rho A s^{\frac{\alpha+2}{2}}\|B_{1,\tau}\|_\infty^\alpha}\Big]P(ds).
  \end{align}
The function that maps $x\ge 0$ in $x^{\frac{\alpha}{2}}$ is concave as
$\alpha\in[0,2)$. Then, for every $\beta\ge 0$ we have
\begin{equation}
    \|B_1\|_\infty^\alpha=\big(\|B_1\|_\infty^2\big)^{\frac{\alpha}{2}}\le \beta^{\frac{\alpha}{2}}+\frac{\alpha}{2}\beta^{\frac{\alpha}{2}-1}\big(\|B_1\|_\infty^2-\beta\big)
    \le \beta^{\frac{\alpha}{2}}+\frac{\alpha}{2}\beta^{-\frac{2-\alpha}{2}}\|B_1\|_\infty^2.
\label{mom_1}
\end{equation}
By taking $\beta:=s^{\frac{2+\alpha}{2-\alpha}}(\alpha\rho
A/\lambda)^{\frac{2}{2-\alpha}}$ in (\ref{mom_1}) with some $s\ge 0$
we realize that
\begin{align}
  \nonumber
  \rho As+ \rho A s^{\frac{\alpha+2}{2}}\|B_1\|_\infty^\alpha &\le \rho As+
  (\alpha/\lambda)^{\frac{\alpha}{2-\alpha}}(\rho A)^{\frac{2}{2-\alpha}}s^{\frac{2+\alpha}{2-\alpha}}+\frac{\lambda}{2}\|B_1\|_\infty^2\\
  \nonumber
  &\le \eta\Big(1+s^{\frac{2+\alpha}{2-\alpha}}\Big)+\frac{\lambda}{2}\|B_1\|_\infty^2, 
\end{align}
where we have set $\eta:=\rho
A+(\alpha/\lambda)^{\frac{\alpha}{2-\alpha}}(\rho
A)^{\frac{2}{2-\alpha}}$ for brevity. Thus
\begin{equation*}
  \Ex\big[e^{\rho M_1}\big]\le e^{\eta}\,\Ex\Big[e^{\eta S_1^{\frac{2+\alpha}{2-\alpha}}}\Big]\,
  \Ex\Big[e^{\frac{\lambda}{2}\|B_1\|_\infty^2}\Big].
\end{equation*}
Fubini's theorem and (\ref{tail_norm}) give
\begin{align}
  \nonumber
  \Ex\Big[e^{\frac{\lambda}{2}\|B_1\|_\infty^2}\Big]&=\Ex\bigg[1+\frac{\lambda}{2}\int_0^{+\infty}e^{\frac{\lambda}{2}x}\mathds{1}_{\{\|B_1\|_\infty>\sqrt{x}\}}dx\bigg]\\
  \nonumber
  &=1+\frac{\lambda}{2}\int_0^{+\infty}e^{\frac{\lambda}{2}x}\,\prob\big[\|B_1\|_\infty>\sqrt{x}\big]dx<+\infty.
\end{align}
Fubini's theorem and the hypothesis $\limsup_{s\uparrow+\infty}
s^{-\frac{2+\alpha}{2-\alpha}}\ln\prob[S_1>s]=-\infty$ give
\begin{align}
  \nonumber
 \Ex\Big[e^{\eta S_1^{\frac{2+\alpha}{2-\alpha}}}\Big]&=\Ex\bigg[1+\eta\int_0^{+\infty}e^{\eta s}\mathds{1}_{\big\{S_1>s^{\frac{2-\alpha}{2+\alpha}}\big\}}ds\bigg]\\
  \nonumber
  &=1+\eta\int_0^{+\infty}e^{\eta s}\,\prob\big[S_1>s^{\frac{2-\alpha}{2+\alpha}}\big]ds<+\infty.
\end{align}

\section{Proof of Lemma \ref{lem:varphi_time}}
\label{proof:varphi_time}

The function $\Phi$ that maps $(\zeta,k)\in\Rl^2$ in
\begin{equation*}
    \Phi(\zeta,k):=\Ex\big[e^{\zeta S_1+kX_1}\big]=\int_0^1\frac{\Ex[e^{(\zeta+kx)S_1}]}{\pi\sqrt{x(1-x)}}\,dx
\end{equation*}
is finite and analytic throughout $\Rl^2$ if
$\ell:=\lim_{s\uparrow+\infty}-\frac{1}{s}\ln\prob[S_1>s]=+\infty$,
since $\Ex[e^{\zeta S_1}]<+\infty$ for all $\zeta\in\Rl$ in such
case. If instead $\ell<+\infty$, then $\Ex[e^{\zeta S_1}]<+\infty$ for
$\zeta<\ell$ and $\Ex[e^{\zeta S_1}]=+\infty$ for $\zeta>\ell$. In
this second case $\Phi$ is finite and analytic in the region of pairs
$(\zeta,k)$ with $\zeta<\ell\wedge(\ell-k)$, whereas
$\Phi(\zeta,k)=+\infty$ for $\zeta>\ell\wedge(\ell-k)$. It follows
that $\varphi(k)=\sup\{\zeta\in\Rl:\Phi(\zeta,k)\le
1\}\le\ell\wedge(\ell-k)$ for all $k$. Thus, $\varphi(k)\le-0\vee k$
when $\ell=0$, which gives part $(iii)$ of the lemma as we already
know that $\varphi(k)\ge-0\vee k$.

We are going to identify $\varphi$, and then to prove parts $(i)$ and
$(ii)$ of the lemma.  If $\ell=+\infty$, then for each $k\in\Rl$ the
value $\varphi(k)$ of $\varphi$ turns out to be the unique real number
$\zeta$ that solves the equation $\Phi(\zeta,k)=1$, and the function
$\varphi$ is analytic throughout $\Rl$ by the analytic implicit
function theorem. If $0<\ell<+\infty$ and
$\Phi(\ell\wedge(\ell-k),k)>1$, then $\varphi(k)$ is the unique real
number $\zeta$ such that $\Phi(\zeta,k)=1$, whereas if
$0<\ell<+\infty$ and $\Phi(\ell\wedge(\ell-k),k)\le 1$, then
$\varphi(k)=\ell\wedge(\ell-k)$. These arguments show that the case
$0<\ell<+\infty$ requires to investigate the behavior of
$\Phi(\zeta,k)$ along the curve $\zeta=\ell\wedge(\ell-k)$ in order to
assess $\varphi(k)$. The identity
\begin{equation*}
  \Phi(\ell\wedge(\ell-k),k)=\int_0^1\frac{\Ex[e^{(\ell-|k|x)S_1}]}{\pi\sqrt{x(1-x)}}\,dx
\end{equation*}
states that $\Phi(\ell\wedge(\ell-k),k)$ is a convex function of
$|k|$.  Let us demonstrate that there exist two positive constants
$C_-$ and $C_+$ such that for all $k\ne 0$
\begin{equation}
  C_-\frac{\Lambda-e^{\frac{\ell}{|k|}}}{\sqrt{|k|}}\le\Phi(\ell\wedge(\ell-k),k)\le C_+\Ex\bigg[\frac{C_+e^{\ell S_1}}{\sqrt{1+|k|S_1}}\bigg]\le
  C_+\frac{\Lambda}{\sqrt{|k|\wedge 1}}
\label{bound_Phi_curve}
\end{equation}
with
\begin{equation*}
\Lambda:=\Ex\bigg[\frac{e^{\ell S_1}}{\sqrt{1+S_1}}\bigg].
\end{equation*}
The limit
\begin{equation*}
\lim_{s\uparrow+\infty}\sqrt{s}\int_0^1\frac{e^{-s x}}{\pi\sqrt{x(1-x)}}\,dx=\int_0^{+\infty}\frac{e^{- x}}{\pi\sqrt{x}}\,dx
\end{equation*}
shows that there exist two constants $C_->0$ and $C_+>0$ such that
\begin{equation}
\frac{C_-}{\sqrt{s}}\le\int_0^1\frac{e^{-s x}}{\pi\sqrt{x(1-x)}}\,dx\le\frac{C_+}{\sqrt{1+s}},
\label{lower_upper_int}
\end{equation}
the lower bound being valid for all $s>1$ and the upper bound being
valid for all $s\ge0$. Then, Fubini's theorem yields for $k\ne 0$
\begin{equation*}
  \Phi(\ell\wedge(\ell-k),k)=\int_0^1\frac{\Ex[e^{(\ell-|k|x)S_1}]}{\pi\sqrt{x(1-x)}}\,dx\ge
  \Ex\bigg[\frac{C_-e^{\ell S_1}\mathds{1}_{\{|k|S_1>1\}}}{\sqrt{|k|S_1}}\bigg]\ge C_-\frac{\Lambda-e^{\frac{\ell}{|k|}}}{\sqrt{|k|}}
\end{equation*}
and
\begin{equation*}
  \Phi(\ell\wedge(\ell-k),k)=\int_0^1\frac{\Ex[e^{(\ell-|k|x)S_1}]}{\pi\sqrt{x(1-x)}}\,dx\le
  \Ex\bigg[\frac{C_+e^{\ell S_1}}{\sqrt{1+|k|S_1}}\bigg].
\end{equation*}

Bounds (\ref{bound_Phi_curve}) together with the identity
$\Phi(\ell\wedge(\ell-k),k)=\Ex[e^{\ell S_1}]$ for $k=0$ show that
$\Phi(\ell\wedge(\ell-k),k)=+\infty$ for all $k\in\Rl$ if
$\Lambda=+\infty$. Then, as in the case $\ell=+\infty$, also in the
case $0<\ell<+\infty$ and $\Lambda=+\infty$ it turns out that
$\varphi(k)$ for every $k\in\Rl$ is the unique real number $\zeta$
that solves the equation $\Phi(\zeta,k)=1$, and the function $\varphi$
is analytic throughout $\Rl$ by the analytic implicit function
theorem. If instead $0<\ell<+\infty$ and $\Lambda<+\infty$, then
bounds (\ref{bound_Phi_curve}) tell us that
$\Phi(\ell\wedge(\ell-k),k)$ is a finite convex function of $|k|\ne
0$, which goes to zero as $|k|$ is sent to infinity by the dominated
converge theorem. Convexity and finiteness imply continuity and the
monotone converge theorem states that $\lim_{|k|\downarrow
  0}\Phi(\ell\wedge(\ell-k),k)=\Ex[e^{\ell S_1}]>1$. Thus, if
$0<\ell<+\infty$ and $\Lambda<+\infty$, then there exists $\lambda>0$
with the property that $\Phi(\ell\wedge(\ell-k),k)>1$ for
$|k|<\lambda$, $\Phi(\ell\wedge(\ell-k),k)=1$ for $|k|=\lambda$, and
$\Phi(\ell\wedge(\ell-k),k)<1$ for $|k|>\lambda$. The number $\lambda$
solves the equation
\begin{equation*}
  1=\Phi(\ell,-\lambda)=\Ex\big[e^{\ell S_1-\lambda X_1}\big].
\end{equation*}
In conclusion, if $0<\ell<+\infty$ and $\Lambda<+\infty$, then
$\varphi(k)$ is the unique real number $\zeta$ that solves the
equation $\Phi(\zeta,k)=1$ for $|k|<\lambda$, whereas
$\varphi(k)=\ell\wedge(\ell-k)$ for $|k|\ge\lambda$. The analytic
implicit function theorem gives that $\varphi$ is analytic on the open
interval $(-\lambda,\lambda)$, but $\varphi$ may be not differentiable
at the boundary points $-\lambda$ and $\lambda$.

\subsection{Proof of part $\boldsymbol{(i)}$}

If $\ell=+\infty$ or $0<\ell<+\infty$ and $\Lambda=+\infty$, then
$\varphi$ is an analytic function throughout $\Rl$ that satisfies
$\Ex[e^{\varphi(k)S_1+kX_1}]=\Phi(\varphi(k),k)=1$ and
\begin{equation*}
  \varphi'(k)=-\frac{\frac{\partial\Phi}{\partial k}(\varphi(k),k)}{\frac{\partial\Phi}{\partial
      \zeta}(\varphi(k),k)}= -\frac{\Ex[X_1e^{\varphi(k)S_1+kX_1}]}{\Ex[S_1e^{\varphi(k) S_1+kX_1}]}
\end{equation*}
for all $k\in\Rl$. Let us verify that
$\lim_{k\downarrow-\infty}\varphi'(k)=0$ and
$\lim_{k\uparrow+\infty}\varphi'(k)=-1$. The symmetry
$\varphi(-k)=\varphi(k)+k$ makes it sufficient to address
$\lim_{k\uparrow+\infty}\varphi'(k)=-1$ only. We have
$\lim_{k\uparrow+\infty}\varphi'(k)\ge-1$ since $X_1\le S_1$. We show
that $\lim_{k\uparrow+\infty}\varphi'(k)\le-1$ by
contradiction. Assume that $\lim_{k\uparrow+\infty}\varphi'(k)>-1$, so
that there exists $\eta>0$ such that $\varphi'(k)>-1+2\eta$ for all
sufficiently large $k$. Let $\delta>0$ be such that
$\prob[S_1>\delta]>0$. Since $\varphi(0)=0$, concavity gives
$\varphi(k)\ge\varphi'(k)k\ge(-1+2\eta)k$ for all sufficiently large
$k$. Then, under the conditions $S_1>\delta$ and $X_1>(1-\eta)S_1$,
for all sufficiently large $k$ we find $\varphi(k)S_1+kX_1\ge
(-1+2\eta)k S_1+(1-\eta)kS_1=\eta kS_1\ge\delta\eta k$, so that
\begin{align}
  \nonumber
  1=\Ex\big[e^{\varphi(k)S_1+kX_1}\big]&\ge e^{\delta\eta k}\,\prob\big[S_1>\delta,\,X_1>(1-\eta)S_1\big]\\
\nonumber
  &=e^{\delta\eta k}\,\prob[S_1>\delta]\bigg(1-\frac{2}{\pi}\arcsin\sqrt{1-\eta}\bigg).
\end{align}
This bound is a contradiction since the r.h.s.\ goes to infinity as
$k$ is sent to infinity.

\subsection{Proof of part $\boldsymbol{(ii)}$}

In the case $0<\ell<+\infty$ and $\Lambda<+\infty$ the function
$\varphi$ is analytic on the open interval $(-\lambda,\lambda)$, and
for each $k\in(-\lambda,\lambda)$ fulfills
$\Ex[e^{\varphi(k)S_1+kX_1}]=\Phi(\varphi(k),k)=1$ and
\begin{equation}
  \varphi'(k)= -\frac{\Ex[X_1e^{\varphi(k)S_1+kX_1}]}{\Ex[S_1e^{\varphi(k) S_1+kX_1}]}.
\label{varphi_derivative}
\end{equation}
We have $\varphi(k)=\ell\wedge(\ell-k)$ for
$k\notin(-\lambda,\lambda)$.  Let us compute the left derivative
$\varphi_-'(k)$ and the right derivative $\varphi_+'(k)$ of $\varphi$
at the points $k=-\lambda$ and $k=\lambda$, which exist because
$\varphi$ is concave and finite on the whole $\Rl$. Once again, the
symmetry $\varphi(-k)=\varphi(k)+k$ tells us that it suffices to
address only the instance $k=\lambda$. As $\varphi(k)=\ell-k$ for
$k\ge\lambda$ we find $\varphi_+'(\lambda)=-1$ and, by concavity,
$\varphi_-'(\lambda)\ge\varphi_+'(\lambda)=-1$ (see \cite{Rockbook},
Theorem 24.1). In order to calculate the left derivative we observe at
first that bounds (\ref{lower_upper_int}) can be used
as before to obtain for $k>0$
\begin{equation*}
  C_-\frac{\Xi-k^{-\frac{1}{2}}e^{\frac{\ell}{k}}}{\sqrt{k}}\le\Ex\big[S_1e^{(\ell-k)S_1+k X_1}\big]
  =\int_0^1\frac{\Ex[S_1e^{(\ell-k x)S_1}]}{\pi\sqrt{x(1-x)}}\,dx\le C_+\frac{\Xi}{\sqrt{k}}
\end{equation*}
with
\begin{equation*}
\Xi:=\Ex\big[\sqrt{S_1}e^{\ell S_1}\big].
\end{equation*}
Thus, $\Ex[S_1e^{(\ell-k)S_1+k X_1}]=+\infty$ or
$\Ex[S_1e^{(\ell-k)S_1+k X_1}]<+\infty$ depending on
whether $\Xi=+\infty$ or $\Xi<+\infty$.

Let us prove that $\varphi_-'(\lambda)\le 1$ when $\Xi=+\infty$. This
means that $\varphi_-'(\lambda)=1=\varphi_+'(\lambda)$ if
$\Xi=+\infty$, so that the function $\varphi$ is differentiable at
$\lambda$, and hence throughout $\Rl$.  For every $k\in(0,\lambda)$ we
have by concavity
$\varphi(\lambda)\le\varphi(k)+\varphi'(k)(\lambda-k)$, so that
\begin{equation}
  \frac{\varphi(k)-\varphi(\lambda)}{k-\lambda}\le\varphi'(k)=-\frac{\Ex[X_1e^{\varphi(k) S_1+kX_1}]}{\Ex[S_1e^{\varphi(k) S_1+kX_1}]}.
\label{upper_left_der}
\end{equation}
Pick $\epsilon\in(0,1)$ and $s>0$. The bound (\ref{upper_left_der})
gives for all $k\in(0,\lambda)$
\begin{align}
  \nonumber
  \frac{\varphi(k)-\varphi(\lambda)}{k-\lambda}&\le-\frac{\Ex[X_1e^{\varphi(k) S_1+kX_1}]}{\Ex[S_1e^{\varphi(k) S_1+kX_1}]}\\
  \nonumber
  &\le-\frac{\Ex[X_1e^{\varphi(k) S_1+kX_1}\mathds{1}_{\{X_1>(1-\epsilon)S_1\}}]}{\Ex[S_1e^{\varphi(k) S_1+kX_1}]}\\
  \nonumber
  &\le-(1-\epsilon)\frac{\Ex[S_1e^{\varphi(k) S_1+kX_1}\mathds{1}_{\{X_1>(1-\epsilon)S_1\}}]}{\Ex[S_1e^{\varphi(k) S_1+kX_1}]}\\
  \nonumber
  &=(1-\epsilon)\bigg\{-1+\frac{\Ex[S_1e^{\varphi(k) S_1+kX_1}\mathds{1}_{\{X_1\le(1-\epsilon)S_1\}}]}{\Ex[S_1e^{\varphi(k) S_1+kX_1}]}\bigg\}\\
      \nonumber
      &\le(1-\epsilon)\bigg\{-1+\frac{\Ex[S_1e^{\varphi(k)S_1+(1-\epsilon)kS_1}]}{\Ex[S_1e^{\varphi(k) S_1+kX_1}\mathds{1}_{\{S_1\le s\}}]}\bigg\}.
\end{align}
Recalling that $\varphi(k)\le\ell-k$ for $k\ge 0$ we see that
$S_1e^{\varphi(k)S_1+(1-\epsilon)kS_1}\le S_1e^{(\ell-\epsilon
  k)S_1}\le S_1e^{(\ell-\epsilon\lambda/2)S_1}$ for
$k\in(\lambda/2,\lambda)$ with
$\Ex[S_1e^{(\ell-\epsilon\lambda/2)S_1}]<+\infty$. Then, by sending
$k$ to $\lambda$, the dominated convergence theorem shows that
\begin{equation*}
\varphi_-'(\lambda)\le(1-\epsilon)\bigg\{-1+\frac{\Ex[S_1e^{(\ell-\epsilon \lambda)S_1}]}{\Ex[S_1e^{(\ell-\lambda) S_1+\lambda X_1}\mathds{1}_{\{S_1\le s\}}]}\bigg\}.
\end{equation*}
The monotone convergence theorem yields
$\lim_{s\uparrow+\infty}\Ex[S_1e^{(\ell-\lambda) S_1+\lambda
    X_1}\mathds{1}_{\{S_1\le s\}}]=+\infty$ if $\Xi=+\infty$.  Thus,
by sending $s$ to infinity we realize that
\begin{equation*}
\varphi_-'(\lambda)\le-1+\epsilon,
\end{equation*}
and the arbitrariness of $\epsilon$ implies $\varphi_-'(\lambda)\le-1$.

To conclude the proof of part $(ii)$ of the lemma it remains to verify
that
\begin{equation*}
\varphi_-'(\lambda)=\lim_{k\uparrow\lambda}\varphi'(k)=-\frac{\Ex[X_1e^{\varphi(\lambda) S_1+\lambda X_1}]}{\Ex[S_1e^{\varphi(\lambda) S_1+\lambda X_1}]}
\end{equation*}
when $\Xi<+\infty$.  The equality
$\varphi_-'(\lambda)=\lim_{k\uparrow\lambda}\varphi'(k)$ is a general
property of concave functions (see \cite{Rockbook}, Theorem 24.1).
The limit is an application of the dominated converge theorem to
(\ref{varphi_derivative}). If fact, since $X_1\le S_1$, for all
$k\in(\lambda/2,\lambda)$ we have
\begin{align}
  \nonumber
  S_1e^{\varphi(k) S_1+k X_1}\le S_1e^{(\ell-k) S_1+k X_1}&=S_1e^{\ell S_1+k(X_1-S_1)}\\
\nonumber
  &\le S_1e^{\ell S_1+(\lambda/2)(X_1-S_1)}=S_1e^{(\ell-\lambda/2) S_1+(\lambda/2)X_1}
\end{align}
with $\Ex[S_1e^{(\ell-\lambda/2) S_1+(\lambda/2)X_1}]<+\infty$ when
$\Xi<+\infty$.

\section{Proof of Lemma \ref{lem:varphi_area}}
\label{proof:varphi_area}

Part $(iii)$ of the lemma is trivial.  In order to address parts $(i)$
and $(ii)$ consider the function $\Phi$ that maps the pair
$(\zeta,k)\in\Rl^2$ in
\begin{equation*}
\Phi(\zeta,k):=\Ex\big[e^{\zeta S_1+kX_1}\big]=\Ex\big[e^{\zeta S_1+\frac{1}{6}k^2S_1^3}\big].
\end{equation*}
We have $\varphi(k)=\sup\{\zeta\in\Rl:\Phi(\zeta,k)\le 1\}$ with the
manifest symmetry $\varphi(-k)=\varphi(k)$. If $r=+\infty$, then
$\Phi$ is finite and analytic throughout $\Rl^2$. If $0<r<+\infty$,
then $\Phi$ is finite and analytic in the region of pairs
$(\zeta,k)\in\Rl^2$ with $|k|<\sqrt{6r}$, whereas
$\Phi(\zeta,k)=+\infty$ for $|k|>\sqrt{6r}$.

\subsection{Proof of part $\boldsymbol{(i)}$}

When $r=+\infty$, $\varphi(k)$ is for each $k\in\Rl$ the unique real
number $\zeta$ that solves the equation $\Phi(\zeta,k)=1$. It turns
out that $\varphi$ is analytic throughout $\Rl$ by the analytic
implicit function theorem. As $1=\Phi(\varphi(k),k)=\Ex[e^{\varphi(k)
    S_1+\frac{1}{6}k^2S_1^3}]$, we find for all $k$
\begin{equation*}
  \varphi'(k)=-\frac{\frac{\partial\Phi}{\partial k}(\varphi(k),k)}{\frac{\partial\Phi}{\partial\zeta}(\varphi(k),k)}
  =-\frac{k}{3}\frac{\Ex[S_1^3e^{\varphi(k) S_1+\frac{1}{6}k^2S_1^3}]}{\Ex[S_1e^{\varphi(k) S_1+\frac{1}{6}k^2S_1^3}]}.
\end{equation*}
Let us show by contradiction that
$\lim_{k\uparrow+\infty}\varphi'(k)=-\infty$. The limit
$\lim_{k\uparrow-\infty}\varphi'(k)=+\infty$ will follow by the
symmetry of $\varphi$. Assume that
$\eta:=\lim_{k\uparrow+\infty}\varphi'(k)>-\infty$. By concavity we
have $0=\varphi'(0)\ge\varphi'(k)\ge \eta$ and
$0=\varphi(0)\le\varphi(k)-\varphi'(k)k$, so that $\varphi(k)\ge \eta
k$ for $k>0$. We find $\varphi(k) s+(1/6)k^2s^3\ge \eta k
s+(1/6)k^2s^3\ge\sqrt{12|\eta|^3k}$ for $k>0$ and
$s>\sqrt{12|\eta|/k}$. Thus, for all $k>0$
\begin{equation*}
  1=\Ex\big[e^{\varphi(k)S_1+\frac{1}{6}k^2S_1^3}\big]\ge e^{\sqrt{12|\eta|^3k}}\,\prob\big[S_1>\sqrt{12|\eta|/k}\big],
\end{equation*}
which is a contradiction since
$\lim_{k\uparrow+\infty}\prob[S_1>\sqrt{12|\eta|/k}]=1$.

\subsection{Proof of part $\boldsymbol{(ii)}$}

If $0<r<+\infty$, then $\varphi(k)$ for $|k|<\sqrt{6r}$ is the unique
real number $\zeta$ that satisfies $\Phi(\zeta,k)=1$, whereas
$\varphi(k)=-\infty$ for $|k|>\sqrt{6r}$. The function $\varphi$ turns
out to be analytic on the open interval $(-\sqrt{6r},\sqrt{6r})$ by
the analytic implicit function theorem.  For all
$k\in(-\sqrt{6r},\sqrt{6r})$ we have $\Ex[e^{\varphi(k)
    S_1+kX_1}]=\Ex[e^{\varphi(k)
    S_1+\frac{1}{6}k^2S_1^3}]=\Phi(\varphi(k),k)=1$ and
\begin{equation*}
  \varphi'(k)=-\frac{\Ex[X_1e^{\varphi(k) S_1+kX_1}]}{\Ex[S_1e^{\varphi(k) S_1+kX_1}]}=
  -\frac{k}{3}\frac{\Ex[S_1^3e^{\varphi(k) S_1+\frac{1}{6}k^2S_1^3}]}{\Ex[S_1e^{\varphi(k) S_1+\frac{1}{6}k^2S_1^3}]}.
\end{equation*}
We study the derivatives of $\varphi$ when the boundary points
$-\sqrt{6r}$ and $\sqrt{6r}$ are approached. Thanks to the symmetry of
$\varphi$, it is sufficient to investigate the limit
$\lim_{k\uparrow\sqrt{6r}}\varphi'(k)$, which exists by concavity.

Concavity of $\varphi$ gives $k\varphi'(k)\le\varphi(k)$ for all
$k\in(0,\sqrt{6r})$ as $\varphi(0)=0$. Concavity and upper
semicontinuity of $\varphi$ imply
$\lim_{k\uparrow\sqrt{6r}}\varphi(k)=\varphi(\sqrt{6r})=:\xi$ (see
\cite{Rockbook}, Corollary 7.5.1).  Thus,
$\lim_{k\uparrow\sqrt{6r}}\varphi'(k)=-\infty$ when $\xi=-\infty$.

Let us show that $\lim_{k\uparrow\sqrt{6r}}\varphi'(k)=-\infty$ even
if $\xi>-\infty$ and $\Lambda<1$. To begin with, we claim that
$\lim_{k\uparrow\sqrt{6r}}\Ex[S_1e^{\varphi(k)S_1+\frac{1}{6}k^2S_1^3}]=+\infty$
in this case. If fact, under the assumption $\xi>-\infty$ the
inequality $e^x\ge 1+x$ valid for all $x$ gives for
$k\in(0,\sqrt{6r})$
\begin{align}
  \nonumber
  \Lambda:=\Ex\big[e^{\xi S_1+rS_1^3}\big]&\ge\Ex\Big[e^{\varphi(k)S_1+\frac{1}{6}k^2S_1^3}e^{\{\xi-\varphi(k)\}S_1}\Big]\\
  \nonumber
  &\ge\Ex\big[e^{\varphi(k)S_1+\frac{1}{6}k^2S_1^3}\big]+\big\{\xi-\varphi(k)\big\}\Ex\big[S_1e^{\varphi(k)S_1+\frac{1}{6}k^2S_1^3}\big]\\
  \nonumber
  &=1+\big\{\xi-\varphi(k)\big\}\Ex\big[S_1e^{\varphi(k)S_1+\frac{1}{6}k^2S_1^3}\big].
\end{align}
If  $\Lambda<1$, then $\varphi(k)-\xi>0$ and
\begin{equation*}
\Ex\big[S_1e^{\varphi(k)S_1+\frac{1}{6}k^2S_1^3}\big]\ge \frac{1-\Lambda}{\varphi(k)-\xi}.
\end{equation*}
By sending $k$ to $\sqrt{6r}$ from below we get
$\lim_{k\uparrow\sqrt{6r}}\Ex[S_1e^{\varphi(k)S_1+\frac{1}{6}k^2S_1^3}]=+\infty$
since we have seen that $\lim_{k\uparrow\sqrt{6r}}\varphi(k)=\xi$.  At
this point, we pick $s>0$ and use the fact that $1\le s^{-2}S_1^2$
when $S_1>s$ to state for $k\in(0,\sqrt{6r})$ the bound
\begin{align}
  \nonumber
  \varphi'(k)&=-\frac{k}{3}\frac{\Ex[S_1^3e^{\varphi(k)S_1+\frac{1}{6}k^2S_1^3}]}{\Ex[S_1e^{\varphi(k)S_1+\frac{1}{6}k^2S_1^3}]}\\
  \nonumber
  &=-\frac{k}{3}\frac{\Ex[S_1^3e^{\varphi(k)S_1+\frac{1}{6}k^2S_1^3}]}{\Ex[S_1e^{\varphi(k)S_1+\frac{1}{6}k^2S_1^3}\mathds{1}_{\{S_1\le s\}}]+
    \Ex[S_1e^{\varphi(k)S_1+\frac{1}{6}k^2S_1^3}\mathds{1}_{\{S_1>s\}}]}\\
  &\le-\frac{k}{3}\frac{\Ex[S_1^3e^{\varphi(k)S_1+\frac{1}{6}k^2S_1^3}]}{\Ex[S_1e^{\varphi(k)S_1+\frac{1}{6}k^2S_1^3}\mathds{1}_{\{S_1\le s\}}]+
    s^{-2}\Ex[S_1^3e^{\varphi(k)S_1+\frac{1}{6}k^2S_1^3}]}.
  \label{bound_derivative_area}
\end{align}
By sending $k$ to $\sqrt{6r}$ we obtain
$\lim_{k\uparrow\sqrt{6r}}\varphi'(k)\le-\sqrt{2r/3}\,s^2$ as a
consequence of the fact that
$\lim_{k\uparrow\sqrt{6r}}\Ex[S_1^3e^{\varphi(k)S_1+\frac{1}{6}k^2S_1^3}]=+\infty$
because
$\lim_{k\uparrow\sqrt{6r}}\Ex[S_1e^{\varphi(k)S_1+\frac{1}{6}k^2S_1^3}]=+\infty$.
The arbitrariness of $s$ implies
$\lim_{k\uparrow\sqrt{6r}}\varphi'(k)=-\infty$.

The case $\xi>-\infty$, $\Lambda=1$, and $\Xi:=\Ex[S_1^3e^{\xi
    S_1+rS_1^3}]=+\infty$ is immediate. Since the function that maps
$x>0$ in $x/(a+x)$ is non-decreasing for any $a\ge 0$, starting from
(\ref{bound_derivative_area}) we find for every $s>0$ and $\sigma>0$
\begin{equation*}
  \varphi'(k)\le-\frac{k}{3}\frac{\Ex[S_1^3e^{\varphi(k)S_1+\frac{1}{6}k^2S_1^3}\mathds{1}_{\{S_1\le \sigma\}}]}{\Ex[S_1e^{\varphi(k)S_1+\frac{1}{6}k^2S_1^3}\mathds{1}_{\{S_1\le s\}}]+
    s^{-2}\Ex[S_1^3e^{\varphi(k)S_1+\frac{1}{6}k^2S_1^3}\mathds{1}_{\{S_1\le \sigma\}}]}.
\end{equation*}
Thanks to the constraints on $S_1$, the dominated convergence theorem
applies and gives
\begin{equation}
\lim_{k\uparrow\sqrt{6r}}\varphi'(k)\le-\sqrt{\frac{2r}{3}}\frac{\Ex[S_1^3e^{\xi S_1+rS_1^3}\mathds{1}_{\{S_1\le \sigma\}}]}{\Ex[S_1e^{\xi S_1+rS_1^3}\mathds{1}_{\{S_1\le s\}}]+
  s^{-2}\Ex[S_1^3e^{\xi S_1+rS_1^3}\mathds{1}_{\{S_1\le \sigma\}}]}.
\label{bound_derivative_area_1}
\end{equation}
As $\lim_{\sigma\uparrow+\infty}\Ex[S_1^3e^{\xi
    S_1+rS_1^3}\mathds{1}_{\{S_1\le \sigma\}}]=\Xi=+\infty$, by
sending $\sigma$ to infinity we realize that
$\lim_{k\uparrow\sqrt{6r}}\varphi'(k)\le-\sqrt{2r/3}\,s^2$. As before,
the arbitrariness of $s$ implies
$\lim_{k\uparrow\sqrt{6r}}\varphi'(k)=-\infty$.

Finally, let us discuss the case $\xi>-\infty$, $\Lambda=1$, and
$\Xi<+\infty$. By sending first $\sigma$ to infinity and then $s$ to
infinity in (\ref{bound_derivative_area_1}) we find under the
hypothesis $\Xi<+\infty$
\begin{equation*}
\lim_{k\uparrow\sqrt{6r}}\varphi'(k)\le-\sqrt{\frac{2r}{3}}\frac{\Ex[S_1^3e^{\xi S_1+rS_1^3}]}{\Ex[S_1e^{\xi S_1+rS_1^3}]}.
\end{equation*}
Let us demonstrate the opposite bound. Once again, we use the
inequality $e^x\ge 1+x$ valid for all $x\in\Rl$ to write down for
$k\in(0,\sqrt{6r})$ the bound
\begin{align}
  \nonumber
  1=\Ex\big[e^{\varphi(k)S_1+kX_1}\big]&=\Ex\big[e^{\xi S_1+\sqrt{6r}X_1}e^{\{\varphi(k)-\xi\}S_1+\{k-\sqrt{6r}\}X_1}\big]\\
\nonumber
&\ge \Lambda+\big\{\varphi(k)-\xi\big\}\,\Ex\big[S_1e^{\xi S_1+\sqrt{6r}X_1}\big]\\
\nonumber
&+\big\{k-\sqrt{6r}\big\}\Ex\big[X_1e^{\xi S_1+\sqrt{6r}X_1}\big]\\
\nonumber
&=1+\big\{\varphi(k)-\xi\big\}\,\Ex\big[S_1e^{\xi S_1+\sqrt{6r}X_1}\big]\\
&+\big\{k-\sqrt{6r}\big\}\Ex\big[X_1e^{\xi S_1+\sqrt{6r}X_1}\big].
\label{aaa1}
\end{align}
Concavity gives $\xi:=\varphi(\sqrt{6r})\le
\varphi(k)+\varphi'(k)(\sqrt{6r}-k)$ for $k\in(0,\sqrt{6r})$. By
combining this bound with (\ref{aaa1}) we realize that for
$k\in(0,\sqrt{6r})$
\begin{equation*}
\varphi'(k)\ge -\frac{\Ex[X_1e^{\xi S_1+\sqrt{6r}X_1}]}{\Ex[S_1e^{\xi S_1+\sqrt{6r}X_1}]},
\end{equation*}
so that
\begin{equation*}
  \lim_{k\uparrow\sqrt{6r}}\varphi'(k)\ge -\frac{\Ex[X_1e^{\xi S_1+\sqrt{6r}X_1}]}{\Ex[S_1e^{\xi S_1+\sqrt{6r}X_1}]}=
  -\sqrt{\frac{2r}{3}}\frac{\Ex[S_1^3e^{\xi S_1+rS_1^3}]}{\Ex[S_1e^{\xi S_1+rS_1^3}]}.
\end{equation*}

\section{Proof of Lemma \ref{lem:auxxx}}
\label{proof:auxxx}

Fubini's theorem gives for each $k\ge
0$
\begin{align}
  \nonumber
  \Ex\Big[e^{k\int_0^1|B_{1,\tau}|\,d\tau}\Big]&=\Ex\bigg[1+\int_0^{+\infty}k e^{k x}\mathds{1}_{\{\int_0^1|B_{1,\tau}|\,d\tau>x\}}dx\bigg]\\
  &=1+\int_0^{+\infty}k e^{k x}\,\prob\bigg[\int_0^1|B_{1,\tau}|\,d\tau>x\bigg]\,dx.
  \label{auxxx1}
\end{align}
On the other hand, it is known \cite{Janson} that
$\lim_{x\uparrow+\infty}xe^{\frac{3x^2}{2}}\prob\big[\int_0^1|B_{1,\tau}|\,d\tau>x\big]=\sqrt{2/3\pi}$,
so that
\begin{equation}
\prob\bigg[\int_0^1|B_{1,\tau}|\,d\tau>x\bigg]\le C\,\frac{e^{-\frac{3x^2}{2}}}{1+x}
\label{auxxx2}
\end{equation}
for all $x>0$ with some constant $C>0$. By combining (\ref{auxxx1})
with (\ref{auxxx2}) one can prove through simple manipulations that
there exists a constant $L>0$ such that for all $k\ge 0$
\begin{equation*}
e^{\frac{1}{6}k^2}=\Ex\Big[e^{k\int_0^1B_{1,\tau}\,d\tau}\Big]\le\Ex\Big[e^{k\int_0^1|B_{1,\tau}|\,d\tau}\Big]\le L\,e^{\frac{1}{6}k^2}.
\end{equation*}
The equality on the left follows from the fact that
$\sqrt{3}\int_0^1B_{1,\tau}\,d\tau$ is distributed as a standard
Gaussian variable.

\section{Proof of Lemma \ref{lem:varphi_Aarea}}
\label{proof:varphi_Aarea}

Let $P:=\prob[S_1\in\cdot\,]$ be the probability measure induced by
$S_1$ on $[0,+\infty)$ and consider the function $\Phi$ that maps the
  pair $(\zeta,k)\in\Rl^2$ in
  \begin{align}
    \nonumber
    \Phi(\zeta,k):=\Ex\big[e^{\zeta S_1+kX_1}\big]&=\int_{[0,+\infty)}e^{\zeta s}\,\Ex\Big[e^{k\int_0^s|B_{1,\tau}|\,d\tau}\Big]P(ds)\\
      &=\int_{[0,+\infty)}e^{\zeta s}\,\Ex\Big[e^{ks^{3/2}\int_0^1|B_{1,\tau}|\,d\tau}\Big]P(ds).
        \label{Phi_int_esplicita}
\end{align}
The mapping $\Phi$ gives $\varphi$ through the formula
$\varphi(k)=\sup\{\zeta\in\Rl:\Phi(\zeta,k)\le 1\}$. To begin with, we
identify the effective domain of $\Phi$.

Formula (\ref{Airy_exp}) entails that there exists a constant $C>0$
such that $\Ex[e^{-s\int_0^1|B_{1,\tau}|\,d\tau}]\le
Ce^{-\nu_1s^{2/3}}$ for all $s\ge 0$. It follows that for every
$\zeta\in\Rl$ and $k<0$
\begin{equation}
\Phi(\zeta,k)\le C\,\Ex\big[e^{(\zeta-\nu_1|k|^{2/3})S_1}\big].
\label{bound_Phi_1}
\end{equation}
At the same time, (\ref{Airy_exp}) states that there exists $s_o>0$
such that $\Ex[e^{-s\int_0^1|B_{1,\tau}|\,d\tau}]\ge
e^{-\nu_1s^{2/3}}$ for all $s>s_o$. Thus, for every $\zeta\in\Rl$ and
$k<0$ we also have
\begin{equation}
\Phi(\zeta,k)\ge\Ex\Big[e^{(\zeta-\nu_1|k|^{2/3})S_1}\mathds{1}_{\{|k|S_1^{3/2}>s_o\}}\Big].
\label{bound_Phi_2}
\end{equation}
Bound (\ref{bound_Phi_1}) shows that $\Phi$ is finite and analytic in
the region of pairs $(\zeta,k)$ such that $\zeta<\ell+\nu_1|k|^{2/3}$
and $k<0$. Bound (\ref{bound_Phi_2}) implies that
$\Phi(\zeta,k)=+\infty$ if $\ell<+\infty$, $k<0$, and
$\zeta>\nu_1|k|^{2/3}+\ell$.

To complete the picture, we observe that Lemma \ref{lem:auxxx} gives
for all $\zeta\in\Rl$ and $k\ge 0$
\begin{equation}
\Ex\big[e^{\zeta S_1+\frac{1}{6}k^2S_1^3}\big]\le \Phi(\zeta,k)\le L\,\Ex\big[e^{\zeta S_1+\frac{1}{6}k^2S_1^3}\big]
\label{bound_Phi_3}
\end{equation}
with some positive constant $L$.  Bearing in mind that $\ell=+\infty$
when $r>0$, these bounds tell us that $\Phi$ is finite and analytic in
the region of pairs $(\zeta,k)\in\Rl^2$ such that $k<\sqrt{6r}$,
whereas $\Phi(\zeta,k)=+\infty$ if $r<+\infty$ and
$k>\sqrt{6r}$.

In the sequel we shall need the estimate
\begin{equation*}
  \Ex\big[X_1e^{\zeta S_1+kX_1}\big]\le L\,\Ex\Big[S_1^{3/2}\big(1+kS_1^{3/2}\big)e^{\zeta S_1+\frac{1}{6}k^2S_1^3}\Big],
\end{equation*}
which holds for all $\zeta\in\Rl$ and $k\ge 0$. To verify this bound,
we appeal to Lemma \ref{lem:auxxx} to get for any $s\ge 0$ and
$\lambda>0$
\begin{equation*}
  \Ex\bigg[\int_0^1|B_{1,\tau}|\,d\tau \,e^{s\int_0^1|B_{1,\tau}|\,d\tau}\bigg]\le \Ex\bigg[\frac{e^{(s+\lambda)\int_0^1|B_{1,\tau}|\,d\tau}}{\lambda}\bigg]\le
  \frac{Le^{\frac{1}{6}(s+\lambda)^2}}{\lambda}.
\end{equation*}
The choice $\lambda=2/(1+s)$ yields
\begin{equation*}
  \Ex\bigg[\int_0^1|B_{1,\tau}|\,d\tau \,e^{s\int_0^1|B_{1,\tau}|\,d\tau}\bigg]\le L(1+s)e^{\frac{1}{6}s^2}.
\end{equation*}
This way, for every $\zeta\in\Rl$ and $k\ge 0$ we find
\begin{align}
  \nonumber
  \Ex\big[X_1e^{\zeta S_1+kX_1}\big]&=\int_{[0,+\infty)}s^{3/2}e^{\zeta s}\,\Ex\bigg[\int_0^1|B_{1,\tau}|\,d\tau \,e^{ks^{3/2}\int_0^1|B_{1,\tau}|\,d\tau}\bigg]P(ds)\\
    \nonumber
    &\le L\int_{[0,+\infty)}s^{3/2}\big(1+ks^{3/2}\big)e^{\zeta s+\frac{1}{6}k^2s^3}P(ds)\\
      \nonumber
      &=L\,\Ex\Big[S_1^{3/2}\big(1+kS_1^{3/2}\big)e^{\zeta S_1+\frac{1}{6}k^2S_1^3}\Big].
\end{align}

In order to prove part $(i)$, $(ii)$, and $(iii)$ of the lemma it is
convenient to study the properties of $\varphi$ on the negative
semiaxis first.

\subsection{The function $\boldsymbol{\varphi}$ on $\boldsymbol(-\infty,0)$}
\label{sec:aux}

Bound (\ref{bound_Phi_2}) and the monotone convergence theorem show
that $\lim_{\zeta\uparrow\ell+\nu_1|k|^{2/3}}\Phi(\zeta,k)=+\infty$
for $k<0$ if either $\ell=+\infty$ or $\ell<+\infty$ and $\Ex[e^{\ell
    S_1}]=+\infty$.  Assume that $\Phi(\ell+\nu_1|k|^{2/3},k)>1$ for
all $k<0$ when $\ell<+\infty$ and $\Ex[e^{\ell S_1}]<+\infty$, as in
part $(iii)$ of the lemma. Then, irrespective of the values of $\ell$
and $\Ex[e^{\ell S_1}]$ we have
$\lim_{\zeta\uparrow\ell+\nu_1|k|^{2/3}}\Phi(\zeta,k)>1$ for $k<0$. As
a consequence, for each $k<0$, $\varphi(k)$ is the unique real number
$\zeta$ that solves the equation $\Phi(\zeta,k)=1$. It turns out that
$\varphi$ is analytic on $(-\infty,0)$ by the analytic implicit
function theorem. The identity $\Ex[e^{\varphi(k)
    S_1+kX_1}]=\Phi(\varphi(k),k)=1$ yields for any $k<0$
\begin{equation*}
  \varphi'(k)=-\frac{\frac{\partial\Phi}{\partial k}(\varphi(k),k)}{\frac{\partial\Phi}{\partial\zeta}(\varphi(k),k)}
  =-\frac{\Ex[X_1e^{\varphi(k) S_1+kX_1}]}{\Ex[S_1e^{\varphi(k) S_1+kX_1}]}.
\end{equation*}
We point out that $\varphi(k)<\ell+\nu_1|k|^{2/3}$ for all $k<0$.

Since $X_1$ is non-negative we have
$\lim_{k\downarrow-\infty}\varphi'(k)\le 0$, where the limit exists by
concavity. Let us prove that
$\lim_{k\downarrow-\infty}\varphi'(k)=0$. The concavity of $\varphi$
yields $0=\varphi(0)\le\varphi(k)-\varphi'(k)k$ for every $k<0$. Thus,
we find $\lim_{k\downarrow-\infty}\varphi'(k)=0$ if we demonstrate
that $\limsup_{k\downarrow-\infty}\varphi(k)/k\ge 0$. The latter is
immediate when $\ell<+\infty$ as $\varphi(k)<\ell+\nu_1|k|^{2/3}$ for
all $k<0$. To address the case $\ell=+\infty$, let $s_o>0$ be such
that $\Ex[e^{-s\int_0^1|B_{1,\tau}|\,d\tau}]\ge e^{-\nu_1s^{2/3}}$ for
all $s>s_o$ and recall that there exists a number $\delta>0$ such that
$\prob[S_1>\delta]>0$.  Pick $k<0$ such that $|k|\delta^{3/2}>s_o$. If
$s>\delta$, then $|k|s^{3/2}>s_o$ and we find
\begin{align}
  \nonumber
  1=\Ex\big[e^{\varphi(k)S_1+kX_1}\big]&\ge \int_{(\delta,+\infty)}e^{\varphi(k)s}\,\Ex\Big[e^{-|k|s^{3/2}\int_0^1|B_{1,\tau}|\,d\tau}\Big]P(ds)\\
  \nonumber
  &\ge\int_{(\delta,+\infty)}e^{\{\varphi(k)-\nu_1|k|^{2/3}\}s}P(ds)\\
  \nonumber
  &=\Ex\big[e^{\{\varphi(k)-\nu_1|k|^{2/3}\}S_1}\big|S_1>\delta\big]\prob\big[S_1>\delta\big].
\end{align}
Jensen's inequality allows us to conclude that
\begin{align}
  \nonumber
  1&\ge\Ex\big[e^{\{\varphi(k)-\nu_1|k|^{2/3}\}S_1}\big|S_1>\delta\big]\prob\big[S_1>\delta\big]\ge e^{\{\varphi(k)-\nu_1|k|^{2/3}\}\Ex[S_1|S_1>\delta]}\,\prob\big[S_1>\delta\big].
\end{align}
Notice that $\Ex[S_1|S_1>\delta]<+\infty$ since $\ell=+\infty$.  Thus,
for any $k<0$ such that $|k|\delta^{3/2}>s_o$ we have
\begin{equation*}
\varphi(k)\le\nu_1|k|^{2/3}-\frac{\ln\prob[S_1>\delta]}{\Ex[S_1|S_1>\delta]},
\end{equation*}
which implies $\limsup_{k\downarrow-\infty}\varphi(k)/k\ge 0$.

\subsection{Proof of part $\boldsymbol{(i)}$}

If $r=+\infty$, then $\Phi$ is analytic throughout $\Rl^2$. In this
case, $\varphi(k)$ is for each $k\in\Rl$ the unique real number
$\zeta$ that solves the equation $\Phi(\zeta,k)=1$. It turns out that
$\varphi$ is analytic throughout $\Rl$ by the analytic implicit
function theorem. As $1=\Phi(\varphi(k),k)=\Ex[e^{\varphi(k)
    S_1+kX_1}]$, we find for all $k$
\begin{equation*}
  \varphi'(k)=-\frac{\frac{\partial\Phi}{\partial k}(\varphi(k),k)}{\frac{\partial\Phi}{\partial\zeta}(\varphi(k),k)}
  =-\frac{\Ex[X_1e^{\varphi(k) S_1+kX_1}]}{\Ex[S_1e^{\varphi(k) S_1+kX_1}]}.
\end{equation*}

We know from the previous section that
$\lim_{k\downarrow-\infty}\varphi'(k)=0$, and one can show that
$\lim_{k\uparrow+\infty}\varphi'(k)=-\infty$ by means of the same
arguments we used for the area. Indeed, (\ref{bound_Phi_3}) gives
$1=\Phi(\varphi(k),k)\ge \Ex[e^{\varphi(k) S_1+\frac{1}{6}k^2S_1^3}]$
for any $k>0$.

\subsection{Proof of part $\boldsymbol{(ii)}$}

Assume that $0<r<+\infty$. In this case, the function $\Phi$ is
analytic in the region of pairs $(\zeta,k)\in\Rl^2$ with
$k<\sqrt{6r}$, whereas $\Phi(\zeta,k)=+\infty$ if $k>\sqrt{6r}$. As a
consequence, $\varphi(k)$ for $k<\sqrt{6r}$ is the unique real number
$\zeta$ that satisfies $\Phi(\zeta,k)=1$, whereas $\varphi(k)=-\infty$
for $k>\sqrt{6r}$. The function $\varphi$ turns out to be analytic on
$(-\infty,\sqrt{6r})$ by the analytic implicit function theorem.  For
all $k\in(-\infty,\sqrt{6r})$ we have $\Ex[e^{\varphi(k)
    S_1+kX_1}]=\Phi(\varphi(k),k)=1$ and
\begin{equation*}
  \varphi'(k)=-\frac{\Ex[X_1e^{\varphi(k) S_1+kX_1}]}{\Ex[S_1e^{\varphi(k) S_1+kX_1}]}.
\end{equation*}

In Section \ref{sec:aux} we proved that
$\lim_{k\downarrow-\infty}\varphi'(k)=0$. We now investigate the
derivatives of $\varphi$ when the boundary point $\sqrt{6r}$ is
approached. Basically, the arguments are the ones we used for the
area, and the limit $\lim_{k\uparrow\sqrt{6r}}\varphi'(k)=-\infty$ is
proved exactly in the same way when $\xi=-\infty$. Let us show that
$\lim_{k\uparrow\sqrt{6r}}\varphi'(k)=-\infty$ even if $\xi>-\infty$
and $\Lambda<1$ or $\xi>-\infty$ and $\Lambda=1$ and $\Xi=+\infty$.
Recalling that $\sqrt{3}\int_0^1B_{1,\tau}\,d\tau$ is distributed as a
Gaussian variable with mean $0$ and variance $1$, we can write down for
$k\in(0,\sqrt{6r})$ the bound
\begin{align}
  \nonumber
  \frac{k}{3}\Ex\big[S_1^3e^{\varphi(k) S_1+\frac{1}{6}k^2S_1^3}\big]&=\int_{[0,+\infty)}\frac{ks^3}{3}e^{\varphi(k)s+\frac{1}{6}k^2s^3}P(ds)\\
    \nonumber
  &=\int_{[0,+\infty)}\Ex\bigg[s^{3/2}\int_0^1B_{1,\tau}\,d\tau\,e^{\varphi(k)s+ks^{3/2}\int_0^1 B_{1,\tau}\,d\tau}\bigg]P(ds)\\
\nonumber
&\le\int_{[0,+\infty)}\Ex\bigg[\int_0^s|B_{1,\tau}|\,d\tau\,e^{\varphi(k)s+k\int_0^s |B_{1,\tau}|\,d\tau}\bigg]P(ds)\\
  \nonumber
  &=\Ex\big[X_1e^{\varphi(k) S_1+kX_1}\big].
\end{align}
At the same time, Lemma \ref{lem:auxxx} yields
\begin{align}
  \nonumber
  \Ex\big[S_1e^{\varphi(k) S_1+kX_1}\big]&=\int_{[0,+\infty)}se^{\varphi(k)s}\,\Ex\Big[e^{ks^{3/2}\int_0^1 B_{1,\tau}\,d\tau}\Big]P(ds)\\
    &\le L\,\Ex\big[S_1e^{\varphi(k) S_1+\frac{1}{6}k^2S_1^3}\big].
\label{XYXYXY}
\end{align}
It follows that for every $k\in(0,\sqrt{6r})$
\begin{equation*}
  \varphi'(k)=-\frac{\Ex[X_1e^{\varphi(k) S_1+kX_1}]}{\Ex[S_1e^{\varphi(k) S_1+kX_1}]}\le
  -\frac{k}{3L}\frac{\Ex[S_1^3e^{\varphi(k) S_1+\frac{1}{6}k^2S_1^3}]}{\Ex[S_1e^{\varphi(k) S_1+\frac{1}{6}k^2S_1^3}]}.
\end{equation*}
Starting from this bound, one can demonstrate that
$\lim_{k\uparrow\sqrt{6r}}\varphi'(k)=-\infty$ if $\xi>-\infty$ and
$\Lambda<1$ or $\xi>-\infty$ and $\Lambda=1$ and $\Xi=+\infty$ by
resorting to the strategies we devised for the area. It is only needed
to check that
$\lim_{k\uparrow\sqrt{6r}}\Ex[S_1e^{\varphi(k)S_1+\frac{1}{6}k^2S_1^3}]=+\infty$
when $\xi>-\infty$ and $\Lambda<1$. Actually, the inequality $e^x\ge
1+x$ valid for all $x$ gives for $k\in(0,\sqrt{6r})$
\begin{align}
  \nonumber
  \Lambda:=\Ex\big[e^{\xi S_1+\sqrt{6r}X_1}\big]&\ge\Ex\Big[e^{\varphi(k)S_1+kX_1}e^{\{\xi-\varphi(k)\}S_1}\Big]\\
  \nonumber
  &\ge\Ex\big[e^{\varphi(k)S_1+kX_1}\big]+\big\{\xi-\varphi(k)\big\}\Ex\big[S_1e^{\varphi(k)S_1+kX_1}\big]\\
  \nonumber
  &=1+\big\{\xi-\varphi(k)\big\}\Ex\big[S_1e^{\varphi(k)S_1+kX_1}\big].
\end{align}
If  $\Lambda<1$, then $\varphi(k)-\xi>0$ and
\begin{equation*}
\frac{1-\Lambda}{\varphi(k)-\xi}\le\Ex\big[S_1e^{\varphi(k)S_1+kX_1}\big]\le L\,\Ex\big[S_1e^{\varphi(k)S_1+\frac{1}{6}k^2S_1^3}\big],
\end{equation*}
the last bound being (\ref{XYXYXY}). Since
$\lim_{k\uparrow\sqrt{6r}}\varphi(k)=\xi$ as $-\varphi$ is convex and
lower semicontinuous (see \cite{Rockbook}, Corollary 7.5.1), from here
we get
$\lim_{k\uparrow\sqrt{6r}}\Ex[S_1e^{\varphi(k)S_1+\frac{1}{6}k^2S_1^3}]=+\infty$.

Finally, let us show that if $\xi>-\infty$ and $\Lambda=1$ and
$\Xi<+\infty$, then
\begin{equation*}
\lim_{k\uparrow\sqrt{6r}}\varphi'(k)=-\frac{\Ex[X_1e^{\xi S_1+\sqrt{6r}X_1}]}{\Ex[S_1e^{\xi S_1+\sqrt{6r}X_1}]}.
\end{equation*}
We stress that $\Ex[X_1e^{\xi S_1+\sqrt{6r}X_1}]<+\infty$ when
$\Xi:=\Ex[S_1^3e^{\xi S_1+rS_1^3}]<+\infty$. In fact, we have seen at
the beginning that $\Ex[X_1e^{\zeta S_1+kX_1}]\le
L\,\Ex[S_1^{3/2}(1+kS_1^{3/2})e^{\zeta S_1+\frac{1}{6}k^2S_1^3}]$ for
all $\zeta\in\Rl$ and $k\ge 0$ with some positive constant $L$.  As
for the area, for all $k\in(0,\sqrt{6r})$ we have
\begin{equation*}
\varphi'(k)\ge -\frac{\Ex[X_1e^{\xi S_1+\sqrt{6r}X_1}]}{\Ex[S_1e^{\xi S_1+\sqrt{6r}X_1}]}.
\end{equation*}
The problem is to prove that
\begin{equation}
\lim_{k\uparrow\sqrt{6r}}\varphi'(k)\le-\frac{\Ex[X_1e^{\xi S_1+\sqrt{6r}X_1}]}{\Ex[S_1e^{\xi S_1+\sqrt{6r}X_1}]}.
\label{ultimo_sforzo}
\end{equation}
The limit exists by concavity. Pick $k\in(0,\sqrt{6r})$ and
$\lambda>0$ and observe that
\begin{align}
  \nonumber
  \varphi'(k)&=-\frac{\Ex[X_1e^{\varphi(k) S_1+kX_1}]}{\Ex[S_1e^{\varphi(k) S_1+kX_1}]}\\
  &\le-\frac{\Ex[X_1e^{\varphi(k) S_1+kX_1}]}{\Ex[S_1e^{\varphi(k) S_1+kX_1}\mathds{1}_{\{X_1\le \lambda S_1\}}]+\lambda^{-1}\Ex[X_1e^{\varphi(k) S_1+kX_1}]}.
\label{ultimo_sforzo_1}
\end{align}
We shall use the fact that $\xi\le\varphi(k)\le 0$ for
$k\in(0,\sqrt{6r})$ as $X_1\ge 0$.  Since the function that maps $x>0$
in $x/(a + x)$ is non-decreasing for any $a\ge 0$ we can change
(\ref{ultimo_sforzo_1}) with
\begin{equation*}
  \varphi'(k)\le-\frac{\Ex[X_1e^{\xi S_1+kX_1}]}{\Ex[S_1e^{\varphi(k) S_1+kX_1}\mathds{1}_{\{X_1\le\lambda S_1\}}]
      +\lambda^{-1}\Ex[X_1e^{\xi S_1+kX_1}]}.
\end{equation*}
On the other hand, the conditions
$\sqrt{S_1}\int_0^1|B_{1,\tau}|\,d\tau\le\lambda$ and
$\int_0^1|B_{1,\tau}|\,d\tau>\lambda$ imply $S_1\le 1$, so that
\begin{align}
  \nonumber
  \Ex\Big[S_1e^{\varphi(k) S_1+kX_1}\mathds{1}_{\{X_1\le\lambda S_1\}}\Big]
  &=\Ex\Big[S_1e^{\varphi(k) S_1+kS_1^{3/2}\int_0^1|B_{1,\tau}|\,d\tau}\mathds{1}_{\{\sqrt{S_1}\int_0^1|B_{1,\tau}|\,d\tau\le\lambda \}}\Big]\\
  \nonumber
  &\le E_\lambda(k)+\epsilon_\lambda
\end{align}
with
\begin{equation*}
E_\lambda(k):=\Ex\Big[S_1e^{\varphi(k) S_1+kS_1^{3/2}\int_0^1|B_{1,\tau}|\,d\tau}\mathds{1}_{\{\int_0^1|B_{1,\tau}|\,d\tau\le\lambda\}}\Big]
\end{equation*}
and
\begin{equation*}
\epsilon_\lambda:=\Ex\Big[e^{\sqrt{6r}\int_0^1|B_{1,\tau}|\,d\tau}\mathds{1}_{\{\int_0^1|B_{1,\tau}|\,d\tau>\lambda\}}\Big].
\end{equation*}
We stress that $\lim_{\lambda\uparrow+\infty}\epsilon_\lambda=0$ by
the dominated converge theorem. By putting the pieces together we find
\begin{equation*}
  \varphi'(k)\le-\frac{\Ex[X_1e^{\xi S_1+kX_1}]}{E_\lambda(k)+\epsilon_\lambda
      +\lambda^{-1}\Ex[X_1e^{\xi S_1+kX_1}]}.
\end{equation*}
At this point, we notice that the dominated convergence theorem gives
$\lim_{k\uparrow\sqrt{6r}}E_\lambda(k)=E_\lambda(\sqrt{6r})$ since
$\Ex[e^{\zeta S_1^{3/2}}]<+\infty$ for all $\zeta\in\Rl$ due to the
fact that $r>0$. The monotone converge theorem yields
$\lim_{k\uparrow\sqrt{6r}}\Ex[X_1e^{\xi S_1+kX_1}]=\Ex[X_1e^{\xi
    S_1+\sqrt{6r}X_1}]$. It follows that
\begin{equation*}
  \lim_{k\uparrow\sqrt{6r}}\varphi'(k)\le-\frac{\Ex[X_1e^{\xi S_1+\sqrt{6r}X_1}]}{E_\lambda(\sqrt{6r})+\epsilon_\lambda
      +\lambda^{-1}\Ex[X_1e^{\xi S_1+\sqrt{6r}X_1}]}.
\end{equation*}
This bound demonstrates (\ref{ultimo_sforzo}) by sending $\lambda$ to
infinity as $\lim_{\lambda\uparrow+\infty}\epsilon_\lambda=0$ and
\begin{equation*}
\lim_{\lambda\uparrow+\infty}E_\lambda(\sqrt{6r})=\Ex\Big[S_1e^{\xi S_1+\sqrt{6r}S_1^{3/2}\int_0^1|B_{1,\tau}|\,d\tau}\Big]=\Ex\big[S_1e^{\xi S_1+\sqrt{6r}X_1}\big].
\end{equation*}

\subsection{Proof of part $\boldsymbol{(iii)}$}

In the light of Section \ref{sec:aux}, it remains to verify that
\begin{equation*}
    \lim_{k\uparrow 0}\varphi'(k)=\begin{cases}
    -\frac{\Ex[X_1]}{\Ex[S_1]}=-\sqrt{\frac{8}{9\pi}}\frac{\Ex[S_1^{3/2}]}{\Ex[S_1]} & \mbox{if }\Ex[S_1^{3/2}]<+\infty,\\
    -\infty & \mbox{if }\Ex[S_1^{3/2}]=+\infty.
    \end{cases}
\end{equation*}
The case $\ell>0$ is immediate. Assume that $\ell>0$ and fix a real
number $\zeta<\ell$. Since $\lim_{k\uparrow 0}\varphi(k)=0$ as
$-\varphi$ is convex and lower semicontinuous (see \cite{Rockbook},
Corollary 7.5.1), there exists $\delta>0$ such that
$\varphi(k)\le\zeta$ for $k\in(-\delta,0)$.  Then, $(X_1\vee
S_1)e^{\varphi(k)S_1+kX_1}\le (X_1\vee S_1)e^{\zeta S_1}$ for
$k\in(-\delta,0)$ with $\Ex[(X_1\vee S_1)e^{\zeta S_1}]<+\infty$ as
$\zeta<\ell$. The dominated converge theorem gives
\begin{equation*}
  \lim_{k\uparrow 0}\varphi'(k)=-\lim_{k\uparrow 0}\frac{\Ex[X_1e^{\varphi(k) S_1+kX_1}]}{\Ex[S_1e^{\varphi(k) S_1+kX_1}]}=-\frac{\Ex[X_1]}{\Ex[S_1]}.
\end{equation*}

To address the case $\ell=0$ and $\Ex[S_1^{3/2}]<+\infty$, we notice
at first that the identity $\Ex[e^{\varphi(k)S_1+kX_1}]=1$ for $k<0$,
combined with the inequality $e^x\ge 1+x$ valid for all $x\in\Rl$,
gives $\varphi(k)\Ex[S_1]+k\Ex[X_1]\le 0$ with
$\Ex[X_1]=\sqrt{8/9\pi}\,\Ex[S_1^{3/2}]<+\infty$. On the other hand,
the concavity of $\varphi$ shows that $0=\varphi(0)\le
\varphi(k)-\varphi'(k)k$ for $k<0$. Thus,
$\varphi'(k)\ge-\Ex[X_1]/\Ex[S_1]$ for all $k<0$. Let us prove that
$\lim_{k\uparrow 0}\varphi'(k)\le -\Ex[X_1]/\Ex[S_1]$. Since for every
$k<0$
\begin{equation*}
  \varphi'(k)=-\frac{\Ex[X_1e^{\varphi(k) S_1+kX_1}]}{\Ex[S_1e^{\varphi(k) S_1+kX_1}]}\le -\frac{\Ex[X_1e^{kX_1}]}{\Ex[S_1e^{\varphi(k) S_1+kX_1}]}
\end{equation*}
as $\varphi(k)\ge 0$ and since $\lim_{k\uparrow
  0}\Ex[X_1e^{kX_1}]=\Ex[X_1]$ by the monotone converge theorem, we
find $\lim_{k\uparrow 0}\varphi'(k)\le -\Ex[X_1]/\Ex[S_1]$ if
$\lim_{k\uparrow 0}\Ex[S_1e^{\varphi(k) S_1+kX_1}]=\Ex[S_1]$. The
monotone converge theorem gives $\lim_{k\uparrow
  0}\Ex[e^{-|k|s^{3/2}\int_0^1|B_{1,\tau}|\,d\tau}]=1$ for all $s\ge
0$.  At the same time, for $k<0$ we have the two bounds $\varphi(k)\le
\nu_1|k|^{2/3}$ as $\ell=0$ and
$\Ex[e^{-|k|s^{3/2}\int_0^1|B_{1,\tau}|\,d\tau}]\le
Ce^{-\nu_1|k|^{2/3}s}$ for any $s\ge 0$. The latter allow us to invoke
the dominated convergence theorem to state that
\begin{align}
  \nonumber
  \lim_{k\uparrow 0}\Ex\big[S_1e^{\varphi(k) S_1+kX_1}\big]&=\int_{[0,+\infty)}\lim_{k\uparrow 0}se^{\varphi(k)s}\,\Ex\Big[e^{-|k|s^{3/2}\int_0^1|B_{1,\tau}|\,d\tau}\Big]P(ds)\\
\nonumber
    &=\int_{[0,+\infty)}sP(ds)=\Ex[S_1].
\end{align}

Finally, we discuss the case $\Ex[S_1^{3/2}]=+\infty$. Pick
$\lambda>0$ and observe that, similarly to (\ref{ultimo_sforzo_1}),
for all $k<0$
\begin{align}
  \nonumber
  \varphi'(k)&=-\frac{\Ex[X_1e^{\varphi(k) S_1+kX_1}]}{\Ex[S_1e^{\varphi(k) S_1+kX_1}]}\\
\nonumber
&\le-\frac{\Ex[X_1e^{\varphi(k) S_1+kX_1}]}{\Ex[S_1e^{\varphi(k) S_1+kX_1}\mathds{1}_{\{X_1\le \lambda S_1\}}]+\lambda^{-1}\Ex[X_1e^{\varphi(k) S_1+kX_1}]}.
\end{align}
Since $\varphi(k)\ge 0$ for $k<0$ and since the function that maps
$x>0$ in $\frac{x}{a+x}$ is non-decreasing for any $a\ge 0$, we find
\begin{equation*}
  \varphi'(k)\le-\frac{\Ex[X_1e^{kX_1}]}{\Ex[S_1e^{\varphi(k) S_1+kX_1}\mathds{1}_{\{X_1\le \lambda S_1\}}]+\lambda^{-1}\Ex[X_1e^{kX_1}]}.
\end{equation*}
Now we make use of the known fact \cite{Janson} that there exists a
constant $D>0$ such that for all $x\ge 0$
\begin{equation*}
  \prob\bigg[\int_0^1|B_{1,\tau}|\,d\tau\le x\bigg]\le De^{-\frac{4\nu_1^3}{27 x^2}}.
\end{equation*}
Then, for all $k<0$ in a neighborhood of the origin we have
$\varphi(k)\le 2\nu_1^3/27\lambda^2$ as $\lim_{k\uparrow 0}\varphi(k)=0$ and
\begin{align}
  \nonumber
  \Ex\Big[S_1e^{\varphi(k) S_1+kX_1}\mathds{1}_{\{X_1\le \lambda S_1\}}\Big]&\le\Ex\Big[S_1e^{\varphi(k) S_1}\mathds{1}_{\{\sqrt{S_1}\int_0^1|B_{1,\tau}|\,d\tau\le \lambda \}}\Big]\\
  \nonumber
  &\le D\,\Ex\Big[S_1e^{\{\varphi(k)-\frac{4\nu_1^3}{27 \lambda^2}\}S_1}\Big]\le D\,\Ex\Big[S_1e^{-\frac{2\nu_1^3}{27 \lambda^2}S_1}\Big]<+\infty.
\end{align}
In conclusion, for every $k<0$ in a neighborhood of the origin
\begin{equation}
  \varphi'(k)\le-\frac{\Ex[X_1e^{kX_1}]}{D\,\Ex\big[S_1e^{-\frac{2\nu_1^3}{27 \lambda^2}S_1}\big]+\lambda^{-1}\Ex[X_1e^{kX_1}]}.
\label{der_phi_origin}
\end{equation}
The monotone converge theorem yields $\lim_{k\uparrow
  0}\Ex[X_1e^{kX_1}]=\Ex[X_1]=\sqrt{8/9\pi}\,\Ex[S_1^{3/2}]=+\infty$. Thus,
by sending $k$ to $0$ from below in (\ref{der_phi_origin}) we get
$\lim_{k\uparrow 0}\varphi'(k)\le-\lambda$. The arbitrariness of
$\lambda$ implies $\lim_{k\uparrow 0}\varphi'(k)=-\infty$.



\end{document}